\documentclass[5p,twocolumn,authoryear,sort&compress]{elsarticle}

\usepackage{lineno}
\modulolinenumbers[5]

\RequirePackage[dvipsnames]{xcolor}
\definecolor{ETHcolor}{RGB}{18,105,176}
\usepackage[
  colorlinks=true,
  hyperindex=true, %
  colorlinks=true,%
  pagebackref=false,%
  plainpages=false,%
  pdfpagelabels,
  ]{hyperref}

\journal{Annual Reviews in Control}

\usepackage[utf8]{inputenc}    %
\usepackage[T1]{fontenc}
\usepackage{amsmath,amsfonts,amssymb,amsthm}
\allowdisplaybreaks
\usepackage{mathtools}
\usepackage{xfrac}

\newtheorem{theorem}{Theorem}

\definecolor{vermilion}{rgb}{0.89, 0.26, 0.2}
\definecolor{tuftsblue}{rgb}{0.28, 0.57, 0.81}
\definecolor{persiangreen}{rgb}{0.0, 0.65, 0.58}

\DeclareMathOperator{\minimize}{minimize}
\DeclareMathOperator{\subjto}{subject~to}

\newcommand{\calX}{\mathcal{X}}

\newcommand{\bbR}{\mathbb{R}}
\newcommand{\calU}{\mathcal{U}}

\newcommand{\tproj}[3]{\Pi_{#1}^{#2}[{#3}]}

\newenvironment{smallbmatrix}
  {\left[\begin{smallmatrix}}
  {\end{smallmatrix}\right]}

  \DeclareMathAlphabet{\anothermymathbb}{U}{bbold}{m}{n}
  \DeclareMathAlphabet{\mymathbb}{U}{BOONDOX-ds}{m}{n}

  \newcommand{\bbI}{\mathbb{I}}

\newcommand{\rmb}{\mathrm{b}}
\newcommand{\rmg}{\mathrm{g}}
\newcommand{\C}{\mathrm{C}}
\renewcommand{\L}{\mathrm{L}}
\newcommand{\G}{\mathrm{G}}
\newcommand{\M}{\mathrm{M}}
\newcommand{\bfz}{\mathbf{z}}
\newcommand{\bfu}{\mathbf{u}}

\newcommand{\bfy}{\mathbf{y}}
\newcommand{\bfP}{\mathbf{P}}
\newcommand{\bfA}{\mathbf{A}}
\newcommand{\bfB}{\mathbf{B}}
\newcommand{\bfC}{\mathbf{C}}
\newcommand{\bfD}{\mathbf{D}}
\newcommand{\calN}{\mathcal{N}}
\newcommand{\calM}{\mathcal{M}}

\DeclareMathOperator{\cl}{cl}
\DeclareMathOperator{\E}{\mathbb{E}}

\usepackage[inline]{enumitem}
\usepackage{graphicx}

\usepackage{thmtools}
\usepackage{nameref}

\declaretheorem[name=Remark,refname={Remark,Remarks},Refname={Remark,Remarks},numberwithin=section,style=remark, qed=\hbox{\small $\blacksquare$}]{remark}
\declaretheorem[name=Example,refname={Example,Examples},Refname={Example,Examples},numberwithin=section,style=remark,qed=\hbox{\small $\blacksquare$}]{example}

\usepackage[capitalize,nameinlink,sort]{cleveref}
\crefname{part}{Part}{Parts}
\crefname{figure}{Figure}{Figures} %
\crefname{assumption}{Assumption}{Assumptions}

\crefname{equation}{}{}
\crefrangelabelformat{equation}{(#3#1#4--#5\crefstripprefix{#1}{#2}#6)}
\crefmultiformat{equation}{(#2#1#3}{, #2#1#3)}{#2#1#3}{#2#1#3}

\renewcommand{\theenumi}{(\roman{enumi})}
\renewcommand{\theenumii}{(\alph{enumii})}

\makeatletter
\renewcommand\p@enumii{\theenumi}
\renewcommand\p@enumiii{\theenumi\theenumii}
\makeatother

\crefname{enumi}{}{}
\crefname{enumii}{}{}
\crefname{enumiii}{}{}

\usepackage[mode=buildnew]{standalone}
\usepackage{graphicx}
\usepackage{soul}
\usepackage{subcaption}
\usepackage{pgfplots}
\usepackage{wrapfig}
\usepackage{tikz}
\usepackage{tikzscale}
\usepackage[nooldvoltagedirection]{circuitikz}

\pgfplotsset{compat=1.16}
\DeclareGraphicsExtensions{.eps,.pdf,.png,.jpg}

\usepgfplotslibrary{colorbrewer}
\usepgfplotslibrary{fillbetween}
\usetikzlibrary{arrows,arrows.meta}
\usetikzlibrary{calc}
\usetikzlibrary{fit}
\usetikzlibrary{shapes,shapes.arrows}
\usetikzlibrary{cd}
\usepgfplotslibrary{groupplots}

\tikzstyle{sum} = [draw, thick, inner sep=0mm, circle, minimum size=1mm] %
\tikzstyle{gainright} = [draw, thick, isosceles triangle, minimum height = 3mm, isosceles triangle apex angle=60]
\tikzstyle{gainleft} = [draw, thick, isosceles triangle, minimum height = 3mm, inner sep = 0mm, isosceles triangle apex angle=60, shape border rotate=180]
\tikzstyle{gainup} = [draw, thick, isosceles triangle, minimum height = 8mm, inner sep = 0mm, isosceles triangle apex angle=60, shape border rotate=90]
\tikzstyle{none} = [draw=none]
\tikzstyle{connector} = [->,thick]
\tikzstyle{line} = [thick]
\tikzstyle{block} = [draw, rectangle, thick,minimum height=1em,minimum width=1em]
\tikzstyle{dashedblock} = [draw, rectangle, thick,dashed,minimum height=1em,minimum width=1em]
\tikzstyle{smallsum} = [draw,circle,inner sep=0mm,minimum size=3mm]
\tikzstyle{branch} = [draw,circle,inner sep=0.5mm,fill=black]

\tikzset{%
  saturation block/.style={%
    draw, thick,
    path picture={
      \pgfpointdiff{\pgfpointanchor{path picture bounding box}{north east}}%
        {\pgfpointanchor{path picture bounding box}{south west}}
      \pgfgetlastxy\x\y
      \tikzset{x=\x*.4, y=\y*.4}
      \draw (-.9,0) -- (.9,0) (0,-.9) -- (0,.9);
      \draw (-.9,-.6) -- (-.6,-.6) -- (.6,.6) -- (.9,.6);
      \node[text width=.1cm] at (-.35,.55) {\scriptsize $P_\calU$};
    }
  }
}

\usepackage{acro}
\acsetup{format/first-long = {}}

\DeclareAcronym{pds}{short=PDS, long=projected dynamical system, long-plural=s, tag=abbrev}
\DeclareAcronym{awa}{short=AWA, long=anti-windup approximation, long-plural=s, tag=abbrev}
\DeclareAcronym{pfm}{short=PFM, long=power flow manifold, tag=abbrev}

\DeclareAcronym{lop}{short=LOP, long=linearized output projection, tag=abbrev}
\DeclareAcronym{ode}{short=ODE, long=ordinary differential equation, long-plural=s, short-plural=s, tag=abbrev}

\DeclareAcronym{dvi}{short=DVI, long=differential variational inequality, long-plural-form=differential variational inequalities, short-plural=s, tag=abbrev}

\DeclareAcronym{lti}{short=LTI, long=linear time-invariant, tag=abbrev}

\DeclareAcronym{cq}{short=CQ, long=constraint qualification, tag=abbrev}
\DeclareAcronym{licq}{short=LICQ, long=linear independence constraint qualification, tag=abbrev}
\DeclareAcronym{kkt}{short=KKT, long=Karush-Kuhn-Tucker, tag=abbrev}
\DeclareAcronym{ssosc}{short=SSOSC, long=Strong Second-Order Sufficiency Condition, tag=abbrev}

\DeclareAcronym{osc}{short=osc, long=outer semicontinuous, tag=abbrev}
\DeclareAcronym{isc}{short=isc, long=inner semicontinuous, tag=abbrev}

\DeclareAcronym{las}{short=LAS, long=locally asymptotically stable, tag=abbrev}
\DeclareAcronym{gas}{short=GAS, long=globally asymptotically stable, tag=abbrev}
\DeclareAcronym{spas}{short=SPAS, long=semiglobally practically asymptotically stable, tag=abbrev}

\DeclareAcronym{qp}{short=QP, long=quadratic program, tag=abbrev}
\DeclareAcronym{sqp}{short=SQP, long=sequential quadratic programming, tag=abbrev}

\DeclareAcronym{es}{short=ES, long=extremum seeking, tag=abbrev}
\DeclareAcronym{mpc}{short=MPC, long=model predictive control, tag=abbrev}
\DeclareAcronym{ma}{short=MA, long=modifier adaptation, tag=abbrev}
\DeclareAcronym{rti}{short=RTI, long=real-time iteration, long-plural=s, tag=abbrev}
\DeclareAcronym{lmi}{short=LMI, long=linear matrix inequality, long-plural-form=linear matrix inequalities, tag=abbrev}

\DeclareAcronym{ac}{short=AC, long=alternating current, tag=abbrev}
\DeclareAcronym{acpf}{short=ACPF, long=AC power flow, tag=abbrev}
\DeclareAcronym{acopf}{short=ACOPF, long=AC optimal power flow, tag=abbrev}
\DeclareAcronym{matpower}{short=MATPOWER , long=AC power flow, tag=abbrev}

\DeclareAcronym{qcqp}{short=QCQP, long=quadratically constrained quadratic program, tag=abbrev}

\DeclareAcronym{sdp}{short=SDP, long=semidefinite program, tag=abbrev}

\DeclareAcronym{submp}{short=\textsc{Sub-MP}, long=submodular multiway partition problem, tag=abbrev}
\DeclareAcronym{admm}{short=ADMM, long=alternating direction method of multipliers, tag=abbrev}

\DeclareAcronym{flop}{short=FLOP, long=floating point operation, tag=abbrev}
\DeclareAcronym{nnz}{short=NNZs, long=number of non-zero elements, tag=abbrev}

\definecolor{darkgreen}{rgb}{0.0, 0.5, 0.0}
\usepackage[prependcaption,colorinlistoftodos]{todonotes}

\usepackage{etoolbox}
\makeatletter
\patchcmd{\@addmarginpar}{\ifodd\c@page}{\ifodd\c@page\@tempcnta\m@ne}{}{}
\makeatother
\reversemarginpar

\newcommand{\tb}[0]{\color{blue}}

\renewcommand{\tb}[0]{}

\usepackage[normalem]{ulem}
\newcommand\rout{\bgroup\markoverwith{\textcolor{red}{/}}\ULon} %

\usepackage{xcolor,soul}
\usepackage{framed}
\colorlet{shadecolor}{yellow}
\soulregister\cite7
\soulregister\cref7
\soulregister\ref7
\soulregister\pageref7
\soulregister\eqref7

\hypersetup{
  pdftitle={Optimization Algorithms as Robust Feedback Controllers},%
  pdfauthor={Adrian Hauswirth, Zhiyu He, Saverio Bolognani, Gabriela Hug, and Florian Dörfler},%
}

\begin{document}

\hypersetup{
    linkcolor=ETHcolor,%
    citecolor=ETHcolor,%
    filecolor=ETHcolor,%
    urlcolor=ETHcolor%
}

\begin{frontmatter}

    \title{Optimization Algorithms as Robust Feedback Controllers}
    \author{Adrian Hauswirth, Zhiyu He, Saverio Bolognani\corref{mycorrespondingauthor}, Gabriela Hug, and Florian D\"orfler}
    \address{Department of Information Technology and Electrical Engineering, ETH Zurich, Switzerland}

    \begin{abstract}
        Mathematical optimization is one of the cornerstones of modern engineering research and practice. Yet, throughout all application domains, mathematical optimization is, for the most part, considered to be a numerical discipline. Optimization problems are formulated to be solved numerically with specific algorithms running on microprocessors. An emerging alternative is to view optimization algorithms as dynamical systems. Besides being insightful in itself, this perspective liberates optimization methods from specific numerical and algorithmic aspects and opens up new possibilities to endow complex real-world systems with sophisticated self-optimizing behavior.
        Towards this goal, it is necessary to understand how numerical optimization algorithms can be converted into feedback controllers to enable robust ``closed-loop optimization''.
        {\tb In this article, we focus on recent control designs under the name of ``feedback-based optimization'' which implement optimization algorithms directly in closed loop with physical systems. 
        In addition to a brief overview of selected continuous-time dynamical systems for optimization, our particular emphasis in this survey lies on closed-loop stability as well as the {robust} enforcement of physical and operational constraints in closed-loop implementations.
        To bypass accessing partial model information of physical systems, we further elaborate on fully data-driven and model-free operations.
        We highlight an emerging application in autonomous reserve dispatch in power systems, where the theory has transitioned to practice by now.
        We also provide short expository reviews of pioneering applications in communication networks and electricity grids, as well as related research streams, including extremum seeking and pertinent methods from model predictive and process control, to facilitate high-level comparisons with the main topic of this survey.
        }
    \end{abstract}


\end{frontmatter}

\onecolumn
\tableofcontents
\twocolumn

\section{Introduction}\label{sec:introduction}

Most advances in mathematical optimization in the past decades have been geared towards numerical implementations of iterative algorithms. The common viewpoint is that an optimization problem can be formulated, transformed, reduced, and relaxed, but ultimately the necessary steps to solve the problem rely purely on numerical linear algebra, which can be implemented and run on microprocessors. This offline computed solution is then used to reach and realize some decision.

This paradigm of optimization as a computational problem is almost synonymous with the field of management science and operations research \citep{dantzigLinearProgrammingExtensions1998,hillierIntroductionOperationsResearch2001,bertsimasDataModelsDecisions2004}, which has flourished ever since the inception of linear programming in the mid-20th century. Today, this kind of offline optimization is applied in various disciplines ranging from econometrics over statistical and machine learning \citep{bishopPatternRecognitionMachine2009} to optimal control \citep{bertsekasDynamicProgrammingOptimal2017}.

From a control perspective, solving an optimization problem offline (with known data) and implementing its output as a decision is a \emph{feedforward} approach. {\tb More precisely, the solution of an optimization problem is used as a set-point for a physical system, and feedback controllers are merely required to steer the physical system to this pre-computed state.\footnote{Beyond computing optimal setpoints, optimization has many more uses in control, such as the optimal design and tuning of controllers, but those topics are not the subject of this article.} However, due to discrepancies between the model used for set-point optimization and the real system, the realized system state is generally suboptimal. Since set-point optimization happens offline and without recourse to real-time measurements, this traditional separation may be interpreted as an \emph{open-loop} setup.}

In contrast, in this article, we {\tb review} a \emph{feedback} approach to constrained nonlinear optimization that drives a physical system towards an optimal steady state {\tb by using real-time measurements in the set-point optimization itself. More concretely, we review ideas that consider optimization algorithms as dynamical systems and cast them as closed feedback loops. Equivalently, these methods may be interpreted as ``optimization algorithms that incorporate real-time measurements'' or as ``feedback controllers that mimic optimization algorithms''.}\footnote{{\tb
            The boundary between {\tb ``feedforward'' and ``feedback optimization''} is not always clear cut. For instance, \ac*{mpc} uses feedback to achieve robustness{\tb ,} but relies on an accurate model for the formulation and solution of an optimal control problem at every iteration. {\tb The works reviewed in this article share in common that they do not solve the primary optimization problem numerically, but rather rely on the closed-loop system dynamics to settle to an optimizer. We refer to Appendix~\ref{sec:cl_opt_approaches} for a comparison.} %
            }}
            
{\tb The study of feedback-based optimization of a domain-independent concept is relatively recent, although it is historically rooted in applications such as congestion control in communication networks \citep{lowOptimizationFlowControl1999} and frequency control in power systems \citep{liConnectingAutomaticGeneration2016,molzahnSurveyDistributedOptimization2017}. In these contexts, existing control systems have been interpreted from an optimization perspective that allows for improved redesigns and retrofitting. For educational reasons, we review these two applications in Appendix~\ref{app:historical}, which have inspired much of the presented theory in this article.}
    
On a more abstract level, this type of \emph{closed-loop optimization} has been pursued mainly for three reasons:
\begin{enumerate}
    \item to {\tb approximately solve problems with inaccurate data and time-varying parameters},
    \item {\tb to achieve constraint satisfaction with minimal model-dependence, and} 
    \item to eliminate the need for exogenous setpoints and reference signals, especially in time-varying systems.
\end{enumerate}
These reasons resonate well with the general feedback and feedforward paradigms advocated in control textbooks \citep{doyleFeedbackControlTheory2009,franklinFeedbackControlDynamic2010}, and we further dwell on them hereafter.

\paragraph{Inaccurate problem data}

    {\tb In practice, optimization problems often} lack precise data. Parameters and states based on measurements and statistical inference are inherently inaccurate, and so is the solution of an optimization problem based on such data. The earliest attempts at addressing this issue have resulted in sensitivity analysis for optimization problems \citep{shapiroSensitivityAnalysisNonlinear1988} which asks how a solution changes as problem parameters vary. From a more practical perspective, {robust optimization} \citep{ben-talRobustOptimization2009} and {stochastic programming} \citep{bonnansConvexStochasticOptimization2019} offer ways to incorporate uncertainty in the problem. However, these approaches {\tb (which are primarily used offline and ahead of time)} are conservative {\tb in the sense that} they need to take into account the full set of possible problem instances {\tb even though, eventually, only one scenario will be realized}. {\tb Often, they also entail massive computational costs.}
{\tb In contrast, by optimizing a system in real-time and using feedback, one needs to react only to the actual realization of the underlying disturbance process. %
}

\paragraph{Constraint satisfaction}

{\tb Beyond the mitigation of the effect of modeling inaccuracies on the accuracy of the solution, closed-loop setups render constraint satisfaction partially} \emph{model-free}. A physical system naturally defines {\tb and enforces} a set of constraints.
These constraints may be of a different nature.
For instance, \emph{input saturation} may be due to limited actuator capabilities or to the actions of a low-level controller. 
The violation of \emph{hard physical limits} may trigger the immediate instability, failure, or destruction of the entire system, hence these constraints need to be satisfied at all times.
{\tb For offline optimization, these limits need to be accurately modeled. In a feedback approach, on the other hand, the closed-loop optimization algorithm automatically operates on a reduced set, e.g., a manifold, defined by the physical constraints (e.g., input-output relations, saturation limits, and hard physical limits).
Other constraints (e.g. thermal limits of physical components) are often fairly benign since they can be violated temporarily, but should be satisfied in the long run. These constraints are not enforced automatically and need to be addressed with proper control design.
In a feedback setting, measurements provide unequivocal information about the violation of such constraints, and this information can be exploited during runtime. 
}

\paragraph{Autonomous operation}

    {\tb Finally, c}ompared to standard feedback loops, closed-loop optimization setups can run {without external setpoints} or references. Instead, an economic objective can be directly optimized as long as a cost function can be specified and evaluated. This feature is particularly powerful in combination with the inherent constraint enforcement: whereas in classical control setups pre-computed setpoints have to be feasible (e.g., lie within actuator limits and be compatible with the physics of the plant), convergence to a feasible and optimal state is a defining aspect of feedback-based optimization.
    {\tb Additionally, in the presence of time-varying problem parameters that are hard to estimate and forecast, feedback-based optimization automatically results in the closed-loop system tracking the solution of the time-varying optimization problem.
    Furthermore, closed-loop setups open the possibility of optimizing a physical system in a \emph{data-driven} and \emph{model-free} manner. The input-output data collected in closed loop offers a rich source for representing, learning, and modulating the behaviors of a system.}

\bigskip

    {\tb In the remainder of this introduction, we discuss a simple, yet insightful, academic example of feedback-based optimization to offer a more concrete perspective on the points above and the technical issues that have been addressed by the research community in recent years. In particular, we illustrate a distinctive challenge in feedback-based optimization, which is the certification of the stability of the closed-loop interconnection between an optimization algorithm and the dynamical system. This topic has been the subject of substantial recent work which we review later in this article.
    }

For the sake of a concise presentation, {\tb in our pedagogical examples and in the tutorial sections of this article,} we restrict ourselves to continuous-time systems. {\tb In many cases,} analogous discrete-time or sampled-data results exist or can be derived similarly to the continuous-time case. {\tb Discrete-time models also raise questions about the efficient implementation of feedback-based optimization. Those issues often require domain-specific solutions which are not the main topic of this article.}

\subsection{Illustrative Examples}
The following idealized yet general example of a gradient system interconnected with a physical plant illustrates the idea {\tb at the center of this review} and {\tb some of the associated challenges}. {\tb Many of the assumptions made in the following are for illustrative purposes and are not needed in general.}

\begin{example}\label{ex:simp_grad}
    Consider a dynamic nonlinear plant
    \begin{align}\label{eq:simple_plant}
        \dot \zeta = {f}(\zeta, u) \qquad y = g(\zeta) + d \, ,
    \end{align}
    where $\zeta,u$ and $y$ are the state, input, and output, and $d$ denotes an additive disturbance. The vector field $f(\cdot, \cdot)$ and the map $g(\cdot)$ describe the process and output measurement, respectively.

    We assume that, {\tb for any fixed $u$,}
    the plant is asymptotically stable with fast-decaying transients such that, for every $u$,  there exists a unique steady state $\hat{h}(u)$ such that $0 = f(\hat{h}(u),u)$. Consequently, there also exists a \emph{steady-state map} $h(u) := g(\hat{h}(u))$. We assume $h$ to be continuously differentiable in $u$.\footnote{{\tb Existence of a (differentiable) steady-state map is generally shown using implicit function theorems \citep{dontchevImplicitFunctionsSolution2014}.}}

    We wish to drive the system to a steady state that minimizes a cost $\Phi(y)$ which is a function of the plant output $y$, hence to the solution of
    \begin{align*}
        \min_{u,y} & \quad \Phi(y) \\
        \subjto & \quad y = h(u) + d.    
    \end{align*}
    Given $h$ and $d$, we may equivalently minimize the \emph{reduced} cost $\tilde{\Phi}(u) := \Phi(h(u)+d)$ instead. For this purpose, we consider a simple gradient flow
    \begin{align}\label{eq:simple_grad_flow}
        \dot u = - \nabla \tilde{\Phi}(u)^T = - \nabla h(u)^T \nabla \Phi(h(u)+d)^T\, ,
    \end{align}
    where $\nabla h (u)$ is due to the chain rule applied to $\Phi(h(u)+d)$.

    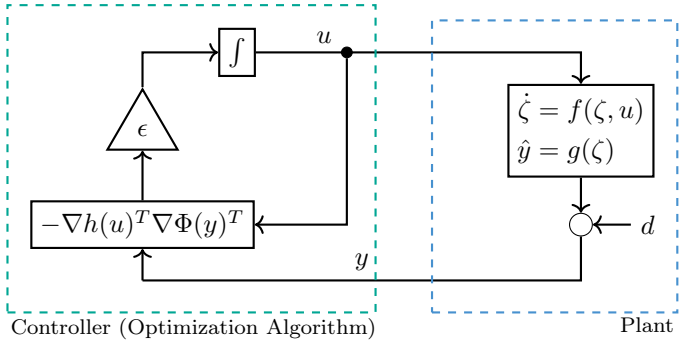
\begin{figure}[bt]
        \centering
        \begin{tikzpicture}
            \matrix[ampersand replacement=\&, row sep=0.3cm, column sep=.45cm] {
                \&  \node[none] (hook) {};\& \node[branch](br2){};\&
                \& \& \node[none] (fit1) {};                                           \\
                \node[gainup] (int) {$\epsilon$};
                \& \&  \& \& \& \&
                \node[block] (plant) {
                    $\begin{aligned}
                            \dot \zeta & = {f}(\zeta, u) \\
                            \hat{y}    & = g(\zeta)
                        \end{aligned}$
                };                                                 \\
                \node[block] (ctrl) {$ - \nabla h(u)^T \nabla \Phi(y)^T$};\& \&  \& \& \& \&
                \node[smallsum] (dist_sum){};                                                 \\
                \node[none](edge) {};    \& \&  \& \&  \& \& \node[none] (fit2) {}; \\
            };
            \node[block,  xshift=-.8cm] at (hook) (int2) {$\int$};

            \draw[connector] (ctrl.north)--(int.south);
            \draw[connector] (int.north)|-(int2.west);
            \draw[line] (int2.east)--(br2.center) node[near end, above]{$u$};
            \draw[connector] (br2.center)-|(plant.north) ;
            \draw[connector] (br2.south)|-(ctrl.east);
            \draw[connector] (plant.south)--(dist_sum.north);
            \draw[line] (dist_sum.south)|-(edge.center) node[near end, above ]{$y$};
            \draw[connector] (edge.center)--(ctrl.south);
            \draw[connector] ([xshift=.5cm]dist_sum.east)--(dist_sum.east) node[at start, right]{$d$};

            \node[fit=(edge)(br2)(ctrl)(int2), draw, dashed, inner sep= 3mm, thick, persiangreen] (fitb1){};
            \node[fit=(plant)(dist_sum)(fit1)(fit2), draw, dashed, inner sep=3mm, thick, tuftsblue] (fitb2) {};

            \node[below right, inner sep = .5mm] at (fitb1.south west) {\footnotesize Controller (Optimization Algorithm)};
            \node[below left, inner sep = .5mm] at (fitb2.south east) {\footnotesize Plant};
        \end{tikzpicture}
        \caption{Simple feedback-based gradient flow}\label{fig:simp_grad}
    \end{figure}

    The gradient flow \eqref{eq:simple_grad_flow} is a closed system {\tb (in the sense that it does not require an exogenous input or define an output)}. However, {\tb by identifying} $h(u)+d$ as the measurable output $y$ {\tb and $u$ as the input}, \eqref{eq:simple_grad_flow} can be easily transformed into an open system and interconnected with the plant \eqref{eq:simple_plant}, as shown in \cref{fig:simp_grad}. This yields the closed-loop dynamics
    \begin{align}\label{eq:ic_grad_system}
        \begin{split}
            \text{plant} & \begin{cases}
                \dot \zeta & =   {f}(\zeta, u) \\
                y          & = g(\zeta) + d
            \end{cases} \\
            \text{controller} & \begin{cases} \dot u & = - \epsilon \nabla h(u)^T \nabla \Phi(y)^T \, ,
            \end{cases}
        \end{split}
    \end{align}
    where $\epsilon > 0$ is a scalar control gain.

    Crucially, the controller in \eqref{eq:ic_grad_system} does not require explicit knowledge of $h$ (nor of $f$ or $g$). Instead, only the cost function gradient $\nabla \Phi(y)$ as well as steady-state input-output sensitivities $\nabla h(u)$ are required. Moreover, the additive disturbance $d$ does not need to be known or explicitly estimated and is fully rejected, i.e., an equilibrium is a critical point of $\Phi(h(u)+d)$, independently of the value of $d$.
    Notice that this is partly possible because the disturbance is additive and therefore it does not affect the input-output sensitivity of the system. We will discuss how to deal with uncertain sensitivities in Section~\ref{subsec:extension}.

    Thanks to the integral control structure of \cref{eq:ic_grad_system}, it can be easily seen that any equilibrium point $(\zeta^\star, u^\star)$ of \eqref{eq:ic_grad_system} is a steady state of the plant and satisfies $\nabla \tilde{\Phi}(u^\star)^T = \nabla h(u^\star)^T \nabla \Phi(h(u^\star)+d)^T = 0$. Therefore, $u^\star$ is a critical point of $\tilde{\Phi}$ (and a minimizer if $\tilde{\Phi}$ is convex).
    
    The systems \eqref{eq:simple_plant}, \eqref{eq:simple_grad_flow}, and \eqref{eq:ic_grad_system} can be understood from a \emph{singular perturbation} viewpoint \citep[Chap. 11]{khalilNonlinearSystems2002}: As $\epsilon \rightarrow 0^+$, the plant behavior is replaced by the algebraic map $h$, and the remaining dynamics \eqref{eq:simple_grad_flow} are the ``slow'' \emph{reduced system}. Conversely, on a fast timescale, on which $u$ and $d$ can be assumed to be constant, the plant dynamics \eqref{eq:simple_plant} are referred to as the ``fast'' \emph{boundary-layer system}.
\end{example}

One of the fundamental points left open by \cref{ex:simp_grad} is \emph{closed-loop stability}. The idea that plant dynamics in \cref{eq:ic_grad_system} need to be fast-decaying relative to the controller is indeed crucial as the following numerical example shows.

\begin{example}\label{ex:num_grad}
    Consider the objective $\Phi(y) = (y^2 - 1)^2$ which is illustrated in the top left panel of \cref{fig:num_grad} and has two isolated minima $\{-1, 1\}$. As plant \eqref{eq:simple_plant}, consider a single-input-single-output second-order plant governed by
    \begin{align*}
        \ddot{\zeta}  + a \dot{\zeta}+ b (\zeta - u) = 0 \,
    \end{align*}
    with $a = 2$, $b = 25$, and $y = \zeta$. There is no disturbance acting on the system. The plant is asymptotically stable, under-damped, and, at steady state, we have $y = \zeta = u$. Hence, the controller in \eqref{eq:ic_grad_system} takes the form
    \begin{align*}
        \dot u = - \epsilon \nabla \Phi(y)^T = - 4 \epsilon y (y^2 - 1) \, .
    \end{align*}

    \Cref{fig:num_grad} shows trajectories of the closed-loop system \cref{eq:ic_grad_system} for the same initial condition, but different values of the gain $\epsilon$, and comparing it to the ``algebraic'' gradient flow \eqref{eq:simple_grad_flow} given by $\dot u =  -\epsilon \nabla \Phi(h(u))^T = - 4 \epsilon u (u^2 - 1)$.

    \begin{figure}
        \centering
        \includegraphics[width=\columnwidth]{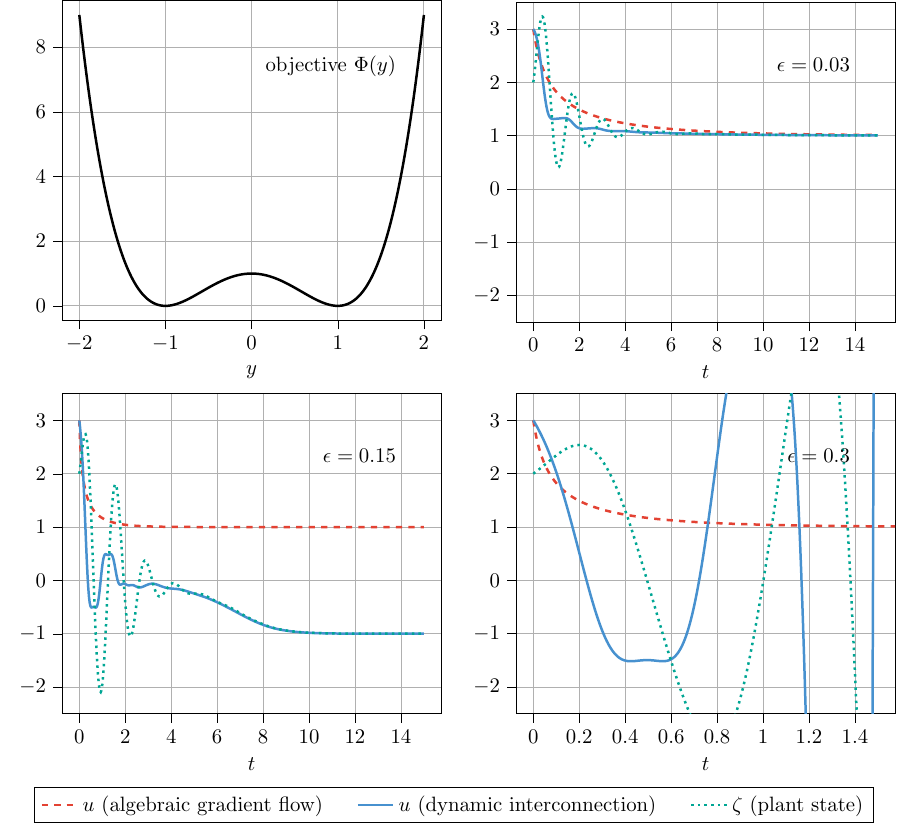}
        \caption{Illustrations for \cref{ex:num_grad} (top left: objective function; remaining panels: system trajectories for different control gains $\epsilon$)}\label{fig:num_grad}
    \end{figure}

    We observe, that for the given initial condition the algebraic gradient trajectory converges to the minimizer at 1. In contrast, the trajectories of the closed-loop system \eqref{eq:ic_grad_system} converge to either one of the two minimizers or diverge, depending on $\epsilon$. In other words, closed-loop stability of \eqref{eq:ic_grad_system} is not guaranteed, and even if it is, convergence may not be to the same minimizer as for \eqref{eq:simple_grad_flow}.
\end{example}

The issue of stability is central to many closed-loop optimization setups and work on this topic, including quantitative stability and robustness requirements, will be discussed in more detail in \cref{sec:cl_stab}.

\subsection{Organization}
The remainder of this article is structured as follows:
\Cref{sec:opt_dyn} {\tb provides a tutorial-like review of} various optimization algorithms {\tb from a dynamical systems perspective. This point of view is crucial for closed-loop optimization.}
    {\tb I}n \cref{sec:ofo}, we then review recent works on constrained feedback-based optimization.
We pay particular attention to different mechanisms to ensure closed-loop stability and to enforce constraints. 
{\tb We further clarify how to leverage data and exploration to facilitate model-free implementations in closed loop.}
    {\tb \Cref{sec:redispatch} presents a case study on the timely problem of optimal autonomous redispatch in electricity grid. This problem requires solving a nonlinear, nonconvex AC optimal power flow problem in closed-loop with the physical power grid and has motivated many implementations and recent advances in closed-loop optimization.}
Finally, \cref{sec:conc} concludes the article, highlights the unresolved problems in this domain, and presents exciting avenues for future research.

{\tb As mentioned before, in \cref{app:historical} we revisit some historical application examples from communication networks and power systems.
In \cref{sec:cl_opt_approaches}, for further contextualization, we {\tb give a high-level description of} other related approaches, such as \emph{extremum seeking}, real-time iterations, and \emph{modifier adaptation}, which perform optimization tasks at a system level within a closed-loop setup. We quickly discuss the main ideas of these well-established domains in order to provide historical context and a very broad overview of the adjacent research landscape.}

\section{Optimization Algorithms as Dynamical Systems}\label{sec:opt_dyn}

In recent years, renewed attention has been paid to the fact that many numerical optimization algorithms, {\tb and specifically first-order iterative algorithms,} can be interpreted as dynamical systems. 
{\tb In this section, we therefore provide an overview of different \emph{optimization dynamics} and review the main properties of these flows.} 
This perspective is essential to bridge the gap between algorithms and their implementation as feedback systems, which will be the focus of \cref{sec:ofo}. {\tb Apart from continuous-time optimization dynamics, discrete-time methods (for optimization or games) are also amenable to closed-loop implementations in a sampled-data setting, see the feedback equilibrium seeking in \citep{belgioioso2021sampled,belgioioso2022online,belgioioso2023tutorial}.}

\subsection{Gradient Flows}
\label{ssec:gradient_methods}

Arguably as old as differential calculus itself are methods that seek out (local) minima of a function {\tb $\Phi(x)$ over an unconstrained domain $x \in \mathbb{R}^n$} by following a descent direction. The most prototypical class of methods are gradient (or ``steepest descent'') schemes, e.g., $\dot{x} = -\nabla \Phi(x)$. Gradient methods exist in a wide variety of forms and contexts. Apart from $\bbR^n$, they can be defined on non-Euclidean manifolds \citep{brockettDynamicalSystemsThat1988,helmkeOptimizationDynamicalSystems1996,absilOptimizationAlgorithmsMatrix2008,zehnderLecturesDynamicalSystems2010} and on infinite-dimensional spaces \citep{luenbergerOptimizationVectorSpace1969,ambrosioGradientFlowsMetric2005}. In the absence of smoothness they can also be generalized to subgradient methods \citep{beckFirstOrderMethodsOptimization2017,clarkeOptimizationNonsmoothAnalysis1990}. On top of that, steepest descent methods are generally available as continuous-time flows, discrete-time algorithms, or stochastic processes \citep{borkarStochasticApproximationDynamical2008}.
For ease of exposition and in line with the rest of this article, we will limit ourselves to continuous-time {\tb deterministic} gradient flows.

At least for continuous-time negative gradient flows, convergence to minimizers appears plausible, if not tautological. However, closer inspection reveals important technical details, summarized in the following theorem {\tb which combines classical results available in most textbooks on calculus, dynamical systems, or optimization}.

\begin{theorem}\label{thm:grad}
    Let $\Phi: \bbR^n \rightarrow \bbR$ be continuously differentiable with locally Lipschitz derivative $\nabla \Phi$ such that, for some $c \in \bbR$, the sublevel set $\Phi^{-1}(c) = \{ x \in \bbR^n \, | \, \Phi(x) \leq c \}$ is compact. Then, the following statements hold for the gradient flow $\dot x = - \nabla \Phi(x)^T$:
    \begin{enumerate}
        \item\label{it:grad1} Trajectories starting in $\Phi^{-1}(c)$ converge to the set of critical points, i.e., points $x^\star$ such that $\nabla \Phi(x^\star) = 0$.

        \item\label{it:grad2} If $\Phi$ is analytic or convex, then every solution starting in $\Phi^{-1}(c)$ converges to a single point. {\tb If $\Phi$ is strongly convex or satisfies the Polyak-Lojasiewicz inequality, convergence is exponential.}

        \item\label{it:grad3} Every asymptotically stable equilibrium is a strict and isolated minimizer, and every local minimizer is stable. If $\Phi$ is analytic or convex, then the set of all (strict) minimizers is equivalent to the set of (asymptotically) stable equilibria.

        \item\label{it:grad4} If $\Phi$ is twice continuously differentiable, the stability of a critical point is partially governed by its Hessian: if $\nabla^2 \Phi(x^\star)$ is positive definite, $x^\star$ is a local minimizer and locally exponentially stable. If $\nabla^2 \Phi(x^\star)$ has at least one negative eigenvalue, $x^\star$ is unstable.
    \end{enumerate}
\end{theorem}

In the absence of local Lipschitz continuity of $\nabla \Phi$, uniqueness of trajectories is not guaranteed and a distinction between \emph{weak} and \emph{strong equilibria} has to be made~\citep{hauswirthProjectedDynamicalSystems2020,cortesDiscontinuousDynamicalSystems2008}.

Point \cref{it:grad1} in \cref{thm:grad} can be seen as a consequence of LaSalle's invariance theorem \citep{khalilNonlinearSystems2002}. Namely, it is immediate that $\Phi$ is non-increasing along gradient trajectories. Compactness of $\Phi^{-1}(c)$ is {\tb sufficient to} preclude unbounded trajectories that escape to the horizon. Further, compactness of $\Phi^{-1}(c)$ together with the continuity of $\Phi$ {\tb is sufficient to guarantee} lower boundedness and thus the existence of a minimizer.

In general, however, trajectories may not converge to a single point but to the entire set of critical points \citep{palisGeometricTheoryDynamical1982}. This pathological behavior is ruled out in \cref{it:grad2} if $\Phi$ is analytic, which guarantees the finite length of trajectories due to Lojasiewicz's inequality \citep{absilConvergenceIteratesDescent2005}, or if $\Phi$ is convex. %
{\tb The conditions that ensure exponential convergence are studied in \cite{bianchinOnlineOptimizationSwitched2020,frasconiMachineLearningKnowledge2016}.}

As indicated in \cref{it:grad3}, the stability of equilibria is related to their optimality, but an equivalence between the two requires additional assumptions. See \cite{absilStableEquilibriumPoints2006} for a proof and counterexamples.
Finally, \cref{it:grad4} follows since, if $\Phi$ is smooth enough, the stability of an equilibrium $x^\star$ can be analyzed by investigating the linearized dynamics $\dot x = - \nabla^2 \Phi(x^\star) x.$

\begin{example}[{\tb Variable Metric and Newton Gradient Flow}]\label{ex:var_metric_grad}
    A  straightforward degree of freedom for gradient flows is the use of a \emph{metric} $Q(\cdot)$ that maps every point $x \in \bbR^n$ to a square symmetric positive-definite matrix $Q(x)$. Under minor technical conditions (e.g., that $Q$ has a uniformly bounded condition number, see \citealt{hauswirthProjectedDynamicalSystems2020}) trajectories of the generalized gradient flow
    \begin{align}\label{eq:metric_grad_flow}
        \dot x = - Q(x) \nabla \Phi(x)^T
    \end{align}
    converge to the set of critical points of $\Phi$. Namely, the use of $Q$, does neither alter the equilibrium points nor their attractivity and stability, but only the trajectories. This feature is illustrated in \cref{fig:sim_grad_metrics1} which shows gradient trajectories for the same nonconvex potential function but for different metrics.

    If $Q$ is constant, \eqref{eq:metric_grad_flow} is equivalent to a Euclidean gradient flow in linearly transformed coordinates. Namely, if $Q = V^T \Lambda V$ is the eigenvalue decomposition of $Q$, we can define the coordinate transformation $y := \sqrt{\Lambda}^{-1} V x$ where $\sqrt{\Lambda}$ denotes the diagonal matrix of square root eigenvalues. Then, it can be easily shown that \eqref{eq:metric_grad_flow} is equivalent to $\dot y = -\nabla \hat{\Phi}(y)^T$ with $\hat{\Phi}(y) := \Phi(V^T \sqrt{\Lambda} y)$.

    If $Q$ is not constant, then, from {\tb a} differential-geometric viewpoint, $Q$ can be interpreted as a matrix representation of a non-Euclidean \emph{Riemannian metric} on $\bbR^n$. Thus, $Q$ endows $\bbR^n$ with a non-flat geometry \citep{leeRiemannianManifoldsIntroduction1997}.

    As a special case of a non-constant metric, assume that $\Phi$ is twice differentiable and strongly convex. Then, $Q$ can be chosen to be the inverse of the Hessian $\nabla^2 \Phi$. This results in a continuous-time version of the classic \emph{Newton method}, also referred to as ``Newton gradient flow'' \citep[Chap. 9.3]{jongenNonlinearOptimizationFinite2001}. However, unlike the iterative Newton method, the continuous-time flow does not exhibit an inherently faster convergence rate compared to other gradient flows. For the inverted Hessian metric, convergence is \emph{isotropic}, i.e., the same from all directions. This property counteracts ill-conditioning of the objective function as illustrated in \cref{fig:sim_grad_newt}.

    \begin{figure}[!t]
        \centering
        \includegraphics[width=\columnwidth]{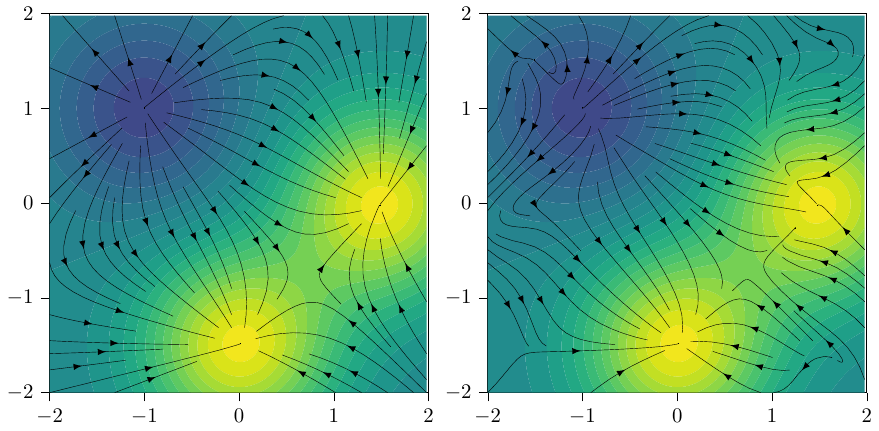}
        \caption{Gradient trajectories for a non-convex objective function. Trajectories under the Euclidean metric (left) and a generic variable metric (right) differ significantly. The critical points (i.e., the minima, maximum, and saddle-point) and their stability properties are unaffected by the choice of metric.}\label{fig:sim_grad_metrics1}
    \end{figure}

    \begin{figure}[!t]
        \includegraphics[width=\columnwidth]{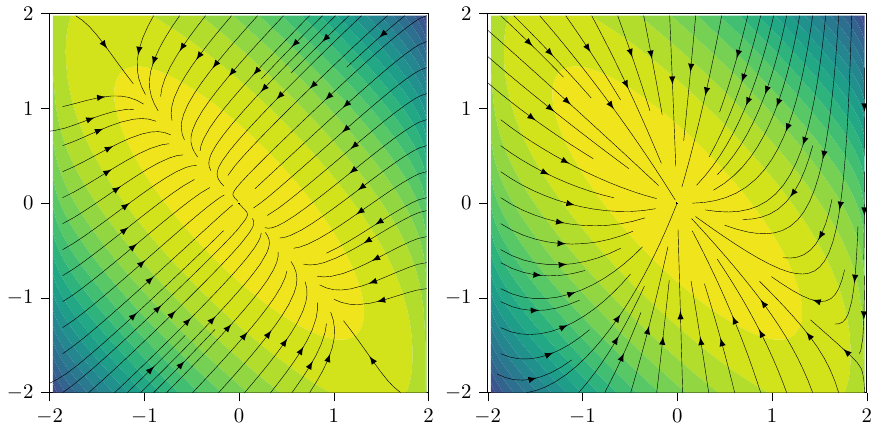}
        \caption{Gradient trajectories for strongly convex, but ill-conditioned objective. The trajectories under the Euclidean metric (left) quickly approach a subspace on which the objective is almost flat, and then converge only slowly to the global optimizer. Trajectories under the ``Newton metric'' (right), approach the global optimizer isotropically, unaffected by the ill-conditioning of the objective.}\label{fig:sim_grad_newt}\label{fig:sim_grad_metrics}
    \end{figure}

    Finally, if $Q$ is sparse, then the \emph{sparsity} pattern induces algebraic structure that can often be exploited for the purpose of a distributed implementation {\tb \citep{bulloLecturesNetworkSystems2022}.}
\end{example}

\subsection{Projected Gradient Flows}\label{sec:proj_grad}

When unilateral (i.e., inequality) constraints restrict the search domain of the optimization problem, it is often possible to rely on gradient flows by introducing suitable projection mechanisms. In a computational context, this is particularly true if the projection onto a given constraint set is easy to evaluate numerically.

Consider the  constrained optimization problem
\begin{align}\label{eq:proj_prob}
    \minimize \quad \Phi(x) \qquad \subjto \quad x \in \calX \, ,
\end{align}
where $\calX \subset \bbR^n$ is closed convex and non-empty, and $\Phi$ is {\tb continuously differentiable}. The classical \emph{projected gradient descent} iteration to solve this problem takes the form
\begin{align}\label{eq:proj_grad_desc}
    x^{k+1} = P_{\calX}\left( x^k - \alpha^k \nabla \Phi(x^k)^T \right) \, ,
\end{align}
where $P_\calX(y) := \arg\min_{x \in \calX} \| x - y \|$ denotes the Euclidean minimum norm projection onto $\calX$, and $\{ \alpha^k \}$ is a sequence of step sizes.

\begin{figure}[tb]
    \centering
    \includegraphics[width=.48\columnwidth]{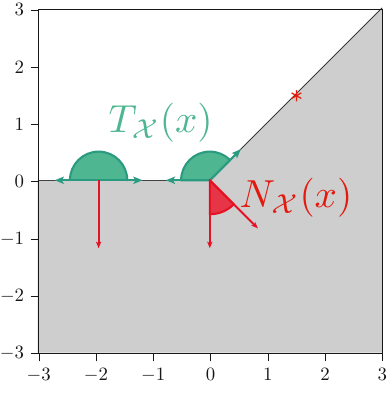}
    \includegraphics[width=.48\columnwidth]{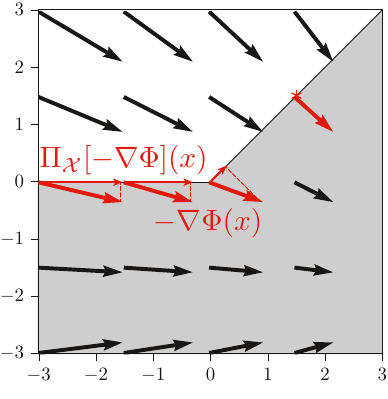}
    \caption{Left: Examples of tangent and normal cones of a set $\calX$ (gray area is infeasible). Right: Projected vector field onto the tangent cone. See \cref{fig:cstr_pgrad} for the resulting projected gradient trajectory.}\label{fig:tgt_cones}
\end{figure}

By choosing infinitesimally small step-sizes, the continuous-time limit of \eqref{eq:proj_grad_desc} is a \emph{projected gradient flow}
\begin{align}\label{eq:proj_grad_flow}
    \dot{x} = \Pi_\calX \left[ -\nabla \Phi(x)^T \right](x),
\end{align}
where $\Pi_\calX[v](x)$ denotes the projection of the vector $v$ onto the tangent cone $T_\calX(x)$ of $\calX$ at $x$, i.e., $\Pi_\calX[v](x) = \arg\min_{w \in T_\calX(x)} \| v - w \|$. As illustrated in \cref{fig:tgt_cones}, the tangent cone is the set of all directions starting at $x$ which point inward into the set $\calX$. For convex sets, it takes the general form $T_\calX(x) = \cl \{ d \, | \, d = \alpha (y -x), y \in \calX, \alpha \geq 0 \}$. 

In contrast to the discrete-time system \eqref{eq:proj_grad_desc}, continuous-time projected gradient flows can be generalized extensively. In particular, convexity of $\calX$ is not required. More precisely, whereas $P_\calX$ requires $\calX$ to be convex to be well-defined, $\Pi_\calX$ is generally well-defined since $T_\calX(x)$ is non-empty and closed convex for a large class of non-convex sets called \emph{(Clarke) regular} (see \citealt{rockafellarVariationalAnalysis2009}).

\begin{figure}[t!]
    \centering
    \begin{subfigure}[t]{.48\columnwidth}
        \centering
        \includegraphics[width=\columnwidth]{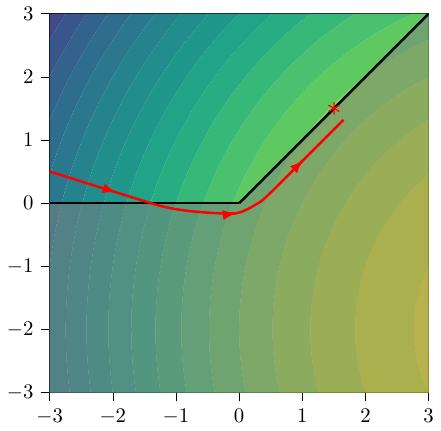}
        \caption{Penalty Function}\label{fig:cstr_pen}
    \end{subfigure}
    \begin{subfigure}[t]{.48\columnwidth}
        \centering
        \includegraphics[width=\columnwidth]{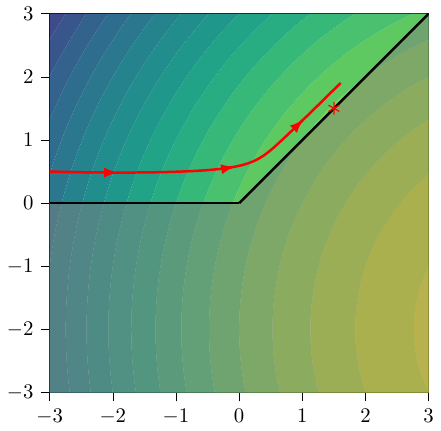}
        \caption{Barrier Function}\label{fig:cstr_bar}
    \end{subfigure} \\
    \begin{subfigure}[t]{.48\columnwidth}
        \centering
        \includegraphics[width=\columnwidth]{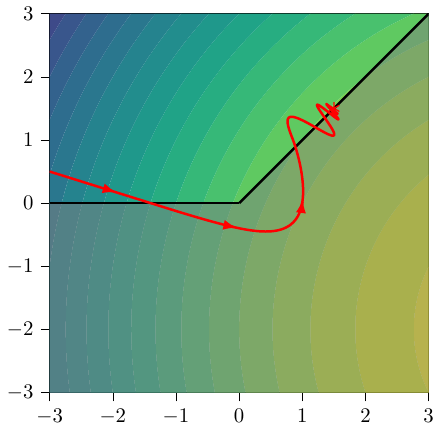}
        \caption{Saddle-Flow}\label{fig:cstr_sadd}
    \end{subfigure}
    \begin{subfigure}[t]{.48\columnwidth}
        \centering
        \includegraphics[width=\columnwidth]{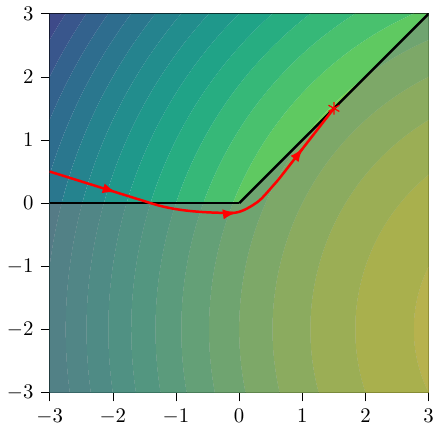}
        \caption{Augmented Saddle-Flow}\label{fig:cstr_asadd}
    \end{subfigure}\\
    \begin{subfigure}[t]{.48\columnwidth}
        \centering
        \includegraphics[width=\columnwidth]{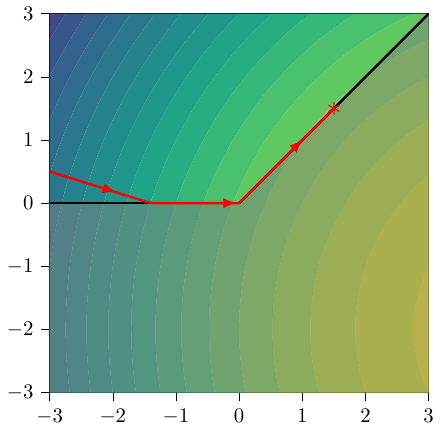}
        \caption{Projected Gradient Flow}\label{fig:cstr_pgrad}
    \end{subfigure}
    \begin{subfigure}[t]{.48\columnwidth}
        \centering
        \includegraphics[width=\columnwidth]{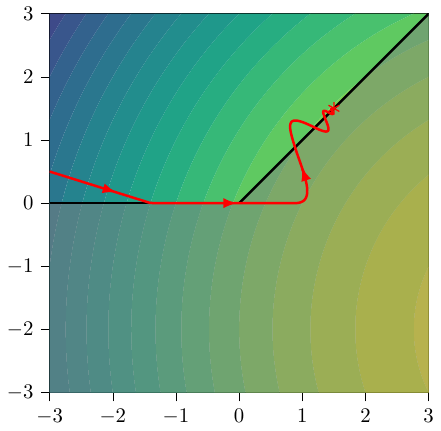}
        \caption{Mixed Saddle-Flow}\label{fig:cstr_msadd}
    \end{subfigure}
    \caption{Behavior of different constraint enforcement mechanisms\\
        All panels show the minimization of a quadratic function subject to two constraints $x_2 \geq 0$ and $x_2 \geq x_1$ (the grayed out area is infeasible).
        Penalty (a) and barrier (b) functions allow for smooth outer and inner approximations of constraints with an unconstrained gradient flow. Saddle-point flows (c) enforce constraints only asymptotically by integrating constraint violation over time. Augmenting saddle-point flows with a penalty term can improve convergence (d) as explained in \cref{ex:sadd_flow}. Projected gradient flows (e) enforce constraints directly by projection, which results in non-smooth trajectories. Individual constraints can also be enforced with a combination of these mechanisms, e.g., as in (f) with a projection for $x_2 \geq 0$ and dualization (saddle-flow) for $x_2 \geq x_1$, as in \cref{ex:sadd_flow}.}\label{fig:cstr_num}
\end{figure}

\Cref{fig:cstr_pgrad} illustrates the qualitative behavior of projected gradient flows. Namely, in the interior of the feasible set, trajectories follow the gradient direction, whereas at the boundary, trajectories follow the steepest \emph{feasible} descent direction. Projected gradient flows are therefore inherently discontinuous systems, and their study requires tools from non-smooth analysis~\citep{aubinDifferentialInclusionsSetvalued1984,hauswirthProjectedDynamicalSystems2020}.
{\tb It is worth remarking that projected gradient flows can be modified in order to recover continuity \citep{wang2000globallyProjectedDynamicalSystems} and in order to ensure that the flow is well defined outside of the feasible region
\citep{jordan2022nonsmoothDynamicalSystems}. 
These properties often come with additional convexity assumptions.
An exception is the modification presented in \cite{cortes2022controlBarrierFunction}, where control barrier function theory is employed to achieve both these properties in a non-convex setting.
Another strategy that handles non-convexity is to leverage anti-windup schemes to smoothly and robustly approximate projected gradient flows \citep{hauswirthAntiWindupApproximationsOblique2020a,hauswirthDifferentiabilityProjectedTrajectories2020,hauswirthRobustImplementationProjected2020}, see \cref{ex:awa}.
}

{\tb Projected gradient flows {\tb extend} properties from their unconstrained counterparts. For instance, similarly to \cref{thm:grad}, trajectories of \eqref{eq:proj_grad_flow} converge to the set of critical points (in this case, Karush-Kuhn-Tucker points of \eqref{eq:proj_prob}), and stability and optimality can be related analogously to \cref{it:grad3} in \cref{thm:grad}~\citep{clarkeOptimizationNonsmoothAnalysis1990,hauswirthProjectedDynamicalSystems2020}.}

\begin{remark}[{\tb Discretization \& Proximal-Point Algorithm}] As mentioned above, the continuous-time model \eqref{eq:proj_grad_flow} can be obtained from \eqref{eq:proj_grad_desc} in the limit as $\alpha \rightarrow 0^+$ where $\alpha$ is a constant step size. Conversely, \eqref{eq:proj_grad_desc} can be interpreted as a \emph{forward Euler} discretization of \eqref{eq:proj_grad_flow}. Analogously, a \emph{backward Euler} scheme for \eqref{eq:proj_grad_flow} yields the implicit form $x^{k+1} = P_{\calX}\left( x^k - \alpha^k \nabla \Phi(x^{k+1})^T \right)$. Comparing optimality conditions, it is easy to show that this discretization is a special case of the more general \emph{proximal-point algorithm}~\citep[Chap. 27.1]{beckFirstOrderMethodsOptimization2017}.
\end{remark}

\begin{remark}[{\tb Non-differentiable objective \& subgradients}]
    Convex projected gradient flows also fall into the category of \emph{subgradient flows}: if $\calX$ is convex, we may consider the minimization of $\Phi(x) + I_\calX(x)$ where $I_\calX$ denotes the indicator function of the set $\calX$. Since $I_\calX$ is not differentiable, instead of a gradient flow, we need to resort to the \emph{subgradient inclusion} $\dot x \in - \nabla \Phi(x)^T - N_\calX(x)$. (In particular, the subgradient of $I_\calX$ is given by the normal cone $N_\calX(x)$ of $\calX$ at $x$, i.e., $\partial I_\calX(x) = N_\calX(x)$.) This differential inclusion is also referred to as \emph{differential variational inequality}~\citep{aubinDifferentialInclusionsSetvalued1984} and can be shown to be equivalent to \eqref{eq:proj_grad_flow} (i.e., admit the same trajectories).
\end{remark}

\begin{remark}[Variable Metric and Projected Newton Gradient Flow]
Analogously to unconstrained gradient flows discussed in \cref{ex:var_metric_grad}, projected gradient flows can also be defined using a variable metric $Q(\cdot)$. 
{\tb This degree of freedom allows, in particular, the definition of \emph{projected Newton gradient flows}.}
This modification, however, requires a generalization of $\Pi_\calX$ in \eqref{eq:proj_grad_flow} to take into account the effects of the metric (see \citealt{hauswirthProjectedDynamicalSystems2020}).
Namely, the variable-metric projected gradient flow takes the form
\begin{align*}
    \dot x = \Pi_\calX^Q \left[- Q(x) \nabla \Phi(x)^T \right] \,,
\end{align*}
where the ``oblique'' projection operator
\begin{align}
    \Pi_\calX^Q [v](x) := \arg \min_{w \in T_{\calX}(x)} (w -v)^T Q(x) (w - v) \label{eq:skew_proj_op}
\end{align}
has to be employed.
\end{remark}

\subsection{Primal-Dual Saddle-Point Flows}\label{sec:saddle}

Simply speaking, under weak technical assumptions, solutions of a constrained optimization problem are saddle-points of the associated \emph{Lagrangian}.
For this reason, dynamical systems that seek out saddle-points rather than extrema of a function are of particular interest for constrained optimization.

As a general yet basic setup, consider a differentiable function $L(x,\mu)$ that is convex in $x$ for every $\mu$, concave in $\mu$ for all $x$, and either strictly convex in $x$ or strictly concave in $\mu$. Then trajectories of the system
\begin{align}\label{eq:basic_sadd}
    \dot x = - \nabla_x L(x, \mu)^T \qquad \dot \mu = \nabla_\mu L(x,\mu)^T
\end{align}
converge to a saddle-point of $L$, i.e., a point $(x^\star, \mu^\star)$ such that $L(x, \mu^\star) \geq L(x^\star, \mu^\star) \geq L(x^\star, \mu)$ for all $x$ and all $\mu$. In particular, \cref{eq:basic_sadd} consists of a gradient descent in the \emph{primal} variables $x$ and a gradient ascent in the \emph{dual} variables $\mu$.

Historically, \emph{saddle-point flows} have primarily been considered in the context of nonlinear circuit analysis \citep{braytonTheoryNonlinearNetworks1964,smaleMathematicalFoundationsElectrical1972}, but their potential for optimization has been observed even before that \citep{arrowStudiesLinearNonlinear1958,koseSolutionsSaddleValue1956}.
Although saddle-point flows have been studied throughout the years \citep{venetsContinuousAlgorithmsSolution1985,blochGeometrySaddlePoint1992}, they have recently become a topic of intense interest due their importance for distributed network optimization (see  \citealt{cortesDistributedCoordinationNonsmooth2019,feijerStabilityPrimalDual2010} and the example in \cref{ex:netw_cong} on control of communication networks).

\begin{example}[{\tb Saddle-Point Flow for Optimization with Equality Constraints}]\label{ex:pd_eqcstr_basic} Consider the following linearly constrained problem
    \begin{align}\label{eq:pd_eqcstr_prob}
        \minimize \quad \Phi(x) \qquad \subjto \quad A x = b
    \end{align}
    where $A \in \bbR^{m \times n}$ and $b \in \bbR^m$, and $\Phi$ is strictly convex. The Lagrangian of \eqref{eq:pd_eqcstr_prob},  given by
    \begin{align*}
        L(x, \mu) := \Phi(x) + \mu^T (A x - b) \,
    \end{align*}
    is strictly convex in $x$ {\tb (since $\Phi$ is strictly convex)} and linear (thus concave) in $\mu$. Consequently, the saddle-point flow \eqref{eq:basic_sadd} is globally convergent and takes the form
    \begin{align}\label{eq:basic_sadd_lin}
        \dot x = - \nabla \Phi(x)^T - A^T \mu \qquad \dot \mu = A x - b \, .
    \end{align}
    Notice that any equilibrium $(x^\star, \mu^\star)$ of \eqref{eq:basic_sadd_lin} satisfies $0 = \nabla \Phi(x^\star)^T + A^T \mu^\star$ and $0 = A x^\star - b$ and thus the Karush-Kuhn-Tucker condition which are (for convex problems) necessary and sufficient for optimality of $x^\star$.
\end{example}

Like gradient flows, saddle-point flows can be modified through the use of positive definite metrics $Q_p, Q_d$ as
\begin{align}\label{eq:basic_sadd_q}
    \dot x = - Q_p \nabla_x L(x, \mu)^T \qquad \dot \mu = Q_d \nabla_\mu L(x,\mu)^T \, .
\end{align}
Often, $Q_p$ and $Q_d$ are chosen to be constant diagonal matrices that speed up or slow down convergence in specific directions. Interestingly, in the limit case $Q_p \gg Q_d$, one recovers a
differential-algebraic system
\begin{align*}
    x^\star(\mu) = \underset{x}{\arg\min}\, L(x, \mu) \qquad \dot{\mu} = Q_d \nabla_\mu L(x^\star(\mu), \mu)^T  \, ,
\end{align*}
which is a continuous-time \emph{dual ascent}. This idea of replacing the primal gradient flow with an explicit (algebraic) minimization can also be applied partially to a subset of variables. {\tb This technique was applied, e.g., in \cite{liConnectingAutomaticGeneration2016} and is illustrated in the example in \cref{ex:opt_freq} for frequency control in power systems.}

Convergence proofs for \eqref{eq:basic_sadd} and its generalization usually rely on the monotonicity of the vector field $[- \nabla_x L(x, \mu) \, \, \nabla_\mu L(x, \mu)]$. This type of argument, however, does not generalize beyond convex-concave saddle-point flows {\tb (see the discussion in \citealt{cherukuri2017saddle})}. Even relaxing strict convexity or concavity requirement of either $x$ or $\mu$, respectively, is problematic since this may lead to oscillations \citep{holdingConvergenceSaddlePoints2014}.
{\tb We will see in \cref{ex:sadd_flow} that augmentation of the Lagrangian is instrumental in recovering stability in these cases.}

Differentiability of the saddle-function $L$ is not generally required. In fact, \eqref{eq:basic_sadd} can be formulated in terms of subgradients if $L$ is not differentiable~\citep{venetsContinuousAlgorithmsSolution1985,goebelStabilityRobustnessSaddlepoint2017}.
In particular, \eqref{eq:basic_sadd} can be generalized to include projections on both $x$ and/or $\mu$ \citep{cherukuriAsymptoticConvergenceConstrained2016,cherukuri2017saddle,cherukuriRoleConvexitySaddlePoint2017,steginkConvergenceProjectedPrimalDual2018,hauswirthLimitBehaviorRole2020}. This possibility is particularly important to deal with inequality constraints, as the following example shows.

\begin{example}[{\tb Projected Saddle-Point Flow for Optimization with Inequality Constraints}]\label{ex:sadd_flow}
We revisit \cref{ex:pd_eqcstr_basic} to show how projected saddle-point flows can be used to deal with inequality constraints and to discuss augmentations based on penalty terms to improve convergence.
    Instead of \eqref{eq:pd_eqcstr_prob}, consider the problem
    \begin{align}\label{eq:sadd_prob}
        \begin{split}
            \minimize \quad & \Phi(x) \\
            \subjto \quad & x \in \calX \\
            & g(x) \leq 0 \, ,
        \end{split}
    \end{align}
    where $\Phi$ and $g$ are convex (but not necessarily strictly convex) and continuously differentiable. Further, let $\calX \subset \bbR^n$ be non-empty and closed convex. We define the \emph{partial} Lagrangian $L: \calX \times \bbR^m_{\geq 0} \rightarrow \bbR$ of \eqref{eq:sadd_prob} as
    \begin{align}
        L(x, \mu) := \Phi(x) + \mu^T g(x) \, ,
    \end{align}
    and note that $\mu$ must lie in the non-negative orthant $\bbR^m_{\geq 0}$ because it is associated with an inequality constraint.

    To find a saddle-point of $L$ on the set $\calX \times \bbR^m_{\geq 0}$ we use the continuous-time projection formalism introduced in \cref{sec:proj_grad} for projected gradient flows. Namely, we consider the \emph{projected saddle-point flow}
    \begin{align}\label{eq:basic_sadd_proj}
    \begin{split}
        \dot   x  = \Pi_{\calX} \big[\overbrace{- \nabla_x L(x, \mu)^T}^{\mathclap{- \nabla \Phi(x)^T - \nabla g(x)^T \mu}} \big](x) \\
        \dot \mu = \Pi_{\bbR_{\geq 0}^m} \big[\underbrace{ \nabla_\mu L(x,\mu)^T}_{ \mathclap{g(x)}}\big](\mu)
    \end{split}
    \end{align}
    where $\Pi_{\calX} \left[ w \right](x) $ and  $\Pi_{\bbR^m_{\geq0}} \left[ v\right](x) $ project $w$ and $v$ onto the tangent cone of $\calX$ and on the non-negative orthant $\bbR^m_{\geq 0}$ at $x$ and $\mu$, respectively. Consequently, trajectories of \eqref{eq:basic_sadd_proj} cannot leave $\calX \times \bbR^m_{\geq 0}$.

    Importantly, two different constraint enforcement mechanisms are at play: on one hand, the constraint $x \in \calX$ is enforced directly by projection, similarly to the projected gradient flow in \cref{sec:proj_grad}; on the other hand, the constraint $g(x) \leq 0$ is enforced by dualization. Namely, the dual variable $\mu$ is updated in response to a constraint violation $g(x) > 0$ and converges to a dual solution of~\eqref{eq:sadd_prob}.

    As shown, e.g., by \cite{goebelStabilityRobustnessSaddlepoint2017}, under weak technical assumptions and if $\Phi$ is strictly convex, trajectories of \eqref{eq:basic_sadd_proj} are guaranteed to converge to a Karush-Kuhn-Tucker point (and thereby to a global optimizer) of \eqref{eq:sadd_prob}.

    If $\Phi$ is non-strictly convex, then trajectories of \eqref{eq:basic_sadd_proj} converge to the optimizer of \eqref{eq:sadd_prob} if, instead of $L$, a penalty augmented Lagrangian of the form
    \begin{align}\label{eq:sadd_primal_aug}
        L^a(x, \mu) := \Phi(x) + \mu^T g(x) + \phi(x)
    \end{align}
    is used, where $\phi$ is a penalty function. In practice, $\phi$ is often chosen to be of the form $\phi(x) := \tfrac{\rho}{2} \| \max\{g(x), 0 \} \|^2$, but other constructions are possible to guarantee asymptotic convergence \citep{hauswirthLimitBehaviorRole2020}.

    If $\Phi$ or $g$ are non-convex, convergence of trajectories to a saddle point is generally not guaranteed. For this reason, \cite{bernsteinOnlinePrimalDualMethods2019,tang2018feedback,dallaneseOptimalPowerFlow2018,tangRunningPrimalDualGradient2018} have used an additional regularization on the dual variables of the form
    \begin{align}\label{eq:sadd_dual_aug}
        L^b(x, \mu) := \Phi(x) + \mu^T g(x) + \phi(x) - \tfrac{\hat{\rho}}{2} \| \mu\|^2 \, ,
    \end{align}
    which renders $L^b$ strongly concave in $\mu$.
    Using this modification and a large enough $\hat{\rho}$, the corresponding saddle-point flow is convergent, even with exponential stability guarantees. However, saddle-points of $L^b$ do generally not coincide with Karush-Kuhn-Tucker points of \eqref{eq:sadd_prob} anymore {\tb (unless the KKT point lies in the interior of the feasible set, in which case $\mu$ is zero and the dual regularization term disappears, see \citealt{bianchinTimeVaryingOptimizationLTI2021}).} In particular, notice that the dual update takes the form
    $\dot{\mu} = \Pi_{\bbR^n_{\geq 0}} \left[ g(x) + \hat{\rho} \mu \right]$
    and thus, at an equilibrium, it holds that $\mu_i^\star (g_i(x^\star) + \hat{\rho} \mu^\star_i) = 0$. In the absence of dual regularization, the same equilibrium condition simplifies to $\mu_i^\star g_i(x^\star) = 0$ which is the \emph{complementary slackness} part of the KKT optimality conditions of \eqref{eq:sadd_prob}.
    {\tb The error between the saddle points of the augmented Lagrangian and the solution of the constrained optimization problem can be bounded as a function of the augmentation term \citep{koshalMultiuserOptimizationDistributed2011}.}
\end{example}

\subsection{Comparison of Constraint-Enforcement Mechanisms}\label{sec:cstr_sum}

{\tb In the previous subsections we have presented several flows that correspond to optimization dynamics. Here, we compare and contrast the behavior of these different mechanisms based on how they handle constraints that define the feasible region of the optimization problem, as this is a crucial aspect in online and closed-loop applications.

Gradient flows, as defined in \cref{ssec:gradient_methods}, do not naturally allow for unilateral (i.e., inequality) constraints.
A relatively easy and widely applicable way to incorporate them is the addition of \emph{penalty} or \emph{barrier} terms to the objective. 
However, both approaches by themselves can enforce constraints only approximately. 

More precisely, for a constraint of the form $g(x) \leq 0$ where $g: \bbR^n \rightarrow \bbR$ is continuously differentiable, a common penalty function is the squared 2-norm of the constraint violation vector, i.e., $\phi(x) = \tfrac{\rho}{2} \|\max \{ g(x), 0 \}\|^2$ where $\rho > 0$ denotes a scaling parameter. Many variations, including different norms on constraint violations are possible. 
The common feature of penalty function lies in the fact that they technically allow for constraint violations and thus produce an ``outer'' approximation of the feasible region, i.e., minimizers of a penalty-augmented cost function $\Phi(x) + \phi(x)$ do not generally satisfy $g(x) \leq 0$.\footnote{A notable exception are so-called \emph{exact penalty methods} that transform a constrained optimization problem into an unconstrained one without changing the location of minimizers, albeit at the expense of smoothness or other technical drawbacks \citep{dipilloExactPenaltyFunctions1989}.}
This is illustrated in \cref{fig:cstr_pen}.

Barrier functions, on the other hand, can be used to enforce constraints strictly, i.e., without allowing for any violation. For this purpose, a barrier function $\psi(\cdot)$ for the constraint $g(x) \leq 0$ needs to be such that for $x \rightarrow x^\star$ with $g(x^\star) = 0$ we have $\psi(x) \rightarrow \infty$. A common example satisfying this condition are negative log-barriers of the form $\psi(x) = - \tfrac{1}{\mu} \log(g(x))$ which are important for interior-point methods for constrained convex programming \citep{nesterovInteriorpointpolynomialalgorithms1994}.
Barrier functions therefore achieve an ``inner'' (i.e., conservative) approximation of the feasible region, as illustrated in \cref{fig:cstr_bar}.

Both approaches are widely applicable (under minor technical assumptions) and do not require convexity of the constraints.
In general, steeper penalties or barriers (i.e., $\rho \rightarrow \infty$ or $\mu \rightarrow \infty$) will lead to more precise results. However, excessively steep augmentations can lead to stability issues when implemented numerically or in closed loop with a dynamical system (see \cref{sec:cl_stab}).

Constraint enforcement by (infinitesimal) projection, as for continuous-time projected gradient flows presented in \cref{sec:proj_grad}, is mathematically well-posed and works in very general settings. 
Furthermore, constraints are represented exactly and they are satisfied at all times (see \cref{fig:cstr_pgrad}).
However, we will see in \cref{sec:cstr_enforc} that producing these continuous-time projected flows on a physical system is not straightforward. In some instances, the projection behavior emerges naturally from physical saturation and anti-windup control. 
When the constraint is not enforced by the plant, careful design of the discrete-time implementation is required in order to maintain computational tractability and preserve the stability properties without introducing additional convexity assumptions.

Dualization of constraints leads to saddle-point flows where dual variables are computed by integrating the constraint violation over time. Hence, transient constraint violations are generally unavoidable, even in the presence of augmentation terms (see \cref{fig:cstr_sadd} and \ref{fig:cstr_asadd}). Tuning can be difficult, especially for nonlinear problems. Suboptimal parameter choices can lead to severely under- or over-damped transients that may venture far outside the feasible domain, which is undesirable in online and closed-loop applications. This problem gets only more challenging for high-dimensional and ill-conditioned problems.

In theory, each constraint (in functional form) can be enforced with one of these mechanisms independently of the other constraints. For example, in \cref{fig:cstr_msadd}, the constraint $x_2 \geq 0$ is enforced by projection whereas $x_2 \geq x_1$ is dualized.
As we will see in \cref{sec:cstr_enforc}, all these constraint enforcement methods can be applied in an online feedback setup. This design freedom is particularly useful in control setups where the real-world nature of constraints can dictate the appropriate enforcement mechanism. For instance, barrier functions may be considered for constraints that may not be violated under any circumstances. Constraints that are naturally enforced by physical saturation, mechanical constraints or similar are best represented by projections. Dualization in combination with a penalty term is particularly helpful to enforce constraints asymptotically and often allows for distributed implementations.}

\subsection{Time-Varying Online Optimization}\label{sec:tv_opt}

A topic that is central to online optimization, in open or closed loop, is the study of problems that vary over time. In recent years, this topic has garnered significant interest because of its relevance to many applications in control, robotics, machine learning, and others \citep{simonettoTimeVaryingConvexOptimization2020,dallaneseOptimizationLearningInformation2020,hazan2022introduction}.
The focus has been the development of algorithms that can track the solution of a non-stationary optimization problem with performance guarantees. Historically, two perspectives can be distinguished.

On the one hand, \cite{zinkevichOnlineConvexProgramming2003, hallOnlineConvexOptimization2015,jadbabaieOnlineOptimizationCompeting2015,besbesNonStationaryStochasticOptimization2015,lesage-landryPredictiveOnlineConvex2020,shiOnlineOptimizationMemory2020,nonhoffOnlineGradientDescent2020} and others frame time-varying optimization as an \emph{iterative learning problem}: at every iteration $k$ an agent chooses an action $x_k$ and subsequently a convex function $\Phi_k$ is revealed. The agent's goal is to minimize her \emph{regret}, i.e., some measure of accumulated suboptimality.

On the other hand, \cite{simonettoDecentralizedPredictionCorrectionMethods2017, rahiliDistributedContinuousTimeConvex2017,rahiliDistributedConvexOptimization2015,tangRunningPrimalDualGradient2018,tang2018feedback,fattahiAbsenceSpuriousLocal2020} are inspired more by control theory and describe time-varying optimization as a \emph{tracking problem} whereby an optimization algorithm defines a time-varying solution map $t \mapsto x^\star(t)$ that needs to be followed as closely as possible by the online optimization scheme.

Both viewpoints share a common base.

First, additional assumptions need to be introduced in order to guarantee well-posedness of the tracking problem. In many cases, convexity is assumed to guarantee the existence of a unique minimum and, in the case of strong convexity, a unique minimizer. Exceptions are \cite{tangRunningPrimalDualGradient2018,tang2018feedback,suboticQuantitativeSensitivityBounds2020,zavalaRealtimeNonlinearOptimization2010,dontchevEulerNewtonContinuation2013,kungurtsevPredictorCorrectorPathFollowingAlgorithm2017,suwartadiSensitivityBasedEconomicNMPC2017}, where results for non-convex optimization problems are presented.

Second, to give meaningful performance guarantees, some sort of ``bounded variation'' in the optimization problem has to be assumed.
A common assumption is that the rate of change of the optimizers is bounded by a known constant (see for example \citealt{maddenBoundsTrackingError2021,bianchinTimeVaryingOptimizationLTI2021,bernsteinOnlinePrimalDualMethods2019}). This assumption is relaxed in \cite{suboticQuantitativeSensitivityBounds2020} which shows how, in special cases, the rate of change of the optimizer can be bounded using information about the objective and the constraint functions only.

{\tb Depending on the application, the time-varying optimization problem might be unconstrained~\citep{popkovGradientMethodsNonstationary2005}, constrained to a stationary set \citep{mokhtariOnlineOptimizationDynamic2016,zinkevichOnlineConvexProgramming2003,hallOnlineConvexOptimization2015,jadbabaieOnlineOptimizationCompeting2015}, or have time-varying constraints \citep{fazlyabInteriorPointMethod2016,fazlyabPredictionCorrectionInteriorPointMethod2018,rahiliDistributedContinuousTimeConvex2017,rahiliDistributedConvexOptimization2015,suboticQuantitativeSensitivityBounds2020,tangRunningPrimalDualGradient2018,tang2018feedback,zhang2021onlineProximalADMM,bianchinTimeVaryingOptimizationLTI2021,hauswirthTimevaryingProjectedDynamical2018,bernsteinOnlinePrimalDualMethods2019}. }

Roughly speaking, algorithms for time-varying optimization can be divided into \emph{running} algorithms that do not incorporate any information about the evolution of the problem \citep{tangRunningPrimalDualGradient2018,tang2018feedback,popkovGradientMethodsNonstationary2005,bernsteinOnlinePrimalDualMethods2019,simonettoDistributedAsynchronousTimevarying2014,simonettoTimevaryingConvexOptimization2017,colombino2020online,bastianelloDistributedInexactProximal2020,bianchinTimeVaryingOptimizationLTI2021} and \emph{predictive} schemes that exploit some knowledge or estimate about the change in the optimization problem \citep{simonettoPredictionCorrectionAlgorithmsTimeVarying2017, lesage-landryPredictiveOnlineConvex2020, fazlyabPredictionCorrectionInteriorPointMethod2018, simonettoDecentralizedPredictionCorrectionMethods2017,bastianelloPrimalDualPredictionCorrection2020,zheng2020implicit}.

The following simple examples generalize \cref{ex:simp_grad} to a time-varying setting and present, just for the sake of illustration, a running and a predictive scheme.

\begin{example}\label{ex:tv_running}
    Consider the same setup as in \cref{ex:simp_grad}. Namely, let a physical plant be characterized by the steady-state input-output map $y = h(u) + d(t)$. We now assume that $d(t)$ is time-varying and (Lebesgue) measurable. Consequently, we wish to track the solution of
    \begin{align}\label{eq:ex_tv_prob}
        \underset{\tb u}{\minimize} \quad \tilde{\Phi}(u,t) := \Phi(h(u) + d(t)) \, .
    \end{align}
    Without modifying the controller we get the (non-autonomous) closed-loop system
    \begin{align}\label{eq:tv_grad}
        \dot u & = - \nabla_u \tilde{\Phi}(u,t)^T \\&=  - \nabla h(u)^T \nabla \Phi(y)^T \qquad \quad y = h(u) + d(t)
    \end{align}
    which is a running algorithm to solve \eqref{eq:ex_tv_prob}.

    Assume that $\tilde{\Phi}$ is $\beta$-strongly convex for every $t$ and consequently has a unique global minimizer $u^\star(t)$ for every $t$. Further, assume that $\| u^\star(t) - u^\star(t') \| \leq \ell \| t - t' \|$ for all $t, t'$. In other words, the unique optimizer is $\ell$-Lipschitz.

    The quantity $\ell$ can sometimes be bounded from problem parameters \citep{suboticQuantitativeSensitivityBounds2020}. For instance, if it is known that $d$ is $\ell_d$-Lipschitz in $t$ and $\Phi$ is $\beta'$-strongly convex, then, the estimate $\ell \leq \beta' \ell_d$ holds.

    Exploiting strong convexity and Cauchy-Schwarz inequality, the distance between $u(t)$ of \eqref{eq:tv_grad} and $u^\star(t)$ can be shown to differentiable for almost all $t$ and satisfy
    \begin{align*}
         & \tfrac{d}{dt} \tfrac{1}{2} \| u(t) - u^\star(t) \|^2                                               \\
         & \qquad \leq \left\langle \dot u(t), u(t) - u^\star(t) \right\rangle + \ell \| u(t) - u^\star(t) \|
        \\
         & \qquad \leq - \beta \| u(t) - u^\star(t) \|^{2} + \ell \| u(t) - u^\star(t) \|
        \, .
    \end{align*}
    Consequently, $\| u(t) - u^\star(t) \|$ is decreasing as long as $\| u(t) - u^\star(t) \| > \ell / \beta$. It follows from standard invariance arguments that, as $t \rightarrow \infty$, $u(t)$ will be $\ell/\beta$-close to $u^\star(t)$. 
\end{example}

\begin{example}
    To illustrate a continuous-time predictive algorithm, consider the same setup as in \cref{ex:tv_running} and assume, in addition, that $\Phi$ is twice continuously differentiable, and the time derivative $\dot{d}(t)$ is available (e.g., can be estimated using finite differences). Then, we may consider the following control law derived from the sensitivity-conditioning approach \citep{picallo2023sensitivity}
    \begin{align*}
        \dot  u & = - \nabla h(u)^T \nabla_u \Phi(y)^T  - \nabla h(u)^T \nabla_{yy}^2 \Phi(y) \dot{d}(t) \\
        y       & = h(u) + d(t)
    \end{align*}
    which reduces to
    \begin{align*}
        \dot  u = - \nabla_u \tilde{\Phi}(u,t)^T  - \tfrac{d}{dt} (\nabla_u \tilde{\Phi})^T(u,t) \, .
    \end{align*}
    One may consider the time-varying Lyapunov function $W(u,t) := \tfrac{1}{2} \| \nabla_u \tilde{\Phi}(u,t) \|^2$. It holds that
    \begin{align*}
        \tfrac{d}{dt}{W}(u(t),t) & = \nabla_u W(u(t),t) \dot u(t) + \nabla_t W(u(t),t) \\
                                 & = \| \nabla_u \tilde{\Phi}(u(t),t) \|^2 < 0 \, .
    \end{align*}
    for all $t$ and all $u(t) \neq u^\star(t)$. This fact can be used to guarantee zero tracking error as $t \rightarrow \infty$.
\end{example}

\section{Online Feedback-Based Optimization}\label{sec:ofo}

We now present the key idea and main topic of this {\tb review} article---the implementation of optimization algorithms in closed loop with physical systems. The unconstrained feedback gradient flow from \cref{ex:simp_grad} in the introduction {\tb can be considered as a minimal viable example} for this approach.
A more general formulation is the problem of driving an asymptotically stable plant with the steady-state map $y=h(u)+d$ to the solution of the optimization problem
    \begin{subequations}\label{eq:basic_prob}
    \begin{align}
        \min_{u,y} & \quad \Phi(y) \\
        \subjto & \quad y = h(u) + d\\
        & \quad u \in \mathcal U \\
        & \quad y \in \mathcal Y,
    \end{align}
    \end{subequations}
where $\mathcal U$ and $\mathcal Y$ are the feasible regions of the input $u$ and the output $y$, respectively.
This formulation can be generalized to include a dependence of the cost $\Phi$ on the input $u$ directly, without any technical complication. 

{\tb Almost any optimization algorithm in the form of a dynamical system can be defined as open systems and interconnected with the physical plant, but not every such control design makes sense in the real world and achieves the aforementioned goal. The resulting interconnection needs to be robust, satisfy the plant operational constraints, and be implementable given the limited model knowledge.

The seminal papers on network congestion control \citep{lowOptimizationFlowControl1999,kellyRateControlCommunication1998} have demonstrated the potential of studying optimization dynamics in closed-loop implementations, see \cref{ex:netw_cong} for a self-contained review. Afterward, \cite{jokicConstrainedSteadystateRegulation2009} and \cite{brunnerFeedbackDesignMultiagent2012} consider a more abstract formulation of steering the outputs of a nonlinear continuous-time dynamical plant to a steady-state that solves a constrained convex optimization problem. These pioneering works inspired a plethora of contributions in different directions.

In this section, we review three fundamental aspects that constitute the core challenges when aiming to achieve a \emph{robust} interconnection of an optimization algorithm with a physical plant:
\begin{itemize}
    \item In \cref{sec:cl_stab}, we present approaches to guarantee \emph{closed-loop stability} in the presence of plant dynamics, either via singular perturbation analysis or linear matrix inequality certificates.
    \item In \cref{sec:cstr_enforc}, we show how partial model information and feedback measurements can be used to \emph{robustly enforce constraints} in closed-loop optimization.
    \item In \cref{subsec:extension}, we discuss strategies for bypassing the need of model information and achieving fully \emph{data-driven and model-free operations}.
\end{itemize}
For each of these themes, we provide a tutorial exposition and review the relevant literature.
}

\subsection{Closed-Loop Stability}\label{sec:cl_stab}

\Cref{ex:num_grad} in the introduction has illustrated that a simple gradient-based controller interconnected with a dynamical system is not necessarily stable, unless the control gain $\epsilon$ is small enough (\cref{fig:num_grad}). In other words, \emph{sufficient timescale separation} between the fast plant behavior and the slow optimization dynamics is generally required. 

For control design, it is convenient to reduce the fast transient plant dynamics to their algebraic steady-state map. We have explicitly applied this design step in \cref{ex:simp_grad} and will do so implicitly in the rest of the article. In this section, however, we will investigate the problems of timescale separation and closed-loop stability in detail, and survey works that have tackled this issue.

We identify and present two main research streams: First, we discuss stability results inspired by singular perturbation analysis, which formalize and quantify stability in terms of timescale separation. This approach is very general and applicable to nonlinear (but asymptotically stable) plant dynamics and non-convex optimization dynamics, but potentially very conservative. A second line of research taps into robust control and provides computational stability certificates in the form of linear matrix inequalities. While these stability guarantees are less conservative, they apply only to linear time-invariant plants and convex optimization dynamics.

\begin{remark}
    Exploiting \emph{passivity} is sometimes a third possibility to certify closed-loop stability \cite[Chap. 6]{khalilNonlinearSystems2002}. This approach requires that the control loop can be identified as a feedback interconnection of passive systems, which is not always possible.
    Saddle-point flows (discussed in \cref{sec:saddle}) are one class of systems that are amenable to passivity arguments to study stability and robustness \citep{vanderschaftRelationPortHamiltonianGradient2011,simpson-porcoInputOutputAnalysis2016}. These ideas have been applied to congestion avoidance in communication networks (\cref{ex:netw_cong}) in \cite{wenUnifyingPassivityFramework2004} and in Chapter 4 of \cite{lowAnalyticalMethodsNetwork2017}. Similarly, passivity has been exploited in power systems applications (like \cref{ex:opt_freq}) in \cite{steginkUnifyingEnergyBasedApproach2017,tripPassivityBasedDesign2017}.
\end{remark}

\begin{remark}
    In the following,  for the sake of keeping the presentation consistent, we maintain the same modeling framework where both the plant and the optimization algorithm are described as continuous time flows. 
    In some applications, for example when the presence of communication delays requires less frequent algorithms iterations, a sampled-data approach may be more appropriate. 
    \cite{belgioioso2021sampled, belgioioso2022online} provide stability results in this hybrid setting.
\end{remark}

\subsubsection{Singular Perturbation Analysis}\label{subsec:singular_perturb}

For the purpose of providing explicit bounds on $\epsilon$, \cite{mentaStabilityDynamicFeedback2018,hauswirthTimescaleSeparationAutonomous2020} pursue a singular perturbation approach and arrive at easy-to-compute sufficient conditions for closed-loop stability. In particular, \cite{hauswirthTimescaleSeparationAutonomous2020} considers general nonlinear (but stable) plant dynamics and investigates a variety of optimization dynamics for convex and non-convex problems including (projected) gradient, saddle-point, and momentum methods.
Using similar techniques, \cite{simpson-porcoAnalysisSynthesisLowGain2020} provides stability guarantees for general low-gain integral controllers that satisfy an infinitesimal contraction property, but with a discussion of feedback-based optimization limited to LTI plant dynamics.
The following example illustrates the main idea behind these stability conditions.

\begin{example}\label{ex:ts_arg}
    We wish to characterize the stability of the feedback loop introduced in \cref{ex:simp_grad}. In particular, we want to formulate conditions on the gain $\epsilon$ in \cref{fig:simp_grad} that guarantee closed-loop stability. For this purpose we pass to the singular perturbation decomposition into reduced and boundary-layer error dynamics illustrated in \cref{fig:simp_grad_sp} \citep{khalilNonlinearSystems2002,kokotovicSingularPerturbationMethods1999}.
    In particular, $\hat{h}$ is defined such that $f(\hat{h}(u), u) = 0$ for all $u$.

    \begin{figure}[bt]
        \centering
        \includegraphics[width=\columnwidth]{Figs/sing_pert_fig.tex}

        \caption{Feedback-based gradient flow from \cref{fig:simp_grad} rearranged and decomposed into reduced and boundary-layer error dynamics}\label{fig:simp_grad_sp}
    \end{figure}

    The resulting reduced dynamics correspond exactly to the simplified model that we have already used in the design of the optimizing controller, where the plant is replaced by its algebraic steady-state map.
    The boundary-layer error dynamics $z := \zeta - \hat{h}(u)$ evolve as $\dot z = f(\zeta, u)$ for any fixed $u$. If these error dynamics are exponentially stable (and other technical assumptions are satisfied), standard converse results guarantee the existence of a Lyapunov function $W$ and parameters $\gamma, \omega > 0$ such that, for any fixed $u$, it holds that
    \begin{align}\label{eq:ts_lyap_params}
        \begin{split}
            \dot W(\zeta - \hat{h}(u)) & \leq - \gamma \| \zeta - \hat{h}(u) \|^2 \\ \| \nabla_u  W(\zeta - \hat{h}(u)) \|  & \leq \omega \| \zeta - \hat{h}(u) \| \, .
        \end{split}
    \end{align}

    A class of Lyapunov function candidates to certify stability of the closed system in \cref{fig:simp_grad_sp} is given by
    \begin{align*}
        V_\delta (u, \zeta) = \delta\tilde{\Phi}(u) + (1 - \delta) W(\zeta - \hat{h}(u))
    \end{align*}
    where $\delta \in (0,1)$ is a convex combination parameter.

    Let $L$ be the Lipschitz constant of $\nabla \hat{\Phi}$. \cite{hauswirthTimescaleSeparationAutonomous2020} have shown that for all
    \begin{align}\label{eq:ts_threshold}
        \epsilon < \epsilon^\star := \frac{\gamma}{\omega L}
    \end{align} the parameter $\delta$ can be chosen such that $V_\delta$ is non-increasing and thus a LaSalle invariance argument guarantees (asymptotic) stability.
\end{example}

\begin{remark}
    Recall from Section~\ref{sec:cstr_sum} that constraints can be incorporated into an optimization problem through the combination of gradient flows with penalty or barrier functions. The bound \eqref{eq:ts_threshold} on $\epsilon$ in \cref{ex:ts_arg} indicates that there is a natural limitation to this constraint enforcement mechanism. Notably, penalty functions cannot be made arbitrarily steep, thus increasing $L$, without at the same time reducing $\epsilon$ by the same fraction. For the same reason log-barrier functions that do not have a Lipschitz gradient have to be applied with great care, as closed-loop stability cannot be guaranteed with these singular perturbation bounds on $\epsilon$.
\end{remark}

The timescale separation argument inspired by singular perturbation analysis, as presented in \cref{ex:ts_arg}, works under very general conditions. In particular, for the case of a gradient flow convexity of $\Phi$ is not generally required. The type of stability proof can also be established for various optimization algorithms interconnected with exponentially stable plants. Examples include but are not limited to Newton flows, projected gradient flows, and saddle-point flows encountered in the previous section. {\tb Moreover, the above timescale separation argument can also be applied in the sampled-data setting that involves the interconnection of continuous-time plants and discrete-time equilibrium-seeking algorithms, see the feedback equilibrium seeking in \citep{belgioioso2021sampled,belgioioso2022online,belgioioso2023tutorial}.}
These conditions can be very conservative, but they are qualitatively tight. Namely, non-examples in \cite{hauswirthTimescaleSeparationAutonomous2020} show how subgradient flows and continuous-time accelerated gradient flows interconnected with dynamical plants are not generally stable. These setups are not amenable to the same type of timescale separation argument as above because important assumptions such as uniform asymptotic stability of the reduced dynamics are not satisfied \citep{povedaInducingUniformAsymptotic2019}.

\subsubsection{LMI Stability Certificates}
Alternatively, closed-loop stability can be certified by applying tools from robust control, namely linear matrix inequalities. Recent results on integral quadratic constraints and their use for the analysis of optimization algorithms \citep{lessardAnalysisDesignOptimization2016,fazlyabAnalysisOptimizationAlgorithms2018} have proven very useful for this purpose and have been applied in \cite{colombino2020online,nelsonIntegralQuadraticConstraint2018,colombino2019towards}. In particular, \cite{nelsonIntegralQuadraticConstraint2018} studies the joint design of stabilizing control, estimator, and optimization dynamics.
Similarly, \cite{lawrenceOptimalSteadystateControl2018,lawrenceLinearConvexOptimalSteadyState2021} consider an output regulation framework and reduce the control design to a stabilization problem for the case of convex problems and LTI systems.

Limiting the applicability of these techniques from robust control is the fact that plant dynamics are generally required to be LTI, and objective functions need to be (strongly) convex. Moreover, these LMI-based conditions are in the form of computational stability certificates and do not directly translate into control design procedures or tuning recommendations.

The following example presents a LMI-IQC based stability test for the feedback-based gradient flow.
\begin{example}\label{ex:lmi}
    Consider the same setup as in \cref{ex:simp_grad}, and, in addition, assume that $\Phi$ is $m$-strongly convex and has a $L$-Lipschitz gradient. Further, assume that
    \begin{align}\label{eq:lti_sys}
        \dot \zeta = f(\zeta, u) = A x + B u \qquad y = g(\zeta) + d = C \zeta + d
    \end{align}
    is an LTI system. Assuming that $A$ is invertible, the steady-state input-output map is given by $h(u) = H u $ with $H = - C A^{-1} B$. Furthermore, \eqref{eq:ic_grad_system} has a unique equilibrium $(\zeta^\star, u^\star)$ such that $\nabla \tilde{\Phi}(u^\star) = 0$ and $0 = A \zeta^\star + B u^\star$.

    \begin{figure}
        \centering
        \begin{tikzpicture}
            \matrix[ampersand replacement=\&, row sep=0.45cm, column sep=.8cm] {
                \& \node[block] (lti) {
                    \begin{tabular}{rl}
                        $\dot{\bfz}$ & $= \bfA \bfz + \bfB \bfu$ \\
                        $\bfy$       & $= \bfC \bfz + \bfD w$
                    \end{tabular}}; \\
                \node[none] (edge1) {}; \&
                \node[block, inner sep=.35cm] (nonlin) {$\Delta_\epsilon$}; \&
                \node[none] (edge2) {};    \\
            };

            \draw[line] (lti.east)-|(edge2.center)  node[near end, right]{$\bfy$};
            \draw[connector] (edge2.center)--(nonlin.east);
            \draw[line] (nonlin.west)--(edge1.center);
            \draw[connector] (edge1.center)|-([yshift=-.3cm]lti.west) node[near start, left]{$\bfu$};
            \draw[connector] ([yshift=.3cm,xshift=-1cm]lti.west)--([yshift=.3cm]lti.west) node[at start, left]{$w$};
        \end{tikzpicture}
        \caption{Feedback-based gradient flow from \cref{fig:simp_grad} rewritten as a feedback interconnection of an LTI system with a nonlinearity $\Delta$}\label{fig:robu_ic}
    \end{figure}

    The closed-loop system \cref{eq:ic_grad_system} can be written in a standard robust control setup as an interconnection of a nonlinearity with an LTI system. This is shown in \cref{fig:robu_ic}, where
    \begin{align*}
        \bfz & := \begin{bmatrix} \zeta \\ u \end{bmatrix} \quad & \bfA & := \begin{bmatrix} A & B \\ 0 & 0 \end{bmatrix} & \quad \bfB & := \begin{bmatrix} 0 & 0 \\ 0 & - H^T \end{bmatrix} \\
        \bfy & := \begin{bmatrix} y \\ u\end{bmatrix} \quad & \bfC & := \begin{bmatrix} 0 & \bbI \\ C & 0 \end{bmatrix} & \quad \bfD & := \begin{bmatrix} \bbI \\ 0 \end{bmatrix}
    \end{align*}
    and
    \begin{align}
        \bfu = \Delta_\epsilon(\bfy) =  \begin{bmatrix} 0 \\ \epsilon \nabla \Phi(y) \end{bmatrix}  \, .
    \end{align}
    Crucially, $\Delta_\epsilon(\cdot)$ satisfies the IQC defined by
    \begin{align}
        \begin{bmatrix}
            \bfy - \bfy^\star \\ \bfu - \bfu^\star
        \end{bmatrix}^T
        \underbrace{\begin{bmatrix}
                -2 \epsilon^2 mL \bbI & \epsilon (L + m) \bbI \\  \epsilon (L + m) \bbI & - 2 \bbI
            \end{bmatrix}}_{\Xi_\epsilon}
        \begin{bmatrix}
            \bfy - \bfy^\star \\ \bfu - \bfu^\star
        \end{bmatrix} \geq 0
    \end{align}
    for all $(y,u)$ such that $y = \Delta(y)$ \citep[Lemma 6]{lessardAnalysisDesignOptimization2016}. Equivalently, we have
    \begin{align}\label{eq:iqc}
        \begin{bmatrix}
            \bfz - \bfz^\star \\ \bfu - \bfu^\star
        \end{bmatrix}^T \bfC^T {\Xi_\epsilon} \bfC
        \begin{bmatrix}
            \bfz - \bfz^\star \\ \bfu - \bfu^\star
        \end{bmatrix} \geq 0
    \end{align}
    for all $(\bfz, \bfu)$ such that $\bfu = \Delta(\bfC \bfz)$.

    Next, consider a Lyapunov function of the form $ V(\bfz) = (\bfz - \bfz^\star)^T \bfP (\bfz - \bfz^\star)$ where $\bfP \succeq 0$ remains to be determined. The time derivative of $V$ along system trajectories can be written as a quadratic form
    \begin{align}\label{eq:time_deriv}
        \dot V(\bfz) =  \begin{bmatrix}
            \bfz - \bfz^\star \\ \bfu - \bfu^\star
        \end{bmatrix}^T
        \begin{bmatrix} \bfA^T \bfP + \mathbf{P A} & \mathbf{P B} \\ \bfB^T \bfP & 0 \end{bmatrix}
        \begin{bmatrix}
            \bfz - \bfz^\star \\ \bfu - \bfu^\star
        \end{bmatrix} \, .
    \end{align}
    Combining \eqref{eq:iqc} and \cref{eq:time_deriv}, it follows that $\dot V(\bfz) < 0 $ for $\bfz \ne \bfz^\star$ (and thus the closed-loop is stable) if the LMI
    \begin{align}
        \begin{bmatrix} \bfA^T \bfP + \bfP \bfA & \mathbf{P B} \\ \bfB^T \bfP & 0 \end{bmatrix} \prec -  \bfC^T {\Xi_\epsilon} \bfC
    \end{align}
    holds for some $\bfP \succeq 0$.
\end{example}

\subsection{Constraint Enforcement in Closed Loop}\label{sec:cstr_enforc}

The capacity and various possibilities to robustly enforce complicated constraints despite model uncertainty is one of the distinguishing features of feedback-based optimization.

In \cref{sec:opt_dyn} we have seen different mechanisms used in continuous-time optimization algorithms to deal with constraints. \Cref{sec:cstr_sum} summarizes these options and highlights how they can be combined. In this subsection, we will explore the closed-loop implementations of these mechanisms. In particular, in Section~\ref{sec:sat_proj} we will discuss the use of saturation to enforce constraints on plant inputs without major computational effort. In Section~\ref{sec:engconstraints}, we will review strategies to handle general output constraints. 

\subsubsection{Input Saturation via Projection and Anti-Windup}\label{sec:sat_proj}

Physical plants are generally subject to limited actuator capabilities, simply due to the fact that any realistic system is bound by the laws of physics and can handle only signals of finite power.
Such actuation limits can often be modeled as input saturation of control signal. At a physical level, such saturation is, for example, the result of mechanical constraints, or it can be observed in electric circuits involving diodes and other semi-conducting devices.
From a systems perspective, low-level controllers that protect devices and subsystems from operating outside of a safe zone of conditions can also be modeled as saturation.
The nature of these constraints implies that they cannot be violated at any point in time. Given the constraint enforcement possibilities from \cref{sec:cstr_sum}, it makes sense to model saturation as a projection onto the set of feasible inputs with the key property that the projection is naturally ``evaluated'' and applied by the physical system.

This property creates the possibility of outsourcing the handling of these constraints from the controller to the plant and thus reduces computational requirements and the need for exact modeling of these constraints, which can be time-varying and/or unknown.

For discrete-time feedback-based optimization schemes such as those in \cite{bernsteinOnlinePrimalDualMethods2019,dallaneseOptimalPowerFlow2018,ganOnlineGradientAlgorithm2016,simpson-porcoLowGainStabilityProjected2022}, this model-free constraint handling by saturation is relatively straightforward (at least for gradient descent with standard Euclidean metric), under the assumption that the saturated control signal is measured and available within one sampling interval. In this case, it is possible to simply add the control increment to the measured saturated control signal of the preceding sampling interval.

For continuous-time methods, as in \cite{hauswirthProjectedGradientDescent2016,changSaddleflowDynamicsDistributed2019}, exploiting input saturation is trickier because a continuous-time integrator in cascade with a saturation element will generally lead to \emph{integrator windup}. For this reason \cite{hauswirthRobustImplementationProjected2020,hauswirthAntiWindupApproximationsOblique2020a,hauswirthDifferentiabilityProjectedTrajectories2020} study the use of anti-windup compensators for feedback-based optimization and rigorously show how these control designs can be used to smoothly and robustly approximate continuous-time projected gradient flows like \cref{eq:proj_grad_flow}. The following example illustrates this idea.

\begin{example}\label{ex:awa}

    Consider the same setup as in \cref{ex:simp_grad}. Namely, we want to drive a plant (with fast decaying dynamics and steady-state map $y = h(u) + d$) to an optimal steady state minimizing the cost function $\tilde{\Phi}(u)$.

    In addition, we assume that the plant is subject to input saturation that acts like a projection onto a non-empty set $\calU$ of admissible inputs. For simplicity, assume that $\calU$ is convex. In order to mitigate the effects of integrator windup, a simple anti-windup scheme with a tuneable gain $\tfrac{1}{K}$ is in place as illustrated in \cref{fig:simp_aw}.

    \begin{figure}[!t]
        \centering
        \includegraphics[width=\columnwidth]{Figs/awa_grad.tex}
        \caption{Anti-windup system \eqref{eq:ic_awa_sys} which approximates a projected gradient flow \eqref{eq:proj_grad_flow_awa} as $K \rightarrow 0^+$}\label{fig:simp_aw}
    \end{figure}

    A state-space representation of this system is given, analogously to \cref{eq:ic_grad_system}, by
    \begin{align}\label{eq:ic_awa_sys}
        \begin{split}
            \text{plant} & \begin{cases}
                \overline{u} & = P_\calU(u)                 \\
                \dot \zeta   & =   {f}(\zeta, \overline{u}) \\
                y            & = g(\zeta) + d
            \end{cases} \\
            \text{controller} & \begin{cases} \dot u & = - \epsilon \nabla h(u)^T \nabla \Phi(y)^T - \tfrac{1}{K}(u - \overline{u})\, .
            \end{cases}
        \end{split}
    \end{align}

    Assuming, as before, that the plant dynamics can be approximated by the algebraic map $y = h(u) + d$ and with $\tilde{\Phi}(u) = \Phi(h(u) + d)$, the system \eqref{eq:ic_awa_sys} reduces to
    \begin{align}\label{eq:aw_grad_flow}
        \dot u = - \epsilon \nabla \tilde{\Phi}(P_\calU(u))^T - \tfrac{1}{K} ( u - P_\calU(u)) \, .
    \end{align}

    The key aspect in \cref{eq:aw_grad_flow} is the fact that $\nabla \tilde{\Phi}$ is evaluated at $P_\calU(u)$ rather than at $u$. In \cite{hauswirthDifferentiabilityProjectedTrajectories2020} it was shown that trajectories $u(t)$ of \eqref{eq:aw_grad_flow} converge in the sense that $P_\calU(u(t))$ converges to the set of KKT points of the optimization problem
    \begin{align}\label{eq:awa_opt_prob}
        \minimize \quad  \tilde{\Phi}(u) \qquad \subjto \quad  u \in \calU \, .
    \end{align}
    In particular, if \eqref{eq:awa_opt_prob} is convex, then convergence of $P_\calU(u(t))$ is to the global minimizer.

    Importantly, this means that upon convergence, the plant in \cref{fig:simp_aw} is at an optimal steady state even though the (unsaturated) control $u$ is not optimal (only the saturated control $P_\calU(u)$ is optimal).

    Furthermore, \cite{hauswirthRobustImplementationProjected2020,hauswirthAntiWindupApproximationsOblique2020a} have shown that \eqref{eq:aw_grad_flow} approximates the projected gradient flow
    \begin{align}\label{eq:proj_grad_flow_awa}
        \dot{u} = \Pi_\calU\left[- \nabla h (\overline{u})^T \nabla \Phi(y)^T \right](\overline{u})
    \end{align}
    as $K \rightarrow 0^+$. This insight extends to other anti-windup schemes which can be shown to approximate more general projected dynamical systems.
\end{example}

\subsubsection{Output Constraints via Dualization and Approximate Projections}
\label{sec:engconstraints}

We now turn to mechanisms that allow us to enforce more general constraints, and in particular those that apply to plant outputs. 

In the following, we focus on two approaches. First, we illustrate how augmented and projected saddle-point flows, as in \cref{ex:sadd_flow}, can be implemented as feedback controllers. This approach is particularly suited for solving convex problems in closed loop and distributed over a network (assuming the problem exhibits a suitable sparsity structure). For non-convex problems, these algorithms can still be applied, but theoretical global convergence guarantees are not generally available.

As a second possibility to enforce output constraints, we discuss a special discretization of projected gradient flows that can be implemented as a feedback controller and comes with strong global convergence guarantees, even for non-convex setups. However, it is less easily amenable to a distributed implementation and, instead, requires the solution of a simple quadratic program at every iteration.

\paragraph{Projected Saddle-Flows as Feedback Controllers}
Consider the optimal steady-state problem \cref{eq:basic_prob}. Assume that the disturbance $d$ is fixed and that the engineering constraints $(u,y) \in \calX$ can be expressed as $\calX := \{ (u,y) \, | \, g(y,u) \leq 0 \}$, where $g: \bbR^n \rightarrow \bbR^m$ denotes a continuously differentiable constraint function.
After eliminating $y$ from \cref{eq:basic_prob}, we are left with the reduced problem
\begin{subequations}\label{eq:basic_prob_red}
    \begin{align}
        \underset{\tb u}{\minimize} \quad & \hat{\Phi}(u) := \Phi(u, h(u)+d)                         \\
        \subjto \quad   & u \in \calU                                              \\
                  & u \in \hat{\calX} := \{ u \, | \, \hat{g}(u) \leq 0 \}\,
    \end{align}
\end{subequations}
with $\hat{g}(u) := g(u, h(u))$.

This problem can be tackled with a projected saddle-point flow with a primal augmentation term $\phi(u) = \tfrac{\rho}{2} \| \max\{ \hat{g}(u) , 0 \} \|^2$ as presented in \cref{ex:sadd_flow}. Namely, we have
\begin{subequations}\label{eq:proj_pd_cl_control}
    \begin{align}
        \dot{u}   & = \Pi_\calU \left[ - \nabla \hat{\Phi}(u)^T - \nabla \hat{g}(u)^T (\mu + \rho \max\{\hat{g}(u), 0 \})\right](u) \\
        \dot{\mu} & = \Pi_{\bbR^m_{\geq 0}} \big[ \hat{g}(u) \big](\mu) \, .
    \end{align}
\end{subequations}
Recall from \cref{ex:sadd_flow} that, in the special case where \eqref{eq:basic_prob_red} is convex, trajectories of \eqref{eq:proj_pd_cl_control} converge to a global minimizer. In particular thanks to the primal augmentation, as in \cref{eq:sadd_primal_aug}. The primal augmentation term reduces oscillations, especially with non-strictly convex objective functions. If the problem is non-convex, a dual augmentation term, as in \cref{eq:sadd_dual_aug}, can be added to guarantee global convergence at the expense of changing equilibrium points away from the true KKT points of \eqref{eq:basic_prob_red}.

We wish to realize \eqref{eq:proj_pd_cl_control} as the closed-loop dynamics of a feedback loop incorporating the physical plant. To achieve this, we replace any evaluation of $h(u) + d$ with the plant output $y$. This yields the controller
\begin{subequations} \label{eq:pd_control}
    \begin{align}
        \dot{u}   & = \Pi^Q_\calU \left[ - \epsilon Q(u) H(u) G(u,y) \right](u) \\
        \dot{\mu} & = \Pi_{\bbR^m_{\geq 0}} \big[ g(u,y) \big](\mu)
    \end{align}
\end{subequations}
where $H(u) := \begin{bmatrix} \bbI & \nabla h(u)^T \end{bmatrix}$ and
\begin{align*}
    G(u,y) := \begin{smallbmatrix}
        \nabla_u \Phi(u,y)^T - \nabla_u g(u,y)^T \left(\mu + \rho \max\{g(u,y), 0 \} \right) \\ \nabla_y \Phi(u,y)^T  - \nabla_y g(u,y)^T \left(\mu + \rho \max\{g(u,y), 0 \} \right) \end{smallbmatrix} \, .
\end{align*}
We have also included in \cref{eq:pd_control} a variable metric $Q(u)$ on the primal variables. This degree of freedom will be exploited in the forthcoming \cref{ex:pw_ex_pd,ex:lop_survey}.

In practice, one can build an estimate $\hat{H}$ of $H(u)$ by using $y$ and $u$, i.e., $\hat{H}(u,y) \approx H(u)$. We defer the discussion on how to perform this estimation to Section~\ref{subsubsec:sens_learn}.

Furthermore, if the plant \eqref{eq:simple_plant} is an asymptotically stable LTI system subject to constant disturbances
\begin{align}
    \dot{\xi} = A \xi + B u + d_1 \qquad y = C \xi + D u + d_2 \, ,
\end{align}
then the steady-state sensitivity matrix $\nabla h(u) = - C A^{-1} B + D$ is constant and independent of $d$, which renders its estimation much easier.

Despite a lack of theoretical convergence guarantees for non-convex problems and closed-loop stability certificates, primal-dual saddle-point algorithms (and their projected, proximal, and augmented variations) have enjoyed great popularity. In particular, saddle-point flows often lend themselves to a distributed implementation. \Cref{ex:netw_cong,ex:opt_freq} are prime examples. In the context of \eqref{eq:pd_control}, a distributed implementation is possible if $\nabla h(u)$ has a fixed sparsity pattern, the cost function is separable, and the constraints $u \in \calU$ and $g(u,y) \leq 0$ are appropriately localized. Moreover, selecting a sparsity-inducing metric $Q(u)$ can also facilitate distributed operations. In this case, the computation of a given component $u_i$ requires only the quantity $\mu_j + \rho \max\{g_j(u,y), 0 \}$ to be communicated from all neighboring nodes $j$ in a network.

\paragraph{Linearized Output Projections}%

As a second possibility to solve \eqref{eq:basic_prob} in closed loop (and thus tackle the reserve dispatch problem from \cref{sec:redispatch}), we highlight the method proposed in \cite{haberleNonconvexFeedbackOptimization2021} and extended in \cite{hauswirthOptimizationAlgorithmsFeedback2020b}. This method is conceptually related to \cite{torrisiProjectedGradientConstraint2018}, in that we seek a descent direction that ensures satisfying the linearized constraint.

We have seen how projected (and augmented) saddle-point flows can be easily turned into feedback controllers. The key mechanism for enforcing the engineering constraints $(u,y) \in \calX$ thereby lies in integrating the constraint violation $g(u,y)$ that is based on plant output measurements. In \cref{ex:pw_ex_pd} we have noted the ease of deriving a discrete controller by simply applying an explicit Euler discretization.

Now, instead, we present a discrete controller that approximates the projected gradient flow
\begin{align}\label{eq:proj_grad_flow_comp}
    \dot{u} = \Pi_{\calU \cap \hat{\calX}} [ - \nabla \hat{\Phi}(u)](u)\, .
\end{align}
The trajectories of \eqref{eq:proj_grad_flow_comp} are guaranteed to converge to the KKT points of the reduced problem \cref{eq:basic_prob_red}, even under non-convexity (but under technical assumptions similar to those in \cref{thm:grad}). We refer to \cite{hauswirthProjectedDynamicalSystems2020} for a discussion of the available convergence guarantees for projected gradient flows.

We have previously seen in \cref{sec:sat_proj} that input constraints $u \in \calU$ are often enforced by input saturation (possibly requiring an anti-windup compensator to avoid integrator windup). The challenge in \eqref{eq:proj_grad_flow_comp} lies in the fact that also the engineering constraints $\calX$ need to be enforced ``by projection'' and in feedback using measurements of the plant output, rather than using a model-based evaluation of $h(u) + d$.

To address this problem, rather than pondering on the continuous-time dynamics \eqref{eq:proj_grad_flow_comp}, we directly consider a discrete-time controller of the form
\begin{equation} \label{eq:lop_up_survey}
    u^+ = u + \alpha\, \Sigma_\alpha(u,y) \,, %
\end{equation}
where $\alpha>0$ is a fixed step-size, $y = h(u) + d$ is the measured system output, and $\Sigma_\alpha(u,y)$ is defined as the solution of
\begin{subequations}\label{eq:sigma_hat}
    \begin{align}
        \underset{w\in\bbR^p}{\minimize} \, & \| w + Q(u) H(u) \nabla \Phi (u,y)^T  \|^2_{Q(u)}
        \\
        \text{subject to} \quad
            & u+ \alpha w \leq \calU \label{eq:lop_in_cstr}                                     \\
            & g(u + \alpha w,  y+ \alpha \nabla h(u) w) \leq 0 \,, \label{eq:lop_out_cstr}
    \end{align}
\end{subequations}
where, as before, $H(u):=\begin{bmatrix}\mathbb{I}_p & \nabla h(u)^T\end{bmatrix}$ and $Q(u)$ is a metric that assigns to every $u \in \calU$ a positive definite matrix.

The feedback law \eqref{eq:lop_up_survey} approximates a projected gradient descent, by computing a descent direction $w$ that is feasible with respect to $\calU$ (see \cref{eq:lop_in_cstr}) and approximately feasible (up to first order) with respect to $\calX$ (see \cref{eq:lop_out_cstr}. In particular, it can be rigorously shown that \eqref{eq:lop_up_survey} approximates \eqref{eq:proj_grad_flow_comp} as $\alpha \rightarrow 0^+$ \citep{hauswirthOptimizationAlgorithmsFeedback2020b}.

Any equilibrium point of \eqref{eq:lop_up_survey} is feasible and a KKT point of \eqref{eq:basic_prob}. Global convergence of \eqref{eq:lop_up_survey} to KKT points of \eqref{eq:basic_prob} is guaranteed under weak technical assumptions and without convexity, and for a small enough step size $\alpha$ \citep{haberleNonconvexFeedbackOptimization2021}.

Note that, computing $\Sigma_\alpha(u,y)$ requires the solution of a quadratic program. However, in comparison to \ac{rti} schemes discussed in \cref{sec:rti}, the computational effort does not scale with a prediction horizon and no explicit model of the plant dynamics is required. Instead, it is enough to estimate $\nabla h(u)$. On the other hand, in general, \eqref{eq:lop_up_survey} does not lend itself to a natural distributed implementation. Although, depending on the problem structure (e.g., a diagonally dominant $\nabla h(u)$ or a sparsity-inducing metric $Q(u)$), one can solve \eqref{eq:sigma_hat} distributedly at every iteration.

\subsection{{\tb Data-Driven and Model-Free Operation}}\label{subsec:extension}
{\tb
The appealing benefits of the aforementioned feedback-based optimization schemes are built on limited model information, i.e., the sensitivity $\nabla h(u)$ of the steady-state input-output map of a physical plant. Due to the chain rule, this information is required during the gradient-based (i.e., first-order) updates of those schemes. In some cases, we can derive a quite accurate sensitivity estimate based on first-principle models and parameters \citep{ortmann2020experimental,ortmann2020fully}. In a broader context, however, the sensitivity may be elusive because of the complex structure of a plant and the lack of key model parameters. As a result, the effective functioning (e.g., stability and optimality) of feedback-based optimization is affected. %

\begin{example}
    Consider a similar setup as in \cref{ex:lmi}. The objective $\Phi(u,y) = \frac{1}{2} (u^{\top}Ru + y^{\top}Qy)$ is strongly convex and has Lipschitz continuous gradients with respect to $u$ and $y$. A discrete-time linear plant is given by
    \begin{equation}\label{eq:discrete_linear_plant}
        \zeta_{k+1} = A\zeta_k + Bu_k, \qquad y_k = C\zeta_k + d.
    \end{equation}
    The steady-state input-output map is $y = Hu+d$, where $H = C(I-A)^{-1}B$ is the sensitivity matrix. Let $\hat{H}$ be an estimation of $H$. We synthetically generate $\hat{H}$ by $\hat{H} = H + \sigma M$, where $\sigma \geq 0$, and the elements of $M$ are drawn from the uniform distribution. \cref{fig:inexact_sens_performance} illustrates the performance of the closed-loop interconnection of the plant \eqref{eq:discrete_linear_plant} and the following feedback-based optimization scheme using $\hat{H}$
    \begin{equation*}
        u_{k+1} = u_k - \eta(Ru_k + \hat{H}^{\top}Qy_{k+1}).
    \end{equation*}
    We use the same step size $\eta$ for different cases of $\sigma$. We observe that feedback-based optimization is robust to small inaccuracies \citep{colombino2019towards}. However, as $\sigma$ further increases (i.e., the error level of the sensitivity estimation increases), the closed-loop behavior exhibits sub-optimality and even instability. Hence, the accuracy of the sensitivity plays an important role in guaranteeing effective closed-loop performance.
    \begin{figure}
        \centering
        \includegraphics[width=.9\columnwidth]{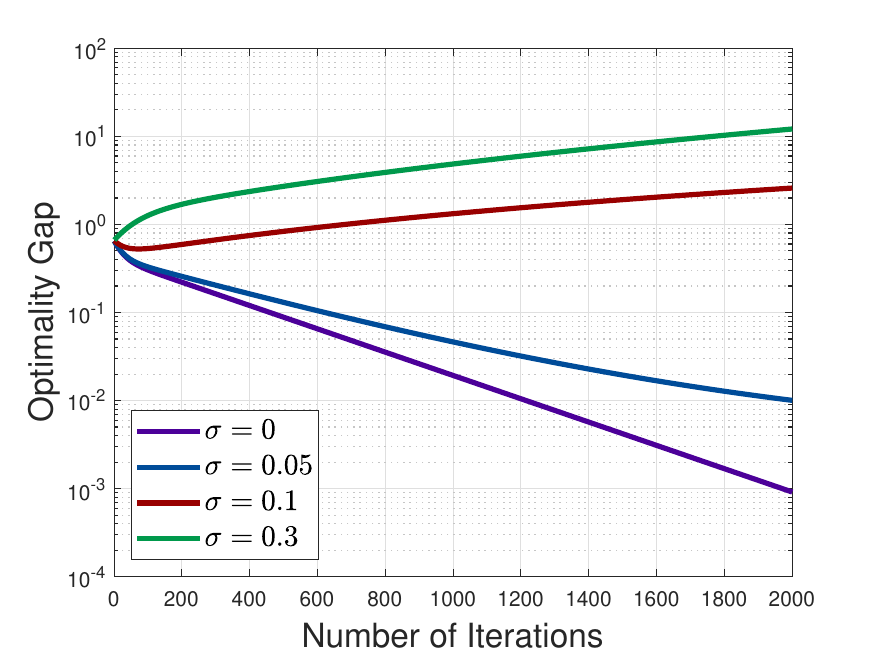}
        \caption{Performance of feedback-based optimization schemes that use inexact sensitivity estimations}
        \label{fig:inexact_sens_performance}
    \end{figure}
\end{example}

Faced with elusive sensitivities, fortunately, we may possess sufficiently rich input-output data collected during historical operations, or we can perform strategic exploration and obtain helpful data on the fly. It is meaningful to tackle the challenges posed by unknown sensitivities by extracting useful patterns from data and exploiting such patterns in the updates. 

We review two main strategies to achieve data-driven and model-free operations. The first strategy consists in learning sensitivities from offline or online data and then plugging them into the first-order updates of the controller. Nonetheless, there also exist scenarios where sensitivities are highly state-dependent and hence hard to learn (e.g., for nonlinear dynamic plants) or even non-existent (e.g., for a non-differentiable steady-state map $h(u)+d$). To handle these issues, the second strategy circumvents learning sensitivities and instead utilizes derivative-free optimization techniques.

\subsubsection{\tb Online Sensitivity Learning}\label{subsubsec:sens_learn}
The knowledge of steady-state input-output sensitivities is a key requirement for the above-mentioned feedback-based optimization schemes. However, such knowledge may not be easily accessible in some applications because of the lack of exact models and accurate parameters. To address this issue, we can utilize offline data or online interactions to perform sensitivity learning. Afterward, we plug these sensitivity surrogates in the gradient-based updates of feedback-based optimization schemes.

We review two approaches to sensitivity learning. The first strategy leverages Willems' fundamental lemma \citep{willems2005note} and uses open-loop input-output trajectories to parametrize the sensitivity matrix of an LTI plant \citep{nonhoff2021data,nonhoff2022online,bianchin2021online}. Provided that finite-length trajectories are informative enough (i.e., they satisfy a persistence-of-excitation condition), the Hankel matrices constructed from such trajectories characterize the whole input-output behavior of an LTI plant (see the review in \cite{markovsky2021behavioral,markovsky2023data}). Remarkably, we bypass the step of identifying system matrices and perform control from raw data in a direct manner. The extension to handle general nonlinear plants remains to be explored.

An alternative strategy learns the sensitivities of nonlinear algebraic maps through recursive least squares estimation. The key idea is to feed the data collected during online implementations to a recursive estimator, thus gradually improving the learning accuracy \citep{picallo2021adaptive,dominguez2023online}.
With the iterative adjustments of feedback-based optimization, we collect a series of inputs and the corresponding measured outputs. The differences between the input-output pairs at neighboring time steps provide rich information on the rates of change, i.e., sensitivities. Since such data arrives in an online fashion, we perform recursive estimation to learn more and more accurate sensitivities. We illustrate this idea through the following example with nonlinear algebraic maps.

\begin{example}\label{eq:exam_sens_learn}
   We revisit \cref{ex:simp_grad} and consider a discrete-time plant described by a differentiable steady-state map $y=h(u)+d$. The objective $\Phi(y)$ is a function of the steady-state output $y$ and is equivalent to the reduced cost $\tilde{\Phi}(u) := \Phi(h(u)+d)$. For the output $y_{k+1}$ at time $k+1$, we analyze its first-order Taylor approximation 
   \begin{equation*}
      y_{k+1} = h(u_{k+1}) + d \approx \underbrace{h(u_k) + d}_{=y_k} + \nabla h(u_k)(u_{k+1}-u_k),
   \end{equation*}
   where we consider a small change of the input $\Delta u_k = u_{k+1} - u_k$. Hence, the map from $\Delta u_k = u_{k+1} - u_k$ to the change of the output $\Delta y_k = y_{k+1} - y_k$ is dominated by the linear sensitivity, i.e., $\Delta y_k \approx \nabla h(u_k) \Delta u_k$, subject to errors due to approximation. Let $s_k = \operatorname{vec}(\nabla h(u_k))$ and $U_{\Delta,k} = \Delta u_k^{\top} \otimes I$, where $\operatorname{vec}$ and $\otimes$ denote the column-wise vector representation and the Kronecker product, respectively. The dynamics of the sensitivity $s_k$ and the measurement equation are given by
   \begin{equation}\label{eq:sensitivity_dynamics}
   \begin{split}
      s_{k+1} &= s_k + w_{p,k}, \\
       \Delta y_k &= U_{\Delta,k} s_k + w_{m,k},
   \end{split}
   \end{equation}
   where $w_{p,k}$ and $w_{m,k}$ represent the change and the approximation error, respectively. \cite{picallo2021adaptive} assume that they are modeled by Gaussian noises with covariance matrices $\Sigma_{p,k}$ and $\Sigma_{m,k}$. To reflect the dependence of $s_k$ on $u_k$, it is further assumed that $\Sigma_{p,k}$ is proportional to $\|\Delta u_k\|^2$. Then, learning sensitivity based on \eqref{eq:sensitivity_dynamics} amounts to the classical problem of state estimation. The online sensitivity learning strategy \citep{picallo2021adaptive} based on the Kalman filter is as follows
   \begin{equation}\label{eq:kf_sens_learn}
   \begin{split}
      \hat{s}_{k+1} & =\hat{s}_k + K_{k}\left(\Delta y_{k}-U_{\Delta,k} \hat{s}_k\right), \\
      \Sigma_{k+1} & =\left(I-K_{k} U_{\Delta,k}\right)\Sigma_{k} + \Sigma_{p,k}, \\
      K_k &= \Sigma_k U_{\Delta, k}^{\top}\left(\Sigma_{m,k}+ U_{\Delta k} \Sigma_k U_{\Delta k}^{\top}\right)^{-1}.
    \end{split}
    \end{equation}
    Afterward, we convert the vector representation $\hat{s}_k$ to the matrix form $\widehat{\nabla} h(u_k)$ and plug it into the updates of the controller. As illustrated by Figure~\ref{fig:closed_loop_sens_learn}, the closed-loop interconnection is
    \begin{align}\label{eq:sensitivity_learn_closed_loop}
    \begin{split}
      \text{plant} & \begin{cases}
        y_k & = h(u_k) + d
      \end{cases} \\
      \text{controller} & \begin{cases} u_{k+1} &= u_k - \eta \widehat{\nabla} h(u_k)^{\top} \nabla \tilde{\Phi}(y_k) + v_k,
      \end{cases}
    \end{split}
    \end{align}
    where $\eta >0$ is the step size, $v_k$ is a random excitation noise, and $\widehat{\nabla} h(u_k)$ is obtained from \eqref{eq:kf_sens_learn}. \cite{picallo2021adaptive} show that the interconnection of \eqref{eq:kf_sens_learn} and the plant ensures that $\widehat{\nabla} h(u_k)$ is an asymptotically unbiased estimate of the sensitivity $\nabla h(u_k)$ with a bounded variance. The bound on the solution accuracy (i.e., the distance of $u_k$ to the optimal point) depends on this variance.

    The strategy \eqref{eq:kf_sens_learn} uses one previous measurement. By contrast, the recursive sensitivity learning strategy in \cite{dominguez2023online} exploits all the historical data by assigning decaying weights to old measurements.

    For the extension to nonlinear dynamic plants, the central requirement of \eqref{eq:sensitivity_learn_closed_loop} is to achieve sufficiently accurate and fast learning of sensitivities. Otherwise, the controller will suffer from sub-optimality due to cumulative errors in gradient steps.
    \begin{figure}[!bt]
    \centering
         \resizebox{\columnwidth}{!}{
    \begin{tikzpicture}
      \matrix[ampersand replacement=\&, row sep=0.3cm, column sep=.55cm] {
        \node[smallsum](sum1){}; \& \node[branch](britr) {}; \\ \\
         
        \& \node[block](delay) {$z^{-1}$}; \& \& \node[branch](brsens) {}; \\

        \node[smallsum](gd){}; \& \& \node[block](senslearn) {$\begin{gathered}
            \text{Sensitivity} \\
            \text{learning}~\eqref{eq:kf_sens_learn}
        \end{gathered}$}; \& \node[block] (plant) {
          $\hat{y}_k = h(u_k)$
        }; \\

        \node[gainup] (gain) {$\eta$}; \& \& \& \node[smallsum](output){}; \\

        \node[block] (gradient) {
          $-\widehat{\nabla}h(u)^{\top} \nabla \tilde{\Phi}(y)$
        }; \& \node[branch,xshift=.5cm](edgeout) {}; \\
      };
    
      \draw[connector] ([xshift=-.5cm]sum1.west)--(sum1.west) node[at start, left]{$v_k$} node[at end, above left] {$+$};
      \draw[connector] ([xshift=.5cm]output.east)--(output.east) node[at start, right]{$d$};
      \draw[connector] (gradient.north)--(gain.south);
      \draw[connector] (gain.north)--(gd.south) node[at end, below right] {$-$};
      \draw[connector] (gd.north)--(sum1.south) node[at end, below right] {$+$};
      \draw[connector] (britr.south)--(delay.north) node[at start, below right] {$u_k$};
      \draw[connector] (delay.south)|-(gd.east) node[at end, above right] {$+$};
      \draw[connector] (sum1.east)-|(plant.north);
      \draw[connector] (plant.south)--(output.north);
      \draw[connector] (output.south)|-(gradient.east) node[at start, below left] {$y_k$};
      \draw[connector] (output.west)-|(senslearn.south) node[at end, below left] {$\Delta y_k$};
      \draw[connector] (brsens.west)-|(senslearn.north) node[at end, above left] {$\Delta u_k$};
      \draw[line] (senslearn.west)-|(edgeout.north) node[at end, above right] {$\widehat{\nabla}h(u_k)$};
    \end{tikzpicture}
    }
    \caption{Feedback-based optimization with sensitivity learning via recursive least squares \eqref{eq:sensitivity_learn_closed_loop}}
    \label{fig:closed_loop_sens_learn}
    \end{figure}
\end{example}

\subsubsection{\tb Model-Free Control via Derivative-Free Optimization}
Confronted with general cases involving nonlinear dynamic plants or non-differentiable steady-state maps, we may not be able to learn sensitivities accurately enough. This restriction may jeopardize the overall performance of feedback-based optimization.

A natural angle of attack is to give up learning sensitivities altogether and pursue entirely model-free operations. To achieve this goal, we consider feedback-based optimization without gradient evaluations deeply rooted in derivative-free optimization. The precursors include the Nelder–Mead simplex method \citep{nelder1965simplex} and evolutionary optimization algorithms \citep{simon2013evolutionary} (see \cite{conn2009introduction} for an extensive review). The paradigms more related to this review are simultaneous perturbation stochastic approximation \citep{spall1992multivariate} and its modern variant, i.e., zeroth-order optimization (ZO) \citep{nesterov2017random,liu2020primer}, as well as Bayesian optimization \citep{frazier2018tutorial,shahriari2016taking}. We provide a more detailed account of model-free feedback-based optimization schemes based on ZO, because they i) handle more general problems involving high-dimensional systems, and ii) ensure closed-loop stability when interconnected with dynamic plants (similar to the timescale separation argument in \cref{subsec:singular_perturb}).

\paragraph{Zeroth-Order Optimization-based Schemes}
The key idea of ZO is to construct stochastic gradient estimates from the evaluations of objective values, and then set update directions as negative gradient estimates. Zeroth-order optimization algorithms enjoy comparable convergence rates as their first-order counterparts, albeit further exhibiting polynomial dependence on the dimension of the problem \citep{flaxman2005online,duchi2015optimal,shamir2017optimal,zhang2022new,chen2022improve}.

Recall that we aim to realize fully model-free feedback-based optimization. The insight from ZO is that we can utilize inputs and real-time outputs to evaluate objective values, construct gradient estimates, and then iteratively update control inputs \citep{chen2020model,he2022model,tang2023zeroth,poveda2017robust}. While constructing gradient estimates, however, we cannot obtain numerous independent and identically distributed samples from the plant, because every input actuation changes the operating condition. Hence, we need to judiciously use samples and select inputs.
In the following example, we illustrate the model-free designs of feedback-based optimization.

\begin{example}\label{example:model_free_scheme}
    Similar to \cref{eq:exam_sens_learn}, we consider a discrete-time plant with the steady-state map $y = h(u) + d$. Let $\tilde{\Phi}(u) := \Phi(h(u) + d)$ denote the reduced cost function. We assume that $\tilde{\Phi}(u)$ is smooth and Lipschitz continuous.

    We first discuss how to construct a stochastic gradient estimate. The key idea is to explore $\tilde{\Phi}(u)$ at perturbed points and scale a random vector with quantities related to objective values to form the estimate.
    Let $v \sim \mathcal{N}(0,\frac{1}{p}I)$ be a $p$-dimensional vector sampled from the normal distribution. Consider the directional derivative $\langle\nabla \tilde{\Phi}(u),v\rangle$ of $\tilde{\Phi}(u)$ along $v$. We use the finite differences of function values to approximate it as follows
    \begin{equation*}
        \langle\nabla \tilde{\Phi}(u),v\rangle \approx \frac{\tilde{\Phi}(u+\delta v)-\tilde{\Phi}(u)}{\delta},
    \end{equation*}
    where $\delta > 0$ is the smoothing radius and is typically chosen sufficiently small. Since $v$ is stochastic and satisfies $\E_v\left[vv^{\top}\right] = \frac{1}{p}I$, we have
    \begin{equation*}
        \E_v\left[p\langle\nabla \tilde{\Phi}(u),v\rangle v\right] = \E_v\left[pvv^{\top}\nabla \tilde{\Phi}(u)\right] = \nabla \tilde{\Phi}(u),
    \end{equation*}
    which implies that $p\langle\nabla \tilde{\Phi}(u),v\rangle v$ is, in expectation, an unbiased estimate of the gradient $\nabla \tilde{\Phi}(u)$. Thus, we replace $\langle\nabla \tilde{\Phi}(u),v\rangle$ with the aforementioned finite differences and obtain the following two-point gradient estimates \citep{duchi2015optimal}
    \begin{equation}\label{eq:grad_est_two_pt}
        \phi^{(2)}(u) = \frac{p}{\delta}(\tilde{\Phi}(u+\delta v)-\tilde{\Phi}(u))v.
    \end{equation}
    The above estimates are still useful in the sense of expectation. In fact, $\E_v[\phi^{(2)}(u)]$ and $\E_v[\tilde{\phi}^{(2)}(u)]$ equal to the gradient of the smooth approximation $\tilde{\Phi}_{\delta}(u)$ of the objective function $\tilde{\Phi}(u)$ \citep{nesterov2017random}, where
    \begin{equation*}
        \tilde{\Phi}_{\delta}(u) = \E_v\left[\tilde{\Phi}(u+\delta v) \right].
    \end{equation*}
    Moreover, the gradients and the values of $\tilde{\Phi}_{\delta}(u)$ and $\tilde{\Phi}(u)$ can be made arbitrarily close by selecting a $\delta$ sufficiently close to $0$.
    
    The gradient estimates in \eqref{eq:grad_est_two_pt} require two evaluations of objective values. We need to actuate the plant twice with $u+\delta v$ and $u$, measure the corresponding steady-state outputs $y$ and $y'$, and then evaluate the objective values $\tilde{\Phi}(y)$ and $\tilde{\Phi}(y')$. Such designs work well for algebraic maps, where steady-state outputs are directly available \citep{chen2020model,poveda2017robust,tang2023zeroth}. Nonetheless, they may not be suitable for real-time decision-making with dynamic plants. In that case, within a short interval, it would be challenging to inject two different inputs and measure the corresponding steady-state outputs. That is, every time interval may only allow one actuation step and return one dynamic measurement.

    The key to real-time operations lies in utilizing one new evaluation of objective values. A close look at $\E_v[\phi^{(2)}(u)]$ in \eqref{eq:grad_est_two_pt} helps to achieve this goal. Indeed,
    \begin{equation*}
      \begin{split}
        \E_v[\phi^{(2)}(u)] &= \E_v \left[\frac{p}{\delta}\tilde{\Phi}(u+\delta v)v - \frac{p}{\delta}\tilde{\Phi}(u)v \right] \\
          &= \E_v \left[\frac{p}{\delta}\tilde{\Phi}(u+\delta v)v\right],
      \end{split}
    \end{equation*}
    where we use the fact that $\E_v[v] = 0$ for $v \sim \mathcal{N}(0,\frac{1}{p}I)$, and we shall design the iterative updates such that $u$ is independent of $v$ (see \eqref{eq:model_free_control} below). The implication is that the one-point estimate
    \begin{equation}\label{eq:grad_est_one_pt}
        \phi^{(1)}(u) = \frac{p}{\delta}\tilde{\Phi}(u+\delta v)v
    \end{equation}
    suffices to be a sensible estimate of the gradient \citep{flaxman2005online}, because its expectation equals that of the two-point estimate $\phi^{(2)}(u)$. Nonetheless, one issue of $\phi^{(1)}(u)$ is its potentially high variance. The main reason is that the magnitude of $\phi^{(1)}(u)$ depends on the objective value $\tilde{\Phi}(u+\delta u)$, whose range of change may be large. To solve this issue, we can subtract the term $\tilde{\Phi}(u+\delta u)$ in $\phi^{(1)}(u)$ with an appropriate value that is independent of $v$, thus reducing the variance of the gradient estimate. Notably, \cite{zhang2022new} proposes the following one-point residual-feedback estimate
    \begin{equation}\label{eq:grad_est_one_pt_residual}
        \tilde{\phi}^{(1)}(u) = \frac{p}{\delta}\left(\tilde{\Phi}(u+\delta v) - \tilde{\Phi}(u'+\delta v')\right)v,
    \end{equation}
    which subtracts \eqref{eq:grad_est_one_pt} with the objective value at the previous point $u'+\delta v'$ during iterations, and $v'\sim \mathcal{N}(0,\frac{1}{p}I)$ is independent of $v$. Furthermore, \cite{chen2022improve} exploits insights into high-pass and low-pass filters in extremum seeking \citep{ariyurRealTimeOptimization2003} and presents a unifying structure for one-point estimates. This structure encompasses residual feedback for variance reduction and a momentum term for acceleration.
    
    Motivated by the above observations, we consider the model-free feedback-based optimization scheme \citep{he2022model} illustrated by Figure~\ref{fig:FO_model_free}. The closed-loop dynamics are as follows
    \begin{figure}[!bt]
    \centering
    \resizebox{\columnwidth}{!}{
    \begin{tikzpicture}
      \matrix[ampersand replacement=\&, row sep=0.3cm, column sep=.55cm] {
        \&\node[smallsum](sum1){};
        \&\&\& \node[none](edgetop) {}; \\

        \node[gainup](gain1) {$\delta$}; \& \node[branch](br2){}; \& \\

        \node[smallsum](sum2){$\times$}; \& \& \node[block](delay) {$z^{-1}$}; \& \node[none](disturbance) {$d$}; \& \node[block] (plant) {
          $\begin{aligned}
             \zeta_{k+1} & = f(\zeta_k, u_k) \\
              y_k & = g(\zeta_k, d)
            \end{aligned}$
        }; \\
        
        \node[none](source){$v_k$}; \& \node[smallsum](gd){}; \&
        \&  \&  \\

        \node[block](rnd){RND}; \& \node[gainup](gain2) {$\eta$}; \& \& \node[dashedblock](costprev) {$\Phi(y_{k-1})$}; \& \node[none](edgebot) {}; \\
        \\

        \node[none](source2){$v_{k-1}$}; \& \node[block](gdestimate){$\phi_{k-1}$}; \& \node[branch](brgdest){};
        \& \node[block] (cost) {$\Phi(y_k)$};
        \&   \\
      };

     \draw[connector] (source.north)--(sum2.south);
     \draw[connector] (sum2.north)--(gain1.south);
     \draw[connector] (gain1.north)|-(sum1.west) node[at end, above] {$+$};
     \draw[connector] (sum1.east)-|(plant.north) node[at start, above right] {$u_k$};
     \draw[connector] (disturbance.east)|-(plant.west);
     \draw[connector] (plant.south)|-(cost.east);
     \draw[connector] (cost.west)|-(gdestimate.east);
     \draw[connector] (source2.east)--(gdestimate.west);
     \draw[connector] (rnd.north)--(source.south);
     \draw[connector] (rnd.south)--(source2.north);
     \draw[connector] (gdestimate.north)--(gain2.south);
     \draw[connector] (gain2.north)--(gd.south) node[at end, below right] {$-$};
     \draw[connector] (delay.south)|-(gd.east) node[at end, above right] {$+$};
     \draw[connector] (gd.north)--(sum1.south) node[at end, below left] {$+$};
     \draw[connector] (br2.east)-|(delay.north) node[at start, above right] {$w_k$};
     \draw[connector,dashed] (costprev.west)-|(brgdest.north);
     \draw[connector,dashed] (edgebot.west)--(costprev.east);
    \end{tikzpicture}
    }
    \caption{Model-free feedback-based optimization \eqref{eq:model_free_control} in a closed loop}
    \label{fig:FO_model_free}
    \end{figure}
    \begin{align}\label{eq:model_free_control}
    \begin{split}
      \text{plant} & \begin{cases}
        \zeta_{k+1} & = f(\zeta_k, u_k) \\
        y_k & = g(\zeta_k, d)
      \end{cases} \\
      \text{controller} & \begin{cases} w_{k+1} &= w_k - \eta \phi_k \\
      \phi_k &= \frac{p}{\delta} \left(\Phi(y_{k+1}) - \Phi(y_k)\right) v_k \\
      u_{k+1} &= w_{k+1} + \delta v_{k+1},
      \end{cases}
    \end{split}
    \end{align}
    where $\eta>0$ is the step size, $k=0,1,\ldots$, and $v_0,v_1,\ldots$ are independent and identically distributed random vectors drawn from $\mathcal{N}(0,\frac{1}{p}I)$. In \eqref{eq:model_free_control}, $w_{k+1}$ is the candidate solution and is iteratively updated based on the gradient estimate $\phi_k$ at the previous solution $w_k$. The estimate $\phi_k$ is similar to \eqref{eq:grad_est_one_pt_residual}. The difference is that due to real-time operations, we substitute the steady-state output $h(u_k,d)$ with the measurement output $y_{k+1}$ after applying the input $u_k$. Consequently, $\Phi(y_{k+1})$ is an approximation of the true objective value $\tilde{\Phi}(u_k) = \Phi(h(u_k,d))$. Afterward, the controller adds an exploration noise $\delta v_{k+1}$ to $w_{k+1}$ and applies the input $u_{k+1}$ to the plant. To establish the stability of the closed-loop interconnection \eqref{eq:model_free_control}, we can extend the timescale separation argument discussed in \cref{subsec:singular_perturb}. A further extension to handle input constraints via Frank-Wolfe type updates is also available \citep{he2022model}.
\end{example}

By and large, model-free methods in \cref{example:model_free_scheme} are less sample-efficient than feedback-based optimization schemes with perfect sensitivities \citep{hauswirthTimescaleSeparationAutonomous2020,bernsteinOnlinePrimalDualMethods2019,colombino2020online}. In other words, they require relatively more actuation steps to reach a certain solution accuracy, see \cref{fig:comparison} for a detailed comparison.
One intuition is that model-free methods drop the useful structural information contained in prior knowledge or approximate models and instead rely on stochastic exploration. This stochasticity may influence the overall convergence rates and cause an increased number of actuations. Another interpretation is that not all history is taken into account at every iteration, which differs from the recursive estimation of sensitivity in \cref{subsubsec:sens_learn}. This issue may be a disadvantage when the sensitivity is almost constant.
\begin{figure}[!tb]
    \centering
    \includegraphics[width=0.9\columnwidth]{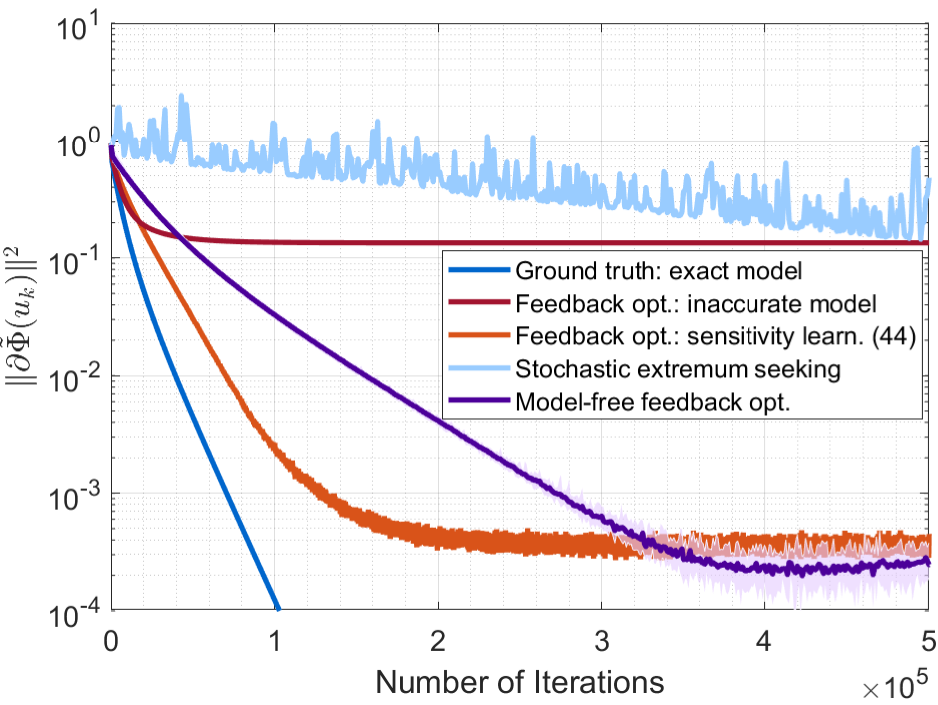}
    \caption{Comparison of the model-free feedback controller in \eqref{eq:model_free_control}, feedback controllers with perfect sensitivities or with sensitivity learning in \eqref{eq:sensitivity_learn_closed_loop}, and a stochastic extremum seeking algorithm in \citep{liu2016stochastic}}
    \label{fig:comparison}
\end{figure}

\paragraph{Bayesian Optimization-based Schemes}
Another line of model-free feedback-based optimization schemes is built on Bayesian optimization \citep{frazier2018tutorial,shahriari2016taking}. The key idea is to maintain probabilistic surrogate models of unknown objectives and constraints based on samples, and then set the next query as the optimal point of an acquisition function. This update mechanism is shown to be sample-efficient and hence attractive for applications where data collection is costly \citep{neumann2020data,simonetto2021personalized}. Another appealing part is that black-box constraints are satisfied with a high probability \citep{krishnamoorthy2023model}. It benefits from the design that balances exploration and exploitation by judiciously expanding a safe region and selecting samples from this region \citep{sui2015safe,krishnamoorthy2022safe}. The closed-loop implementation requires the detection of steady states \citep{krishnamoorthy2023model}, which implies a complete time-scale separation. Furthermore, it is possible to incorporate contextual information to better adapt to time-varying external conditions \citep{xu2023violation}.

In general, Bayesian optimization-based methods are mostly suitable for low-dimensional systems. Otherwise, the computational cost of using Gaussian processes as surrogate models of objectives and constraints is significant \citep{shahriari2016taking}.

\begin{remark}
    The aforementioned model-free schemes utilize exact function evaluations at perturbed candidate solutions. In some applications (e.g., optimization of personalized utility functions with user's feedback \citep{simonetto2021personalized,notarnicola2023distributed}), only infrequent and noisy evaluations are available. To address this issue, we can construct approximations based on parametric estimation with basis functions \citep{cothren2022data} or neural networks \citep{cothren2022online}. Then, we optimize these approximations via iterative descent algorithms. In this line of work, the considered objective function is decoupled and contains two parts that depend on the input and the output. Approximations are separately constructed for these two parts. Therefore, these schemes still require knowing steady-state input-output sensitivities.
\end{remark}
}

\section{{\tb Application Example: Optimal Reserve Dispatch in Electricity Grids}}\label{sec:redispatch}

{\tb So far, we have explored} small educational examples illustrating various individual concepts and special cases. 
In this section, we illustrate the potential of these methods for the emerging problem of optimal real-time operation of future power systems in the face of increasing intermittent renewable energy generation.

Future power systems will face less predictable and more volatile conditions, which require new control and decision protocols to continue operating power grids safely and efficiently.
In a simplified and idealized fashion, this task of optimal real-time grid operation can be cast as tracking the solution of a (time-varying) \ac{acopf} problem \citep{frankOptimalPowerFlow2012a,frankOptimalPowerFlow2012,huneaultSurveyOptimalPower1991}.
This perspective on the problem has motivated and guided an essential part of the recent research on online feedback optimization and has contributed to defining many of the methodological questions that we addressed earlier in this article.

Although solving \ac{acopf} problems in real-time has been proposed for several decades, the idea of driving the state of a power grid to the solution of an \ac{acopf} by interconnecting an optimization dynamic in closed loop with the physical grid is much more recent. 
In the context of distribution grids, the idea of using iterative measurements (rather than a model) to determine a gradient direction and ultimately implement an iterative optimization algorithm was first explored in \cite{Bolognani2013DistributedControlStrategy} and in \cite{Farivar2013}, with a focus on deriving distributed and decentralized control strategies.
The idea of employing gradient-based online algorithms in closed loop with the power grid to drive the system to the solution of a given \ac{acopf} was further explored (with the addition of constraints, time-varying parameters, and more general problem formulations) in \cite{Bolognani2015ReactivePowerFeedbackControl,bernstein2015composable,reyes2015composable,ganOnlineGradientAlgorithm2016,tangRealtimeOptimalPower2017,mazziOnlineOptimizationAlgorithm2018a,Bernstein2019} among others.

A connected line of research, also motivated by power flow optimizations in the electricity grid, is based on the insights of \cite{bolognaniFastPowerSystem2015} that identifies the solutions space as a smooth manifold. Using this perspective, \cite{hauswirthProjectedGradientDescent2016} proposes a (continuous-time) gradient-based controller that results in a differential-algebraic system. This algorithm (in a discretized form) was applied to \ac{acopf} problems in \cite{hauswirthOnlineOptimizationClosed2017}, and formal convergence guarantees of the discretization were established in \cite{haberleNonconvexFeedbackOptimization2021}.
While considering the solution space of the AC power flow equations as a smooth manifold (rather than assuming an input-output map) adds to the conceptual complexity, it also highlights underlying challenges relating to solvability, controllability, and the possibility for bifurcations.

Besides purely gradient-based schemes, saddle-point algorithms (as defined in \cref{sec:saddle}) have also been extensively studied in the context of real-time feedback-based power system optimization. However, the AC OPF problem's non-convexity makes formal convergence guarantees difficult to establish. For this reason,
\cite{dallaneseOptimalPowerFlow2018} resorts to a dual augmentation term, which improves convergence but perturbs the equilibria such that they do not necessarily correspond to critical points of the optimization problem anymore. 

Further research on this topic has then produced several contributions, and we refer to Section IV in \cite{molzahnSurveyDistributedOptimization2017} for a review of these works.
This feedback-optimization approach to the real-time operation of power grids has reached a remarkable level of technological maturity \citep{Kroposki2020} and has been validated both in laboratory experiments \citep{experimentDC,ortmann2020experimental} and in real-world deployments \citep{ortmann2023deployment}. 

In the following, we do not delve into the peculiarities of any of these implementations; instead, we briefly describe one possible version of this problem, we show how the problem can be brought into the general form we studied in the previous sections, we apply the results and algorithms we reviewed, and we discuss why a feedback-based optimization solution is so successful in this context.

\subsection{{\tb Formulation as a Real-Time Optimal Power Flow Problem}}
Consider an AC power transmission network as an undirected graph with a set $\calN$ of buses and a set $\calM$ of transmission lines. We denote by $(l,k) \in \calM$ a line connecting bus $l$ to bus $k$. As decision variables, each bus $l \in \calN$ has an associated voltage magnitude $v_l$, voltage angle $\theta_l$, and active and reactive power generation $p^\G_l, q^\G_l$. By $p^\G, q^\G$, etc. we denote the vectors obtained from stacking all respective components $p^\G_l$ for all $l \in \calN$.

    {\tb For simplicity, we assume that the objective is given as} the cost $\Phi(p^\G) = \sum_{l \in \calN} \Phi_l(p^\G_l)$ for active power {\tb generation} as the minimization objective.%

Hence, a basic \ac{acopf} problem is given by
\begin{subequations}\label{eq:basic_opf}
    \begin{align}
        \underset{v,\theta, p^\G, q^\G}{\minimize} \quad & \Phi(p^\G) \label{eq:acopf_cost}                                                                                        \\
        \underset{\forall l \in \calN}{\subjto} \quad
                                                         & 0 = p_l^\G - p_l^\L  - \sum\nolimits_{l \rightarrow k} p_{lk}(v_l, v_k, \theta_l, \theta_k)\label{eq:acopf_p}           \\
                                                         & 0 = q_l^\G - q_l^\L - \sum\nolimits_{l \rightarrow k} q_{lk}(v_l, v_k, \theta_l, \theta_k)\label{eq:acopf_q}            \\
                                                         & \underline{p}_l \leq  p^\G_l \leq \overline{p}_l \label{eq:acopf_plim}                                                  \\
                                                         & \underline{q}_l \leq  q^\G_l \leq \overline{q}_l \label{eq:acopf_qlim}                                                  \\
                                                         & \underline{v}_l \leq  v_l \leq \overline{v}_l \label{eq:acopf_vlim}                                                     \\
                                                         & i^2_{lk}(v_l, v_k, \theta_l, \theta_k) \leq \overline{i}_{kl}^2 \quad \forall (l,k) \in \calM\, , \label{eq:acopf_ilim}
    \end{align}
\end{subequations}
The inequalities \crefrange{eq:acopf_plim}{eq:acopf_vlim} denote box constraints on power generation and voltages at each bus, and \cref{eq:acopf_ilim} limits the current that flows through the line from $k$ to $\ell$.

\begin{figure}
    \centering
    \includegraphics[width=\columnwidth]{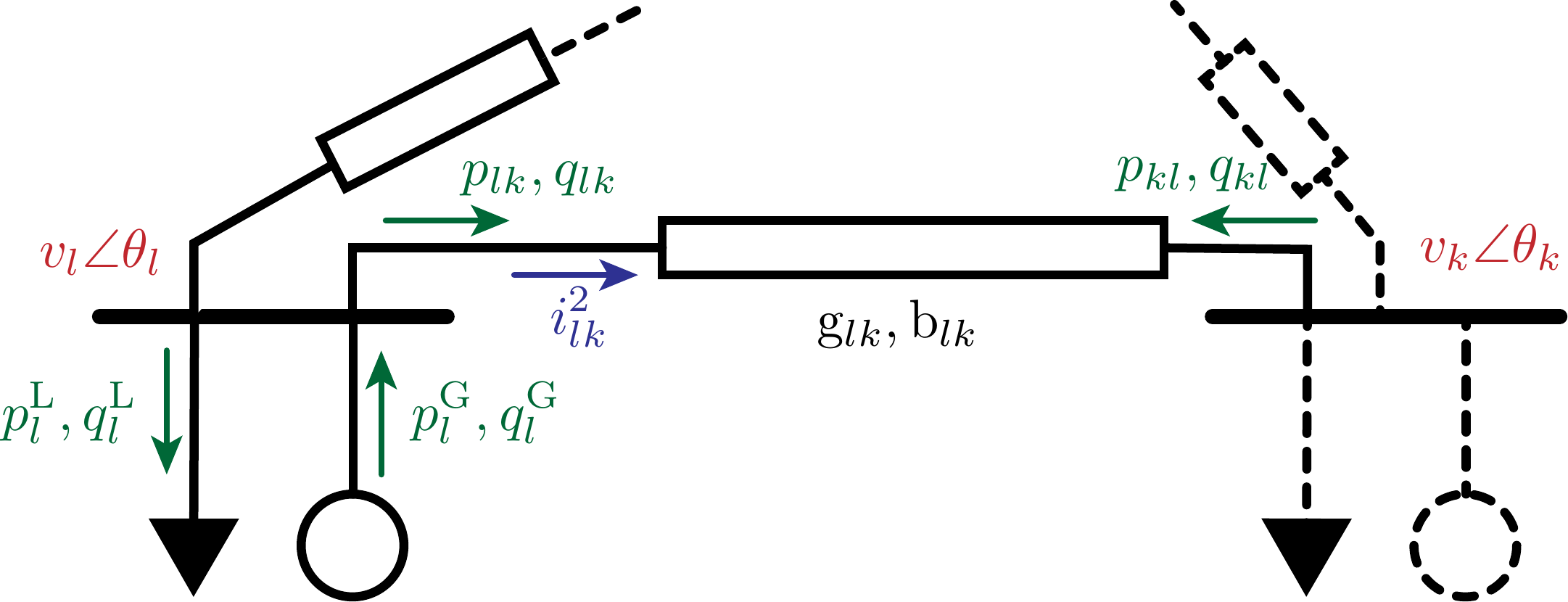}
    \caption{Physical AC power flow quantities}\label{fig:netw_diag}
\end{figure}

The functions $p_{lk}(\cdot), q_{lk}(\cdot)$ and $i^2_{lk}(\cdot)$ denote active and reactive power flow and squared current magnitude from bus $l$ in direction of bus $k$, respectively, as illustrated in \cref{fig:netw_diag}. These nonlinear terms depend on line parameters and constitute the main source of complexity in \ac{acopf} problems. Their explicit expressions (neglecting the so-called \emph{shunt admittances} of the buses) are given by
\begin{align*}
    p_{lk}(v_l, v_k, \theta_l, \theta_k)   & := v_l v_k \left( \rmg_{lk}  \cos(\theta_l - \theta_k) + \rmb_{lk} \sin(\theta_l - \theta_k) \right)  \\
    q_{lk}(v_l, v_k, \theta_l, \theta_k)   & :=  v_l v_k \left( \rmg_{lk}  \sin(\theta_l - \theta_k) - \rmb_{lk} \cos(\theta_l - \theta_k) \right) \\
    i^2_{lk}(v_l, v_k, \theta_l, \theta_k) & := (\rmg_{lk}^2 + \rmb_{lk}^2)\left(v_k^2 + v_l^2 - 4 v_k v_l \cos (\theta_k - \theta_l) \right) \,
\end{align*}
where $\rmg_{lk}$ and $\rmb_{lk}$ denote the \emph{conductance} and \emph{susceptance} of the transmission line connecting buses $l$ and $k$. Note that $ i_{lk}^2 = i_{kl}^2$ holds (without shunt elements), but $p_{lk} = -p_{kl}$ and $p_{lk} = -p_{kl}$ are not generally true. For a more comprehensive introduction to AC power flow including alternative formulations the reader is referred to \cite{molzahnSurveyRelaxationsApproximations2019,frankIntroductionOptimalPower2016}.

\Ac{acopf} problems like \eqref{eq:basic_opf} are well-studied and routinely solved (numerically) in practice. Yet, they remain computationally demanding, because of the nonlinear power flow equations \crefrange{eq:acopf_p}{eq:acopf_q} which render the problem non-convex (although convex relaxations can yield global optimality certificates in some cases; \citealt{molzahnSurveyRelaxationsApproximations2019,lowConvexRelaxationOptimal2014a}).

\subsection{\tb Real-Time Operation as Closed-Loop Optimization}\label{sec:redispatch-fb_form}

In an online setting, \eqref{eq:basic_opf} defines the most desirable operating point of the grid (its power flows and its power injections $p^\G,q^\G$) given the actual power consumption $p^\L,q^\L$ (as opposed to predictions used in offline {day/hour-ahead} calculations). The power consumption $p^\L,q^\L$ and other parameters in \cref{eq:basic_opf} are generally time-varying. Hence, rather than solving \eqref{eq:basic_opf} once for a given configuration, a real-time-operation strategy needs to track the \ac{acopf} problem's solution across time.

From a theoretical viewpoint, robustly tracking solutions of \eqref{eq:basic_opf} presents important challenges:
\begin{enumerate}
    \item Equations \crefrange{eq:acopf_p}{eq:acopf_q} are the steady-state equations of the complicated, interconnected nonlinear dynamics of transmission lines, generation units, loads, and their low-level controllers. Although these dynamics can be assumed to be asymptotically stable, attempts at steering the physical system too aggressively towards a solution of \eqref{eq:basic_opf} can destabilize them.

    \item In general, only an approximate model in the form of the steady-state \acl{acpf} equations \crefrange{eq:acopf_p}{eq:acopf_q} can be employed. A dynamical model for a large-scale power system is not available, as different stakeholders own parts of the system, models are proprietary, operating conditions (e.g., which generation units are online) change over time, and internal states are typically inaccessible.
    Steady-state models have significant parameter uncertainty and are affected by disturbances.

    \item Many physical and engineering constraints of different nature have to be satisfied. Some of these limits, like the generation constraints \crefrange{eq:acopf_p}{eq:acopf_q}, are physical constraints that are strictly enforced by lower-level controllers or through saturation. Other constraints, such as the line flow limits \cref{eq:acopf_ilim}, are thermal limits that can be violated temporarily. Yet other constraints are hard limits that must not be violated at any time. This mainly concerns dynamic and voltage stability limits, which we will not consider in this simplified setting. Simply note that the voltage constraints \cref{eq:acopf_vlim} can be understood as guarding against voltage instability \citep{vancutsemVoltageStabilityElectric1998}.

    \item The nonlinear nature, especially under critical operating conditions, of the power flow equations \crefrange{eq:acopf_p}{eq:acopf_q} call for methods that work in the absence of convexity and for ill-conditioned problems.
\end{enumerate}

With this in mind, we argue that the online feedback-based optimization methods that we reviewed are very well suited for this task.
From a practical perspective, ``closing the loop'' on the \ac{acopf} problem \eqref{eq:basic_opf} offers an opportunity to more closely integrate and combine different power system control tasks and thereby increase economic efficiency, resilience, and autonomy. More concretely, a controller tracking the solution of \eqref{eq:basic_opf} makes complex yet economically efficient re-dispatch decisions in response to unscheduled events to guarantee voltage and line flow limits are respected. Thus, \eqref{eq:basic_opf} should be interpreted as a ``residual'' optimization problem around a pre-planned generation schedule, and constraints \crefrange{eq:acopf_plim}{eq:acopf_qlim} do not necessarily quantify the full production capacity of a generation unit, but rather the amount of dispatchable reserves.

In order to express \eqref{eq:basic_opf} as a problem of the form \eqref{eq:basic_prob}, we need to identify inputs, outputs, and disturbances. For this purpose, for each {generation unit}, the active power output and either the reactive power or voltage magnitude are assumed to be controllable\footnote{One distinguishes between so-called \emph{PQ}- and \emph{PV}-buses. See \cite{hauswirthOnlineOptimizationClosed2017} for a more detailed discussion including the role of the \emph{slack bus} in numerical simulations).} and thus make up the control input $u$. The loads $p^\L, q^\L$ define the disturbance $d$. All remaining quantities form the output $y$.

Consequently, \eqref{eq:basic_opf} can be brought into the form \eqref{eq:basic_prob} where the constraints \crefrange{eq:acopf_plim}{eq:acopf_ilim} are assigned to either $\calU$ or $\calX$ according to whether they apply only to controllable variables or not. Under normal operating conditions, the {local existence and differentiability} of the steady-state map $h(\cdot)$ derived from \crefrange{eq:acopf_p}{eq:acopf_q} is guaranteed by the implicit function theorem \citep{bolognaniFastPowerSystem2015}. 

The following two examples illustrate how feedback-based optimization can be applied to the optimal reserve dispatch problem.

    \begin{example}\label{ex:pw_ex_pd}
        To track the time-varying solution of \cref{eq:basic_opf}, we adopt the feedback controller \eqref{eq:pd_control} that implements a projected and primal augmented saddle-point flow.

        \begin{figure}
            \centering
            \includegraphics[width=\columnwidth]{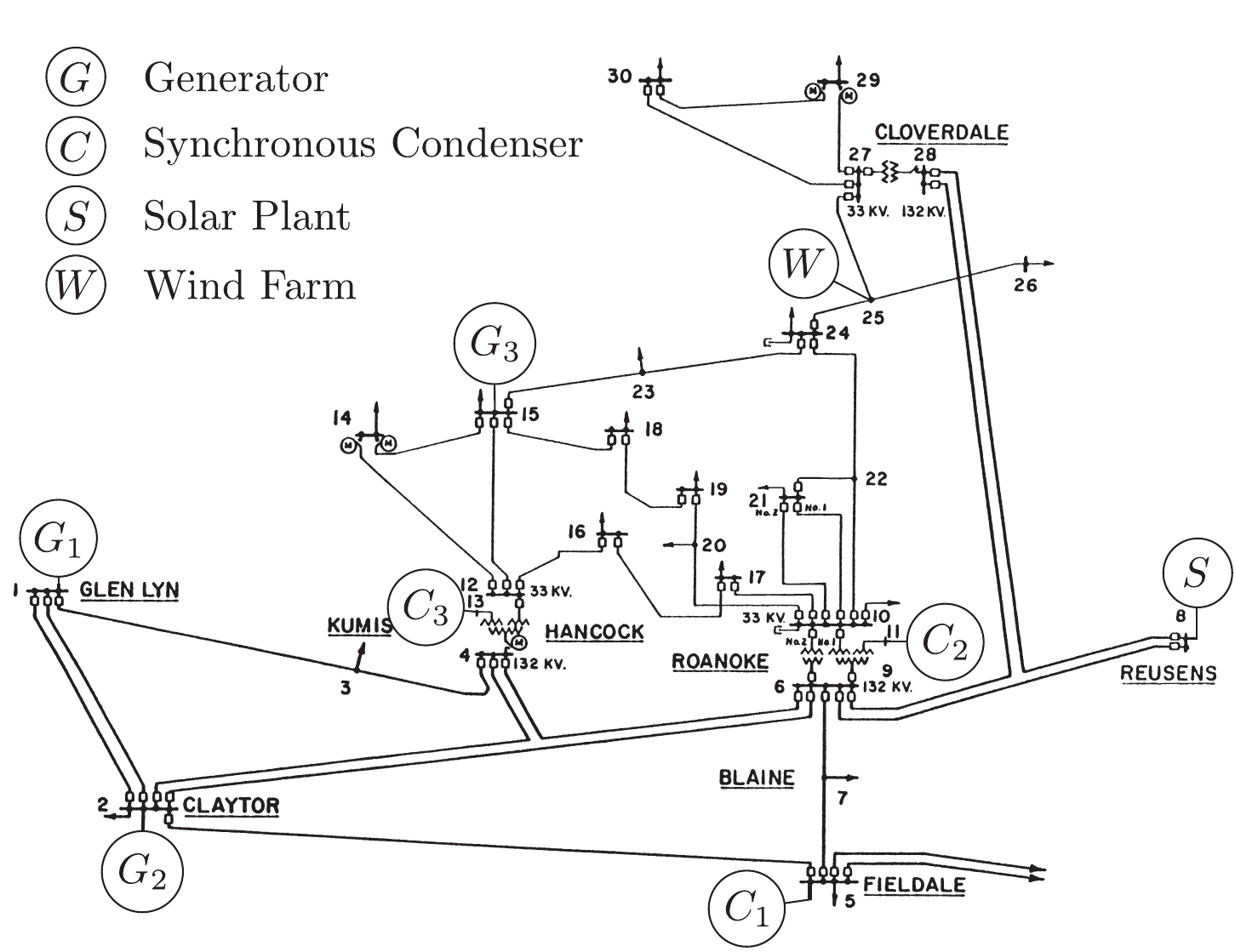}
            \caption{Grid diagram of adapted 30-bus power systems test case used in \cref{ex:pw_ex_pd,ex:lop_survey}.}\label{eq:grid_diag}
        \end{figure}

        As a numerical testbed, we use a modified version of the IEEE 30-bus power systems test case adopted from \cite{hauswirthOnlineOptimizationClosed2017,hauswirthOptimizationAlgorithmsFeedback2020b} and illustrated in \cref{eq:grid_diag}. In particular, a wind and a solar plant provide intermittent generation capacity in addition to three conventional generators. The available solar and wind generation capacity profiles are illustrated in \cref{fig:survey_30bus_profiles}. One of the conventional generators suffers an outage at 4:00 upon which its generation capacity is set to zero for the remainder of the simulation. The goal is to optimally use renewable resources without violating engineering constraints such as voltage limits and thermal line ratings.

        \begin{figure}
            \centering
            \includegraphics[width=\columnwidth]{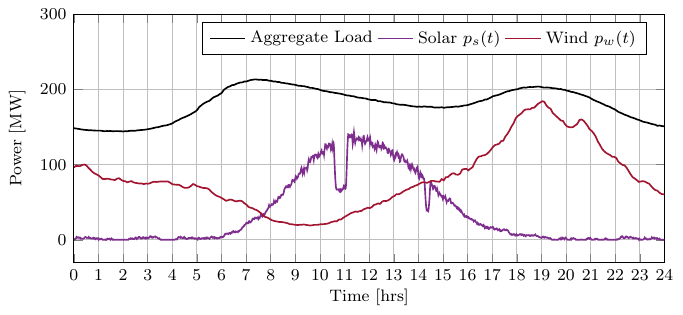}
            \caption{Load and available generation profiles for the modified 30-bus test case in \cref{ex:pw_ex_pd,ex:lop_survey} over a 24h horizon.}
            \label{fig:survey_30bus_profiles}
        \end{figure}

        \begin{figure*}[tb]
            \centering
            \includegraphics[height=3.7cm,trim=0cm 0cm 0cm 0cm,clip]{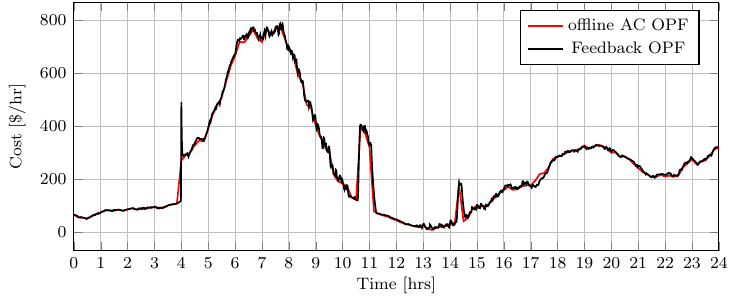}
            \includegraphics[height=3.7cm,trim=.65cm 0cm 0cm 0cm,clip]{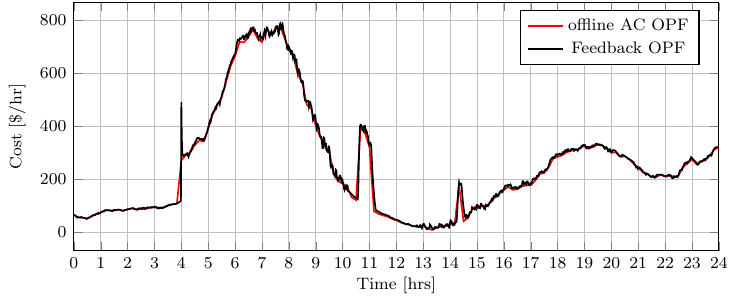}
            \caption{Cost realized by feedback-based optimization in comparison to a-posteriori \ac{acopf} solutions using perfect information. Left: cost achieved by the projected augmented saddle-point flow \eqref{eq:pd_control}; right: cost achieved by the linearized output projection method \cref{eq:lop_up_survey}.}
            \label{fig:survey_30bus_cost}
        \end{figure*}

        \begin{figure*}[tb]
            \centering
            \includegraphics[height=12.2cm,trim=0cm 1.8cm .5cm 0cm,clip]{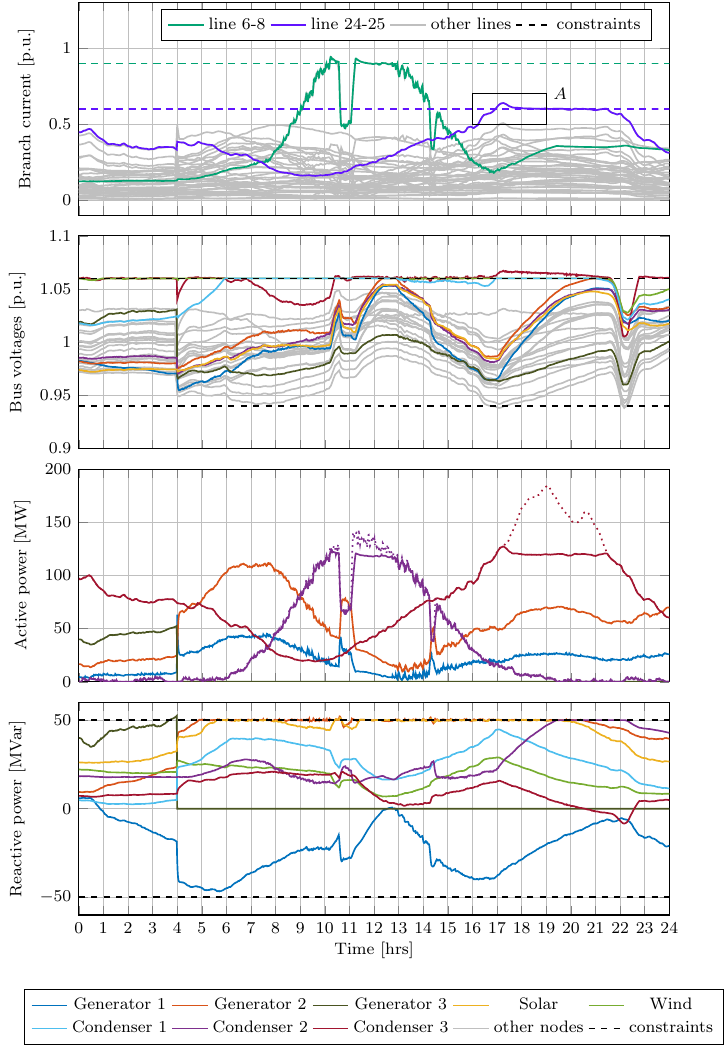}
            \includegraphics[height=12.2cm,trim=.65cm 1.8cm .5cm 0cm,clip]{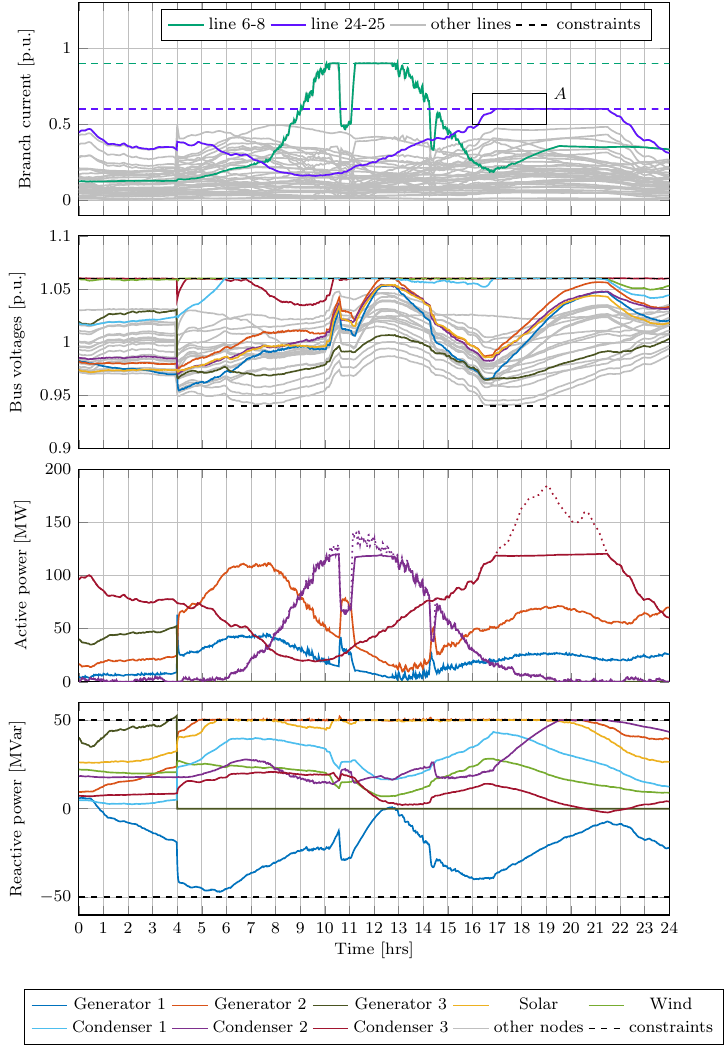}
            \includegraphics[width=1.3\columnwidth,trim=0cm 0cm 0cm 16.5cm,clip]{Figs/arcFigs/30bus_grad_lop_ss_0delay/states_combined}
            \caption{State trajectories realized by feedback-based optimization. Left: trajectories of the projected augmented saddle-point flow \eqref{eq:pd_control}; right: trajectories of the linearized output projection method \cref{eq:lop_up_survey}. See \cref{fig:cstr_viol} for a discussion of the inset A.}\label{fig:survey_30bus_state}
        \end{figure*}

        We implement \eqref{eq:pd_control} as a discretized controller using a simple explicit Euler discretization with a step size corresponding to 1-minute intervals between control actions. 
        The resulting control law is similar to the designs proposed in \cite{dallaneseOptimalPowerFlow2018,tang2018feedback,tangRunningPrimalDualGradient2018,bernsteinOnlinePrimalDualMethods2019}.        
        The physical system is simulated without dynamics, and an off-the-shelf AC power flow solver \citep{zimmermanMATPOWERSteadystateOperations2011} is used to compute a solution to the \ac{acpf} equations \crefrange{eq:acopf_p}{eq:acopf_q} and thus evaluate the input-output map $h(u) + d$. We further compute $\nabla h(u)$ based on input and output measurements and the power flow model in \eqref{eq:basic_opf} (rather than an explicit estimate of $d$). Finally, we choose $Q(u): = {H}(u)^T {H}(u)$ as a metric. This choice of ``implicit'' metric has proven to increase numerical stability in the face of an ill-conditioned input-output map. For a more detailed discussion and motivation, see Chapter 14 in \cite{hauswirthOptimizationAlgorithmsFeedback2020b}.

        The left panel in \cref{fig:survey_30bus_cost} illustrates the (almost perfect) performance achieved by this feedback-based optimization approach in terms of cost compared to the a-posteriori sequential solution of the \ac{acopf} problem (which is based on omniscient and perfect information).

        More importantly, however, the left panel in \cref{fig:survey_30bus_state} shows that, over the entire simulation horizon, constraint violations are very minor and only temporary. The controller achieves this by jointly managing active and reactive power infeed to manage voltage magnitudes and line currents. In particular, both solar and wind generation have to be curtailed to prevent line overloads.

        \begin{figure}
            \centering
            \includegraphics[width=.48\columnwidth]{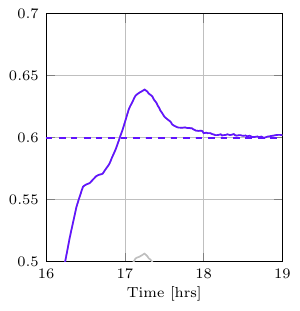}
            \includegraphics[width=.48\columnwidth]{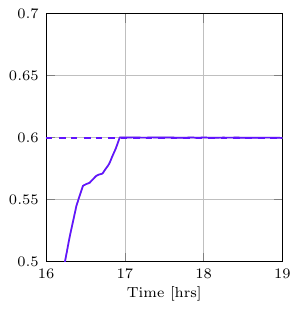}
            \caption{Comparison of line flow limit violation and enforcement (Detail A in \cref{fig:survey_30bus_state}); Left: projected saddle-point flow; Right: linearized output projection;}\label{fig:cstr_viol}
        \end{figure}

    \end{example}

    \begin{example}\label{ex:lop_survey}
        We reconsider the numerical test case from \cref{ex:pw_ex_pd}, but apply the feedback control strategy \cref{eq:lop_up_survey} instead of \eqref{eq:pd_control}. Again, we choose a sampling period (corresponding to the step size $\alpha$) of 1 minute. As before we approximate $\nabla h(u)$ based on input and output measurements and the power flow model, and we choose $Q(u): = {H}(u)^T {H}(u)$ as the metric.

        The right-hand side panels in \cref{fig:survey_30bus_cost,fig:survey_30bus_state} show how this linearized output projection method achieves similar performance to the discretized projected saddle-point flow controller \cref{eq:pd_control}.

        In fact, the constraint satisfaction is even superior, as output constraints such as line flow limits are directly enforced. Whereas for the saddle-point transient constraint violation is generally unavoidable to achieve convergence to non-zero dual variables at an equilibrium, the linearized output projection method predicts (and avoids) constraint violations based on first-order model information. This effect is illustrated in \cref{fig:cstr_viol} where the enforcement of a line flow constraint from the previous simulations is compared. It is readily noticeable that the augmented projected saddle-point flow exhibits transient constraint violation and only asymptotically satisfies the constraint. This observation corresponds to the mechanism illustrated in \cref{fig:cstr_asadd}. In contrast, the linearized output projection exhibits almost no transient constraint violation and this close in behavior to the projected gradient flow it approximates (see also \cref{fig:cstr_pgrad}), at the cost of slightly higher computational complexity (a quadratic program needs to be solved at each iteration).

        {\tb Finally, we highlight that the feedback controller \eqref{eq:lop_up_survey} has gone beyond numerical tests and reached real-world industrial deployment. \cite{ortmann2023deployment} apply the controller \eqref{eq:lop_up_survey} to optimize reactive power flows of a real distribution grid subject to voltage constraints. The deployed controller features compatibility with the existing grid infrastructure and 24/7 uninterruptible operation.}
    \end{example}

\section{Conclusions and Outlook}\label{sec:conc}

In this article, we have surveyed various approaches for solving optimization problems in closed loop with a physical system. As a particular focus, we have interpreted classical numerical optimization algorithms as dynamical systems, and we have studied the conceptually simple idea of interconnecting these algorithms with asymptotically stable plants with well-defined input-output behavior.

Solving optimization problems online and with feedback, rather than offline and in open loop, leads to greater robustness against uncertain problem data and reduces required model information and computational effort because the physical system inherently enforces certain constraints.
The design of such controllers is aided by the recent rediscovery of and interest in the dynamical systems perspective of optimization algorithms. However, in contrast to analyzing existing algorithms with tools from dynamical systems theory, implementing these algorithms as feedback controllers raises novel and unique challenges, such as their closed-loop stability and robustness. Feedback-based optimization calls for smart control designs that meaningfully exploit the physical properties of a system, such as its steady-state response and saturation effects. We have particularly delved into the available guarantees of robustness and constraint satisfaction.

\subsection{Ongoing and Future Research Avenues}

Even though feedback-based optimization has been an active research area for some time, many questions remain unanswered. To complete this article, we propose two worthwhile avenues for future inquiry.

\subsubsection*{General Equilibrium Seeking in Closed Loop}

In this article, we have focused on solving optimization problems and thus steering a physical plant to a state that solves the corresponding KKT conditions. This setup is motivated by many applications that come with steady-state specifications regarding global plant efficiency, running cost, yield, etc.
Instead of optimality conditions, it is also possible to design controllers that seek states that satisfy more general equilibrium conditions, such as solving variational inequalities or generalized equalities. It thus becomes possible, for example, to design feedback controllers in a non-cooperative setting that allow players to converge to a Nash equilibrium \citep{depersisFeedbackControlAlgorithm2019}. The closed-loop structure and the use of online measurements also help to achieve robust tracking of time-varying Nash equilibria due to changing disturbances and problem parameters \citep{belgioioso2021sampled,belgioioso2022online}.

\subsubsection*{Time-Varying Non-Convex Optimization}

As discussed in \cref{sec:tv_opt}, the study of time-varying optimization problems has so far focused on convex problems (often with a strongly convex objective). This offers the convenience of a unique global optimizer that forms a well-defined trajectory over time and can be tracked by appropriate controllers and algorithms.

Solving non-convex problems, such as the optimal reserve dispatch problem from \cref{sec:redispatch}, in closed loop poses a much more challenging problem. Stationary non-convex problems are generally hard to solve due to the existence of multiple local minimizers. Time-varying non-convex problems pose additional problems since the possible appearance and disappearance of minimizers make it generally impossible to track a single global minimizer continuously. Nevertheless, it is a practical necessity to understand how feedback-based optimization performs in non-convex settings and to develop controllers that can cope with such pathologies. Recent papers exploring this direction include \cite{dingEscapingSpuriousLocal2021,linOnlineOptimizationPredictions2020a}.

\section*{Acknowledgements}
The research leading to these results was supported by ETH Z\"urich funds, the Swiss Federal Office of Energy grant \#SI/501708 UNICORN, the Max Planck ETH Center for Learning Systems, and the Swiss National Science Foundation through the NCCR Automation under Grant 180545.

\appendix
\renewcommand{\thesection}{\Alph{section}}

\section{{\tb Historical/Pioneering Applications}}
\label{app:historical}

\subsection{Network Congestion Control}\label{ex:netw_cong}

Arguably one of the largest man-made distributed feedback systems is formed by the congestion management mechanisms at the heart of the Internet~\citep{lowInternetCongestionControl2002}. These protocols manage the allocation of link capacities to individual connections in highly dynamic environments, for heterogeneous agents, and subject to real-world imperfections such as delays, among others~\citep{tangEquilibriumHeterogeneousCongestion2007,wenUnifyingPassivityFramework2004,paganiniCongestionControlHigh2005,vinnicombeStabilityNetworksOperating2002}.
Deterministic, continuous-time flow models have been particularly successful in explaining these protocols, by providing an optimization-based perspective in terms of network utility maximization~\citep{kellyRateControlCommunication1998,lowAnalyticalMethodsNetwork2017,lowOptimizationFlowControl1999,shakkottaiNetworkOptimizationControl2007}, and enabling improved designs~\citep{weiFASTTCPMotivation2006,weiDistributedNewtonMethod2013,zarghamAcceleratedDualDescent2013}.

Given a communication network, let a set of $N$ sources share $M$ links. The set of links used by each to communicate to its destination is collected in the routing matrix $R \in \bbR^{N \times M}$ in the sense that $R_{ij} = 1$ if link $j$ is used by source $i$ and $R_{ij} = 0$ otherwise.
Each link $j$ in the network has an associated congestion measure $\mu_j$ (also referred to as \emph{price}) which may describe queuing delays or packet loss probabilities. Each link $j$ has a finite capacity $c_j$. Each source $i$ has a controllable \emph{source rate} $x_i$, e.g., in the form of its \emph{transmission window size}. The different source rates define a vector of aggregate flows for each line given by $y = R x$ where each row of $R$ corresponds to a link. Conversely, sources are assumed to have access to their respective \emph{aggregate price} $p = R^T \mu$.
In practice, each source can estimate the total price over its own path by measuring delays, estimating the packet loss probability, or detecting the presence of specific ``congestion markers'' on its packets.

The simplest model to describe the \emph{link dynamics}
is given by the \emph{projected gradient} flow
\begin{equation}
    \dot \mu_j = \tproj{\bbR^M_{\geq 0}}{}{y_j - c_j}(\mu_j)
    \label{eq:linkmodel}
\end{equation}
where
\begin{align*}
     & \tproj{\bbR^M_{\geq 0}}{}{y_j - c_j}(\mu_j) = \begin{cases}
        y_j - c_j & \mu_j = 0 \text{ and }y_j > c_j, \\
        0         & \text{otherwise}.
    \end{cases}
\end{align*}
In other words, the price increases when the link is overloaded (and packets are lost) and decreases whenever it is not overloaded, but the price never drops below zero.

\emph{Source controllers} can often be modeled as $\dot x_i =  -\nabla \Phi_i(x_i) - p_i$, where $\Phi_i(x_i)$ is a specific type of cost function, reverse-engineered and depending on the particular protocol. Namely, sources adapt their rate to minimize a given cost function (equivalently, maximize their utility) corrected by the aggregate price.

\begin{figure}[!t]
    \centering
    \begin{tikzpicture}
        \matrix[ampersand replacement=\&, row sep=0.3cm, column sep=.55cm] {
            \& \node[branch](br1) {};
            \& \node[block](r_in) {$R$};    \\
            \node[block](grad) {$\nabla \Phi(x)$};
            \&\node[block](primal) {$\int$};
            \& \& \node[block] (sum) {$\Pi[y - c](\mu)$}; \\
            \& \node[smallsum](sum2){};
            \& \& \node[block](dual) {$\int$}; \\
            \&\& \node[block](r_out) {$R^T$}; \&
            \node[branch] (br3) {}; \&
            \node[none](edgeout) {};       \\
        };

        \draw[line] (primal.north)--(br1.south);
        \draw[connector] (br1.east)--(r_in.west) node[at start, above]{$x$};
        \draw[connector] (r_in.east)-|(sum.north) node[midway, above]{$y$};
        \draw[connector] (sum.south)--(dual.north);
        \draw[connector] (dual.south)|-(r_out.east) node[midway, below]{$\mu$};
        \draw[connector] (r_out.west)-|(sum2.south) node[near end, left]{$-$} node[midway, below] {$p$};
        \draw[connector] (sum2.north)--(primal.south);
        \draw[connector] (br1.center)-|(grad.north) ;
        \draw[connector] (grad.south)|-(sum2.west) node[near end, above]{$-$};

        \draw[line] (dual.south)|-(edgeout.center);
        \draw[connector] (edgeout.center)|-(sum.east) ;
    \end{tikzpicture}
    \caption{Simple network utility maximization}\label{fig:nw}
\end{figure}

Hence, congestion control mechanisms interconnect source controllers and link dynamics into a feedback loop as illustrated in \cref{fig:nw}, where $\Phi(x) := \sum_i \Phi_i(x_i)$. Importantly, each controller is fully distributed and requires only locally available information. Written more compactly, the closed-loop dynamics are defined as
\begin{align*}%
    \dot x =  - \nabla \Phi(x) - R^T \mu  \qquad\quad \dot \mu & =   \tproj{\bbR^M_{\geq 0}}{}{R x - c}(\mu) \, .
\end{align*}
This system defines a \emph{projected saddle-point flow}, as in \cref{ex:sadd_flow}, whose trajectories converge to a solution of
\begin{align}\label{eq:netw_util_prob}
    \minimize \quad \Phi(x) \qquad \subjto \quad R x \leq c \, ,
\end{align}
whenever $\Phi$ is strictly convex.

From a control and optimization perspective, source controllers and link dynamics form a closed feedback loop that implicitly tracks the solution of the underlying utility maximization problem~\citep{wangControlPerspectiveCentralized2011}.

From the viewpoint of closed-loop optimization advocated in this article, one can argue that the links implement the dual dynamics that have to be complemented by controllers that form the primal dynamics. Together they optimize \cref{eq:netw_util_prob}.

\subsection{Optimal Frequency Control}\label{ex:opt_freq}
Recent years have seen renewed interest in the control of power systems because of the new challenges associated with an increasingly dominating number of power electronic devices connected to the electricity grid and the growth of highly intermittent infeed from new renewable energy sources.

Frequency control in AC grids is one control task that has been revisited in this context. In this example, we illustrate how recent work has framed this challenge as a closed-loop optimization problem in order to improve upon classical frequency control schemes.

The reader is referred to \cite{molzahnSurveyDistributedOptimization2017} for a recent survey on this topic and to \cite{zhaoUnifiedFrameworkFrequency2016,simpson-porcoStabilityDistributedAveragingProportionalIntegral2021,chenDistributedAutomaticLoad2020} and references therein for latest results. The ensuing model is based on \cite{liConnectingAutomaticGeneration2016}.

In contrast to the optimal reserve dispatch problem outlined in \cref{sec:redispatch}, we assume a \emph{lossless} AC power transmission network with a set $\calN$ of buses and a set $\calM$  of transmission lines. We consider the power system around an operating point and all quantities are to be interpreted as deviations from their nominal values.

The frequency deviation $\omega_j$ at every bus $j$ is governed by the \emph{swing equation}
\begin{align}\label{eq:f_swing}
    \dot{\omega}_j = \tfrac{1}{M_j}
    \bigg(
    p_j^M - D_j \omega_j - p^{\L}_j -
    \sum\nolimits_{k: j \rightarrow k} p_{jk}
    \bigg)
\end{align}
where $M_j$ is the inertia of the generator and $D_j$ is the damping constant at bus $j$. The mechanical power output of the generator and the power consumed at bus $j$ are denoted by $p^{\M}_j$ and $p^{\L}_j$, respectively. The power flowing out of bus $j$ towards an adjacent bus $k$ is written as $p_{jk}$. The notation $j \rightarrow k$ indicates a transmission line between buses $j$ and $k$, i.e., $(j,k) \in \calM$.

For every $j \rightarrow k$, \emph{line flow dynamics} are linearized and modeled as
\begin{align}\label{eq:f_line}
    \dot{p}_{jk} = \rmb_{jk}( \omega_j - \omega_k)
\end{align}
where $\rmb_{jk}$ is the line susceptance. Note that transmission lines are undirected, and, throughout, the constraint $p_{jk} = - p_{kj}$ holds for every line implicitly.

At steady state and nominal frequency $\omega_j = 0$, the swing equation \eqref{eq:f_swing} expresses the power balance at each node, and \eqref{eq:f_line} yields the so-called \emph{DC-flow} approximation of the AC power flow equations  (assuming nominal voltages, lossless lines, and small angle differences).

The mechanical power output $p^{\M}_j$ is described by a simplified \emph{governor-turbine control} model of the form
\begin{align}\label{eq:f_gov}
    \dot{p}_j^M = - \tfrac{1}{T_j} \left( p^{\M}_j - p^{\C}_j + \tfrac{1}{R_j} \omega_j \right)
\end{align}
where $p^{\C}_j$ is a power adjustment signal and $T_j$ and $R_j$ are constants. In particular, $R_j$ is understood as a generator's participation factor in primary frequency control.

The power adjustment itself is evaluated using the so-called \emph{area control error} (ACE) which combines frequency and line flow deviations. More specifically, we have
\begin{align}\label{eq:f_ace}
    \dot{p}_j^C = - K_j \bigg( B_j \omega_j + \sum\nolimits_{k: j \rightarrow k } p_{jk} \bigg) \, .
\end{align}
with $B_j$ and $K_j$ being constant control parameters.

The role of this type of \emph{automatic generation control} (or \emph{secondary frequency control}) is to react to changes in $p^{\L}_j$, driving frequency deviation to zero, and to make sure that deviations in $p^{\L}_j$ are compensated by the local generator through adjustments of $p^{\M}_j$. Consequently, at steady state, the deviations $P_{ij}$ of the line flows from their nominal values are zero.\footnotemark

\footnotetext{In this example, the area control error is evaluated and regulated on a per-bus basis. In more general settings, multiple buses belong to the same area, and therefore secondary frequency control regulates the aggregate power balance in each area and the power flows between areas.}

Although not immediately obvious, the system \crefrange{eq:f_swing}{eq:f_ace} realizes a (partial) saddle-point flow, as introduced in \cref{sec:saddle}. To see this, we consider the optimization problem
\begin{align}\label{eq:f_opt_prob}
    \begin{split}
        \underset{\tb p^\M, \omega, \{p_{jk}\}_{jk}}{\minimize} \quad & \sum\nolimits_{j\in \mathcal{N}} \left(\tfrac{\beta_j}{2} \left(p^{\M}_j\right)^2 + \tfrac{D_j}{2} \omega_j^2 \right) \\
        \underset{\forall j \in \calN}{\subjto} \quad   & 0 = p^{\M}_j - p^{\L}_j  \\
        & 0 = p^{\M}_j - D_j \omega_j - p^{\L}_j - \sum\nolimits_{k: j \rightarrow k} p_{jk} \, .
    \end{split}
\end{align}
where $\beta_j$ is a positive parameter to be determined.
The solution of this problem guarantees power balance and exact frequency regulation.
We can write out the corresponding (unprojected) primal-dual saddle-point flow from \cref{ex:pd_eqcstr_basic} as
\begin{align}\label{eq:f_pd_flow}
    \begin{split}
        \text{\rotatebox[origin=c]{90}{\tiny primal}}
        & \begin{cases} \dot{\omega}_j             & = - \epsilon_{\omega_j} \left(D_j \omega_j - D_j \lambda_j \right) \\
            \dot{p}^M_j    \mkern-12mu & = - \epsilon_{P^M_i} (\beta_j p^{\M}_j - \mu_j - \lambda_j)        \\
            \dot{p}_{jk}   \mkern-12mu & = - \epsilon_{p_{jk}} (\lambda_j - \lambda_k)
        \end{cases} \\
        \text{\rotatebox[origin=c]{90}{\tiny dual}}
        & \begin{cases}
            \mathrlap{\dot{\mu}_j}\hphantom{\dot{p}^M_j   \mkern-12mu} & = \epsilon_{\mu_j} (p^{\M}_j - p^{\L}_j) \\
            \dot{\lambda}_j                                            & = \epsilon_{\lambda_j} \bigg(
            p^{\M}_j - D_j \omega_j - p^{\L}_j -
            \sum\nolimits_{k: j \rightarrow k} p_{jk} \bigg)
        \end{cases}
    \end{split}
\end{align}
for all $j,k \in \calN$ and $j \rightarrow k$, and
where we use separate positive gains $\epsilon_{(\cdot)}$ for every variable.

By taking the limit $\epsilon_{\omega_j} \rightarrow \infty$, we can replace the corresponding ODE with the algebraic expression $0 = D_j (\omega_j - \lambda_j)$
and thus substitute $\omega_j = \lambda_j$. This results in a so-called \emph{partial} saddle-point flow. In this special case, the $\omega_j$ can hence be interpreted as a primal or dual variable.

We further apply a transformation to $\mu_j$ by defining $p^{\C}_j := K_j \left( M_j \omega_j - \tfrac{1}{\epsilon_{\mu_j}} \mu_j\right)$.
Consequently, we have
\begin{align*}
    \dot{p}^C_j = K_j \left( M_j \dot{\omega}_j - \tfrac{1}{\epsilon_{\mu_j}} \dot{\mu}_j \right) = K_j \big( M_j \dot{\omega}_j - (p^{\M}_j - p^{\L}_j) \big) \, .
\end{align*}

Combining these insights, we can rewrite the partial saddle-point flow derived from \eqref{eq:f_pd_flow} as
\begin{align}\label{eq:f_partial_pd}
    \begin{split}
        \text{\rotatebox[origin=c]{90}{\tiny primal}}
        & \begin{cases}
            \dot{p}^M_j  \mkern-12mu & = - \epsilon_{P^M_i} (\beta_j p^{\M}_j - \mu_j - \omega_j) \\
            \dot{p}_{ij} \mkern-12mu & = - \epsilon_{P_{ij}} \left(\omega_k - \omega_j \right)    \\
        \end{cases} \\
        \text{
            \rotatebox[origin=c]{90}{\tiny primal}
            \rotatebox[origin=c]{90}{\tiny or dual}
        }
        & \begin{cases}
            \mathrlap{\dot{\omega}_j}\hphantom{\dot{p}^M_j   \mkern-12mu}
             & =  \epsilon_{\lambda_j} \bigg(
            p^{\M}_j - D_j \omega_j - p^{\L}_j -
            \sum\limits_{k: j \rightarrow k} p_{jk} \bigg)
        \end{cases} \\
        \text{\rotatebox[origin=c]{90}{\tiny dual}}
        & \begin{cases}
            \mathrlap{\dot{p}^C_j}\hphantom{\dot{p}^M_j   \mkern-12mu}
             & = - K_j \bigg( D_j \omega_j + \sum\limits_{k: j \rightarrow k } p_{jk} \bigg) \, .
        \end{cases}
    \end{split}
\end{align}
By carefully choosing the remaining gains $\epsilon_{(\cdot)}$ and $\beta_j$ (e.g., $\epsilon_{P_{ij}} := \rmb_{ij}$ or $\omega_{\lambda_j} = \tfrac{1}{M_j}$) as well as, under $B_j = D_j$,  \eqref{eq:f_partial_pd} is equal to \crefrange{eq:f_swing}{eq:f_ace}. The choice of $B_j = D_j$ is fragile, since it means that the controller \eqref{eq:f_ace} requires an accurate estimate of the damping $D_j$ at every bus. However, numerical experiments in \cite{liConnectingAutomaticGeneration2016} show that \crefrange{eq:f_swing}{eq:f_ace} is stable even if $B_j = D_j$ does not hold exactly. 

Starting from this observation that automatic generation control can be interpreted from a closed-loop optimization perspective, various works have proposed optimal frequency control schemes that allow for more general objective functions than in \cref{eq:f_opt_prob}, relax the constraint $p^{\M}_j = p^{\L}_j$, or incorporate additional constraints on generation and transmission capacity.

\section{{\tb Related Research Streams}}\label{sec:cl_opt_approaches}

The key idea studied in this article is the interconnection of optimization algorithms with physical systems. The architecture that we discussed (exemplified in \cref{ex:simp_grad}) is by no means the only meaningful way to combine optimization with feedback control. This section provides a pedagogical overview of various other well-established approaches in this direction. Our treatment is deliberately simplified to highlight the fundamental differences.
The reader is referred to respective survey papers for exhaustive reviews of the topics treated in this section.

We focus our attention on the following three methods.

\begin{enumerate}
    \item \emph{\Acl*{es}}, {\color{blue} which is rooted in} adaptive control, emphasizes being completely \emph{model-free} and probing the system with the help of a perturbation signal. However, the approach is best suited for systems with low-dimensional inputs and without complicated engineering constraints.

    \item Emerging from process engineering, the primary merit of the \emph{\acl*{ma}} method is to mitigate the effects of model bias when solving successive optimization problems. This method does not reduce the computational requirements compared to feedforward optimization. 

    \item \emph{\Acl*{rti}} schemes, rooted in model predictive control, aim to solve classical receding horizon problems with limited resources. The focus is the stabilization of a plant under state and input constraints. For this purpose, a full model of the plant dynamics is generally required, and the computational burden scales with the planning horizon.

\end{enumerate}

Many more approaches can be interpreted as optimization in conjunction with feedback control. Several topics such as \emph{iterative feedback tuning} \citep{hjalmarssonIterativefeedbacktuning2002} or \emph{iterative learning control} \citep{bristowSurveyIterativeLearning2006} are concerned with the optimal tuning of controllers by either requiring a sequence of experiments or repetitive operation. Although these methods use measurement data to optimize subsequent control actions, because of their ``episodic'' nature, they are related only remotely to the key idea in this article.

For similar reasons, we do not review recent work on \emph{reinforcement learning}, even though it has produced staggering results, especially in sequential decision-making for games and in robotics. Even though techniques such as \emph{policy gradient} admit a dynamical perspective, RL is generally (and traditionally) framed as optimal control over Markov decision processes \citep{bertsekasDynamicProgrammingOptimal2017,lewisReinforcementLearningApproximate2012}.

\subsection{Extremum Seeking}\label{sec:es}

Arguably, one of the oldest control methods to steer a plant to an extremum of a function rather than tracking a setpoint is \ac{es} (see \citealt{ariyurRealTimeOptimization2003,tanExtremumSeeking19222010} for historical accounts). Its popularity rose in the 1950's and 60's as part of adaptive control \citep{astromAdaptiveControl2008} and later regained momentum with the rigorous theoretical results derived in \cite{wittenmarkAdaptiveExtremalControl1995,krsticStabilityExtremumSeeking2000,tanNonlocalStabilityProperties2006}.

The main idea behind ES is to inject a \emph{dither signal} to locally explore the objective function and ``learn'' its gradient. This dither signal is generally sinusoidal, but other perturbations have been proposed \citep{teelSolvingSmoothNonsmooth2001}. Consequently, the objective can be optimized without recourse to any model information about the plant (and objective) and without any computation aside from the addition and multiplication of the dither signal. The following example illustrates this fact.

\begin{example}\label{ex:es1}
    \begin{figure}[!t]
        \centering
        \begin{tikzpicture}
            \matrix[ampersand replacement=\&, row sep=0.3cm, column sep=.55cm] {
                \&\node[smallsum](sum1){};
                \&\&\& \node[none](edgetop) {}; \\

                \node[gainup](gain1) {$2 a$};
                \&\node[block](gain2) {$\int$};
                \& \& \&
                \node[block] (plant) {
                    $\begin{aligned}
                            \dot \zeta & = {f}(\zeta, \tilde{u}) \\
                            \hat{y}    & = g(\zeta)
                        \end{aligned}$
                };        \\

                \node[none](source){$\sin(\omega t)$};
                \& \node[gainup] (int) {$\tfrac{\epsilon}{a}$}; \& \& \&
                \node[smallsum] (dist_sum){};        \\

                \&\node[smallsum](sum2){$\times$}; \&
                \node[block] (cost) {$-\Phi({y})$};
                \&\& \node[none](edgebot) {};   \\
            };

            \draw[connector] (source.north)|-(gain1.south);
            \draw[connector] (gain1.north)|-(sum1.west) node[at end, above] {$+$};
            \draw[connector] (sum1.east)-|(plant.north) node[near start, below] {$\tilde{u}$};
            \draw[connector] (plant.south)--(dist_sum.north);
            \draw[connector] ([xshift=.5cm]dist_sum.east)--(dist_sum.east) node[at start, right]{$d$};
            \draw[connector] (dist_sum.south)|-(cost.east) node[near end, above] {${y}$};
            \draw[connector] (cost.west)--(sum2.east);
            \draw[connector] (source.south)|-(sum2.west);
            \draw[connector] (sum2.north)--(int.south);
            \draw[connector] (int.north)--(gain2.south);
            \draw[connector] (gain2.north)--(sum1.south) node[near end, left] {$+$} node[midway,right] {${u}$};
        \end{tikzpicture}
        \caption{Simple extremum seeking to minimize $\Phi(h(u))$}\label{fig:es1}
    \end{figure}

    We consider the same problem as in \cref{ex:simp_grad}, i.e., the minimization of $\Phi(y)$ where $y$ is the output of a plant of the form \eqref{eq:simple_plant}. In addition, we assume that the plant is single-input-single-output (and thus all signals are scalar).
    Instead of the gradient scheme in \cref{fig:simp_grad}, we apply the \ac{es} setup in \cref{fig:es1} where a sinusoid perturbs the signal $u$, yielding $\tilde{u}$, which is then fed to the plant.

    As in Example \ref{ex:simp_grad}, we replace the fast plant dynamics by the steady-state map $y = h(u) + d$, and we define the reduced cost function $\tilde{\Phi}(u) := \Phi(h(u) + d)$.
    Consequently, the \emph{reduced} \ac{es} dynamics of the system in \cref{fig:es1} can be expressed in terms of $\dot{u}$ as
    \begin{align*}%
        \dot u = F(u,t) :=  - \tfrac{\epsilon}{a} \underbrace{\Phi\big(h \left( u + 2a \sin(\omega t)\right) + d \big)}_{\tilde{\Phi}(u + 2 a \sin(\omega t))} \sin(\omega t) \, .
    \end{align*}
    The \emph{averaged} dynamics are obtained by integrating $F(u,t)$ from 0 to $T = \tfrac{2\pi}{\omega}$. Namely, using the Taylor expansion of $\tilde{\Phi}$ around $u$, the average control signal is
    \begin{align*}
         & \tfrac{1}{T} \int_0^{T} F(u, t) dt = - \tfrac{\epsilon}{a T} \int_0^{T} \tilde{\Phi} \big( u + 2a \sin(\omega t) \big) \sin(\omega t) dt \\
         & \,\approx -  \tfrac{\epsilon}{aT} \int_0^{T} \big(
        \tilde{\Phi}\left(u\right)
        + 2a \sin(\omega t)  \nabla \tilde{\Phi}(u) \big) \sin(\omega t) dt                                                                         \\
         & \,=  -  \tfrac{\epsilon}{a}  \tfrac{\omega }{2\pi} 2a \nabla \tilde{\Phi}(u)  \int_0^{T}
        \sin(\omega t)^2   dt                \, = - \epsilon \nabla \tilde{\Phi}(u) \, ,
    \end{align*}
    where we have neglected higher-order terms in $\epsilon$. Thus, the \ac{es} scheme approximates the gradient flow \eqref{eq:simple_grad_flow} from \cref{ex:simp_grad}. However, \ac{es} merely requires measurements of $\Phi(y)$ and neither an estimate of $\nabla h$ nor of $\nabla \Phi$.

    In contrast to the design in \cref{ex:simp_grad}, \ac{es} systems generally evolve on three (rather than two) different timescales: the plant dynamics (which have been ignored for this example), the frequency range of the probing signal, and the slow averaged optimization dynamics.
\end{example}

Classically, averaging theory and singular perturbation analysis (for dynamic plants) are used to render the insights from \cref{ex:es1} rigorous \citep{krsticStabilityExtremumSeeking2000,tanNonlocalStabilityProperties2006,guayAdaptiveExtremumSeeking2003}.
More recently, \ac{es} schemes have also been studied with the help of Lie bracket approximations which offer an alternative perspective~\citep{durrLieBracketApproximation2013,durrExtremumSeekingDynamic2017,grushkovskayaClassGeneratingVector2018}.

While original work only considered finding extrema (i.e. minima or maxima) of unconstrained problems, constraints have been incorporated by submanifold constraints \citep{durrExtremumSeekingSubmanifolds2014}, barrier function \citep{dehaanExtremumseekingControlStateconstrained2005,guay2015constrained}, and saddle-point formulations \citep{durrSaddlePointSeeking2013}. Beyond asymptotic constraint satisfaction, \cite{chen2023continuous,chen2021safe} introduce projection maps into ES to robustly handle constraints while achieving model-free optimization.
\ac{es} has also showcased its power in Nash-equilibrium seeking \citep{frihaufNashEquilibriumSeeking2012,stankovicDistributedSeekingNash2012}.
Further, \ac{es} has been studied for stochastic \citep{coitoStochasticExtremumSeeking2005,stankovicExtremumSeekingStochastic2010} and discrete-time setups \citep{stankovicDiscreteTimeExtremum2009,feilingDerivativeFreeOptimizationAlgorithms2018,frihaufFinitehorizonLQControl2013}, and hybrid extensions have been proposed \citep{povedaFrameworkClassHybrid2017,povedaDistributedExtremumSeeking2017,poveda2021robust,bianchinTimeVaryingOptimizationLTI2021}.

\Ac{es} has been applied in the automotive sector \citep{killingsworthHCCIEngineCombustionTiming2009}, process engineering \citep{guayAdaptiveExtremumSeeking2004}, formation flight and obstacle avoidance \citep{binettiFormationFlightOptimization2003,montenbruckExtremumSeekingObstacle2014} and others.
Further applications concern problems in renewable energy such as maximum power point tracking in photovoltaic \citep{ghaffariMultivariableMaximumPower2015} or wind energy systems \cite{ghaffariPowerOptimizationControl2014,krsticExtremumSeekingWind2014}, and Volt-VAR control in power systems \citep{arnoldModelFreeOptimalControl2016}.

Despite strong theoretical guarantees and being model-free, \ac{es} has been confined to relatively low-dimensional, mostly unconstrained, systems. This is due to the fact that plants with multidimensional input require probing signals at different, carefully chosen, frequencies that do not interfere with each other.%

\subsection{Modifier Adaptation}

In the context of \emph{real-time optimization} in process engineering, the notion of \emph{measurement-based optimization} \citep{francoisMeasurementBasedRealTimeOptimization2013,chachuatAdaptationStrategiesRealtime2009} has been used to collect several approaches towards mitigating the effects of model bias in repetitive optimization applications. In the following, we showcase \ac{ma} methods by \cite{marchettiRealtimeOptimizationEstimation2009,marchettiModifierAdaptationMethodologyRealTime2009,gaoIterativeSetpointOptimization2005,francoisUseMeasurementsEnforcing2005}.

Given a model of a physical system, assume that we can solve an optimal steady-state problem like \eqref{eq:basic_prob} numerically. When implementing this solution by setting the appropriate inputs of the stable plant to their pre-computed optimal setpoints, the mismatch between the model estimate (used for computing an optimal state) and the actual plant will invariably lead to a system state that is suboptimal, and possibly violating constraints.

If the optimization of the optimal plant state is performed repeatedly, and at each step the solution is implemented on the physical system, \ac{ma} provides a method to steadily reduce the discrepancy between model-based solution and physical plant by \emph{modifying} the optimal steady-state problem at every iteration by incorporating plant measurements from the previous iteration.
\ac{ma} does not directly ``learn'' or identify a better model of the plant. Instead, \ac{ma} corrects the optimization problem by adding adaption terms to the cost and constraint functions.
The following simple example demonstrates one possible adaptation mechanism.

\begin{example}\label{ex:ma}
    Consider the same setup as in \cref{ex:simp_grad}. Namely, we wish to minimize the function $\tilde{\Phi}(u) := \Phi(h(u)+d)$ where $y= h(u)+d$ is the steady-state input-output map of a plant with fast-decaying dynamics. %
    However, only an approximate model $\tilde{h}$ of $h$ and an estimate $\tilde{d}$ of $d$ is available.

    Therefore, instead of minimizing $\tilde{\Phi}$, we repeatedly solve
    \begin{align}\label{eq:ma}
        \underset{\tb u}{\minimize} \quad \tilde\Phi(\tilde{h}(u) + \tilde{d}) + \lambda_k u \, ,
    \end{align}
    where $\lambda_k$ is a \emph{modifier} at iteration $k$ that is adapted at every iteration based on the outcome of the previous iteration $u^\star_k$. %
    In particular, $\lambda$ is updated according to
    \begin{align*}
        \lambda_{k+1} & = \nabla \tilde{\Phi}(u^\star_{k}) - \nabla \tilde{\Phi}'(u^\star_k) \, ,
    \end{align*}
    where $\nabla \tilde{\Phi}(u^\star_{k})$ needs to be estimated, and $\nabla \tilde{\Phi}'(u^\star_k) := \nabla ( \Phi(\tilde{h}(u^\star_k)+\tilde{d})$ is model-based. The particular structure of $\tilde{\Phi}$, however, lets us write
    \begin{align}
        \nabla \tilde{\Phi}(u^\star_{k}) = \nabla {\Phi}(h(u^\star_{k})+d) \nabla h(u^\star_{k})  \, ,
    \end{align}
    where $h(u^\star_{k}) + d$ is the measured output of the plant. Thus, essentially, only $\nabla h(u^\star_{k})$ needs to be estimated.
    If the scheme converges to some $u^\star$, we can easily verify that $u^\star$ is a critical point of $\tilde{\Phi}(u)$.
\end{example}

\cref{ex:ma} is simplified to the point that a comparison with \cref{ex:simp_grad} is easily possible. However, \ac{ma} methods are easily applied to constrained problems, where modifiers on constraints are introduced analogously \citep{costelloModifierAdaptationConstrained2014,faulwasserEconomicNonlinearModel2018,
    gregoryUseTransientMeasurements2014}.

Clearly, the tricky part about \ac{ma} is the estimation of the (true) plant sensitivities $\nabla \tilde{\Phi}(u^\star_{k})$. This can be achieved with finite differences \citep{mansourComparisonMethodsEstimating2003}. Moreover, \ac{ma} does not reduce the computation burden nor does it aim to reduce the amount of model information required.

\subsection{Model Predictive Control with Incomplete Optimization and Real-Time Iterations}\label{sec:rti}

Historically, \ac{mpc} is an approach to control and stabilize a plant that is subject to input and state constraints. This is achieved by numerically solving an optimal control problem with a finite receding horizon at every sampling time (or, every few sampling instants), but implementing only the first (respectively, few) input(s) of the computed optimal policy before solving the next problem with shifted horizon and based on an updated state measurement.

The high computational requirements have long restricted the application of \ac{mpc} to relatively slow and low-dimensional plants in process engineering. For standard (linear) \ac{mpc}, this issue has led to \emph{explicit} \ac{mpc} \citep{bemporadExplicitLinearQuadratic2002,alessioSurveyExplicitModel2009} which exploits multi-parametric optimization \citep{tondelAlgorithmMultiparametricQuadratic2003} to solve the receding horizon problem ahead of time and implement the controller as a simple lookup table.

More interesting from our perspective are \acp{rti} for \emph{nonlinear \ac{mpc}} \citep{bockDirectMultipleShooting2000,diehlRealtimeOptimizationNonlinear2002}. These methods have emerged as an approximation of \emph{multiple shooting methods} \citep{diehlFastDirectMultiple2006} and have been proposed for various applications in process engineering \citep{diehlEfficientAlgorithmNonlinear2003}, robotics \citep{diehlFastDirectMultiple2006}, and for airborne kites \citep{diehlRealTimeIterationScheme2005}.
The main idea of \acp{rti} is to solve the optimal control problem only approximately at every iteration by performing only a single iteration of the underlying optimization algorithm (which is usually a \ac{sqp} scheme; \citealt{nocedalNumericalOptimization2006,quirynenInexactNewtonTypeOptimization2018,zavalaNonlinearProgrammingStrategies2009}). The first input of the approximate control policy is implemented and the optimization problem for the next sampling period is warm-started at a shifted version of the previous (approximate) solution. In other words, \acp{rti} are a special case of \ac{mpc} methods with incomplete optimization \citep{liao-mcphersonTimedistributedOptimizationRealtime2020,liao2022analysis,figura2020instant,graichenStabilityIncrementalImprovement2010,gruneAnalysisUnconstrainedNMPC2010}. \Cref{ex:rti} below illustrates this procedure.

The underlying idea of \acp{rti} is that the approximation error committed by performing only a single optimization iteration is offset by savings in computation time. In particular, because the receding horizon problem is solved more often, it changes less between samples. This feature allows one to prove stability and convergence of \ac{rti} schemes \citep{diehlStabilizingRealTimeImplementation2007,diehlRealTimeIterationScheme2005,diehlNominalStabilityRealtime2005,zanelliLyapunovFunctionCombined2020}.

However, although they interleave optimization iterations with physical dynamics, \acp{rti} have been developed for stabilization and require an exogenous setpoint as well as a dynamic model of the plant. This is particularly reflected in the assumptions on the state cost function, which are, roughly speaking, required to be quadratic functions centered at the origin (see also \cref{rem:empc} further below on economic \ac{mpc}).

\begin{example}\label{ex:rti}
    We consider a discrete-time plant $\zeta^+ = f(\zeta, u)$ for which the origin is a steady state, i.e., $0 = f(0, 0)$. For simplicity, we do not model any constraints, although \acp{rti} can incorporate them naturally.%

    Consider the receding horizon problem at time~$l$
    \begin{align}\label{eq:realit_opc}
        \begin{split}
            \underset{\{r\}^{K-1}_{1}, \{s\}_{1}^K}{\minimize} \quad &
            \sum_{k=1}^{K-1} \begin{smallbmatrix} s_k \\ r_k\end{smallbmatrix}^T Q \begin{smallbmatrix} s_k \\ r_k\end{smallbmatrix} + s_K^T R s_K \\
            \subjto \quad & s_1 = \overline{\zeta}_l  \\
            &  s_{k+1} = f(s_k, r_k) \quad \forall k \in \{1, \ldots, K-1 \}%
        \end{split}
    \end{align}
    where $K$ denotes the horizon length, $Q,R$ are positive definite stage and terminal cost matrices, and $\overline{\zeta}_l$ denotes the measured plant state at time $l$.
    Let $(\hat{r}^{l}, \hat{s}^{l})$ denote the solution of \eqref{eq:realit_opc} for the sampling instant $l$. Then, the feedback law at $l$ is given by $u[l] = \hat{r}^{l}_1$. In other words, upon solving \eqref{eq:realit_opc}, the first control of the optimal policy is implemented at $l$. This is the key mechanism behind standard \ac{mpc}.

    \Ac{rti} schemes approximate the solution of \eqref{eq:realit_opc} by performing only a single iteration of an \ac{sqp} method. Namely, let $z = (s, r, \lambda)$ and consider the Lagrangian of \eqref{eq:realit_opc} at time $l$ defined as
    \begin{multline*}
        L^l(z) :=  \sum_{k=1}^{K-1} \begin{smallbmatrix} s_k \\ r_k\end{smallbmatrix}^T Q \begin{smallbmatrix} s_k \\ r_k\end{smallbmatrix} + s_K^T R s_K  + \lambda_1 (s_1 - \overline{\zeta}_l) \\+ \sum_{k=1}^{K-1} \lambda_k \left( s_{k+1} - f(s_k, r_k) \right) \, .
    \end{multline*}
    An \ac{sqp} iteration then takes the form
    \begin{align}\label{eq:rti}
        z^+ = z + \Delta z
    \end{align}
    where $\Delta z$ solves the first-order condition
    \begin{align*}
        \nabla L^l(z) + \nabla^2_{zz} L^l(z) \Delta z = 0 \, .
    \end{align*}
    where, in practice, the inverse of the Hessian $\nabla^2_{zz} L^l$ is often approximated.%

    Crucially, the receding horizon problem \eqref{eq:realit_opc} at the next sample $l+1$ is warm-started with the shifted approximate of the previous sample. This procedure leads, under additional assumptions, to local stability of the scheme.
\end{example}

\Cref{ex:rti} elucidates several differences with respect to the optimal steady-state control problem \cref{eq:basic_prob} we wish to solve: First, the purpose of \acp{rti} is primarily to drive a plant to the steady state at the origin, not seeking out an operating point with minimal cost (see \cref{rem:empc} below). Further, a full dynamic model $f$ of the plant dynamics is required, and finally, the computational burden of solving the \ac{sqp} iteration scales not only with the system dimension but also with the prediction horizon. Conversely, compared to \Cref{ex:simp_grad}, while stabilizing the system, RTI also seeks to minimize a quadratic running stage cost.

\begin{remark}\label{rem:empc} Traditionally, linear and nonlinear \ac{mpc} have been considered with the goal of stabilizing a plant and tracking a precomputed setpoint or trajectory. The modern variation of \emph{economic \ac{mpc}} \citep{faulwasserEconomicNonlinearModel2018, ellisTutorialReviewEconomic2014, rawlingsFundamentalsEconomicModel2012} studies the effects of incorporating an economic objective directly and thus not requiring an exogenous setpoint. This idea is very much in line with the topic of this article.

    Economic \ac{mpc}, however, still requires the solution of an optimal control problem at every iteration. Thus, it is computationally very expensive. Moreover, the optimal solution is not a-priori guaranteed to be a steady state of the plant. This makes the stability analysis of economic \ac{mpc} more involved and, in particular, still requires a full model of the plant dynamics.
\end{remark}


\end{document}